%% file: main.tex
\documentclass[11pt]{article}
\usepackage[numbers,longnamesfirst]{natbib}
\usepackage{amsmath,amsfonts,amssymb,amsthm,amscd,mathtools}
\usepackage{algorithm,algorithmic}
\usepackage[justification=centering]{caption}
\usepackage[justification=centering]{subcaption}
\usepackage{xstring}
\usepackage{scalerel}
\usepackage{setspace}
\usepackage{textcase}
\usepackage[dvipsnames]{xcolor}
\usepackage{import}
\usepackage{verbatim}
\usepackage{url}
\usepackage[pdftex, bookmarks=true, colorlinks=true, linkcolor=black, citecolor=black, urlcolor=black]{hyperref} 
\usepackage[capitalise]{cleveref}
\usepackage{graphicx,subcaption}
\usepackage{helvet}
\usepackage{bbm,enumerate,array}
\usepackage{graphics, parskip}
\usepackage{listings}
\usepackage{textcomp}
\usepackage{multirow,tabularx,ltxtable,longtable}
\usepackage[bottom]{footmisc}

\usepackage{tikz,xparse}
\usetikzlibrary{matrix,backgrounds}
\pgfdeclarelayer{myback}
\pgfsetlayers{myback,background,main}
\tikzset{mycolor/.style = {line width=1bp,color=#1}}%
\tikzset{myfillcolor/.style = {fill=#1}}%
\NewDocumentCommand{\highlight}{O{blue!40} m m}{%
  \draw[mycolor=#1] (#2.north west)rectangle (#3.south east);
}

\NewDocumentCommand{\fhighlight}{O{blue!40} m m}{%
  \fill[myfillcolor=#1] (#2.north west)rectangle (#3.south east);
}

\usepackage[margin=1in]{geometry}
\usepackage{fancyhdr}
\graphicspath{{figs/}{more-figs/}} 

\input{my_defs}

\title{Admissible Measurements and Robust Algorithms for Ptychography}

\author{Brian Preskitt \thanks{Pure Storage, Inc.,  e-mail: bpreskitt@purestorage.com} \and Rayan Saab \thanks{Department of Mathematics, University of California San Diego, e-mail: rsaab@ucsd.edu.}}
\date{}

\begin{document}
\maketitle

\begin{abstract}
We study an approach to solving the phase retrieval problem as it arises in a phase-less imaging modality known as ptychography. In ptychography, small overlapping sections of an unknown sample (or signal, say $x_0\in \C^d$) are illuminated one at a time, often with a physical mask between the sample and light source. The corresponding measurements are the noisy magnitudes of the Fourier transform coefficients resulting from the pointwise product of the mask and the sample. The goal is to recover the original signal from such measurements. 

The algorithmic framework we study herein relies on first inverting a linear system of equations to recover a fraction of the entries in $x_0 x_0^*$ and then using non-linear techniques to recover the magnitudes and phases of the entries of $x_0$.
Thus, this paper's contributions are three-fold. First, focusing on the linear part, it expands the theory studying which measurement schemes (i.e., masks, shifts of the sample) yield invertible linear systems, including an analysis of the conditioning of the resulting systems. Second, it analyzes a class of improved magnitude recovery algorithms and, third, it proposes and analyzes algorithms for phase recovery in the ptychographic setting where large shifts --- up to $50\%$ the size of the mask ---  are permitted. 
\end{abstract}

%
%


\import{sections/history/}{history_sec}


\subsection{Notation}
\import{sections/history/}{notation}

\section{Invertible Local Measurement Systems}
\label{ch:span_fam}
\import{sections/meas/}{meas_sec}


\section{Ptychographic Model}
\label{ch:ptychography}
\import{sections/ptychography/}{ptychography_sec}


\appendix
\import{sections/appendix/}{appendix}

\section*{Acknowledgments}
RS was supported in part by the NSF via DMS-1517204. The authors would like to thank Mark Iwen for many stimulating conversations on phase retrieval. 
\bibliographystyle{abbrv}


{\setstretch{1.0}
\bibliography{bibs/dissertation}}

\end{document}

%% file: my_defs.tex
\newtheorem{theorem}{Theorem}
\newtheorem*{thm*}{Theorem}
\newtheorem{corollary}{Corollary}
\newtheorem{lemma}{Lemma}
\newtheorem*{lemma*}{Lemma}
\newtheorem{proposition}{Proposition}
\newtheorem*{proposition*}{Proposition}

\newtheorem*{remark*}{Remark}

\theoremstyle{definition}

\let\mod\undefined
\let\Re\undefined
\let\Im\undefined
\let\vec\undefined

\DeclareMathOperator{\diag}{diag}

\DeclareMathOperator{\Col}{Col}
\DeclareMathOperator{\Row}{Row}
\DeclareMathOperator{\Span}{span}
\DeclareMathOperator{\Nul}{Nul}

\DeclareMathOperator{\supp}{supp}
\DeclareMathOperator{\Tr}{Tr}

\DeclareMathOperator{\sgn}{sgn}

\DeclareMathOperator{\circop}{circ}
\DeclareMathOperator{\mod}{\mathbin{mod}}

\DeclareMathOperator{\Re}{Re}
\DeclareMathOperator{\Im}{Im}

\DeclareMathOperator{\vec}{vec}
\DeclareMathOperator{\mat}{mat}
\DeclareMathOperator{\Skew}{Skew}
\DeclareMathOperator{\SNR}{SNR}

\DeclareMathOperator{\BlkMag}{BlockMag}
\DeclareMathOperator{\BlkVec}{BlockVec}

\DeclarePairedDelimiter{\floor}{\lfloor}{\rfloor}

\DeclarePairedDelimiter{\norm}{\lVert}{\rVert}
\DeclarePairedDelimiter{\abs}{\lvert}{\rvert}
\DeclarePairedDelimiter{\inner}{\langle}{\rangle}
\DeclarePairedDelimiter{\clopen}{[}{)}

\newcommand{\R}{\ensuremath{\mathbb{R}}}
\newcommand{\C}{\ensuremath{\mathbb{C}}}
\newcommand{\N}{\ensuremath{\mathbb{N}}}
\newcommand{\Z}{\ensuremath{\mathbb{Z}}}
\renewcommand{\H}{\ensuremath{\mathcal{H}}}
\newcommand{\sym}{\mathcal{S}}
\newcommand{\Sbb}{\ensuremath{\mathbb{S}}}

\newcommand{\Rn}{\ensuremath{\R^n}}
\newcommand{\Rm}{\ensuremath{\R^m}}
\newcommand{\Rd}{\ensuremath{\R^d}}
\newcommand{\Hd}{\ensuremath{\H^d}}

\newcommand{\Cn}{\ensuremath{\C^n}}
\newcommand{\Cm}{\ensuremath{\C^m}}
\newcommand{\Cd}{\ensuremath{\C^d}}

\newcommand{\Rmxn}{\ensuremath{\R^{m \times n}}}
\newcommand{\Cmxn}{\ensuremath{\C^{m \times n}}}

\newcommand{\Cdxd}{\ensuremath{\C^{d \times d}}}
\newcommand{\Rdxd}{\ensuremath{\R^{d \times d}}}

\newcommand{\kron}{\otimes}
\newcommand{\iid}{\overset{\text{i.i.d.}}{\sim}}
\newcommand{\ee}{\mathrm{e}}
\newcommand{\ii}{\mathrm{i}}

\renewcommand{\th}{\ensuremath{^{\text{th}}}}
\newcommand{\st}{\ensuremath{^{\text{st}}}}
\newcommand{\ow}{\text{otherwise}}

\newcommand{\x}{\ensuremath{\mathbf{x}}}

\newcommand{\y}{\ensuremath{\mathbf{y}}}

\newcommand{\ux}{{x_0}}

\newcommand{\uu}{\underline{u}}

\newcommand{\uX}{{X_0}}

\newcommand{\tbx}{\tilde{x}}

\newcommand{\tF}{\widetilde{F}}
\newcommand{\Lc}{\mathcal{L}}

\newcommand{\Tc}{\mathcal{T}}

\newcommand{\Jc}{\mathcal{J}}
\newcommand{\Rc}{\mathcal{R}}
\newcommand{\Pc}{\mathcal{P}}

\newcommand{\Ac}{\mathcal{A}}
\newcommand{\Nc}{\mathcal{N}}
\newcommand{\CN}{\mathcal{C N}}
\newcommand{\Dc}{\mathcal{D}}
\newcommand{\hz}{\hat{z}}
\newcommand{\hZ}{\hat{Z}}

\newcommand{\bigO}{\mathcal{O}}
\newcommand{\one}{\ensuremath{\mathbbm{1}}}

\newcommand{\dbar}{{\overline{d}}}

\newcommand{\mintheta}{\min\limits_{\theta \in [0, 2\pi]}}

\newcommand{\eit}{\ee^{\ii \theta}}

\newcommand{\conj}[1]{\overline{#1}}
\renewcommand{\subset}{\subseteq}

\def \x {\mathbf x}

\def \y {\mathbf y}

\newcommand{\tX}{\ensuremath{\widetilde{X}}}

\newcommand{\rowmat}[3]{\ensuremath{\begin{bmatrix} #1_{#2} & \cdots & #1_{#3} \end{bmatrix}}}
\newcommand{\rowmatfun}[3]{
  \ensuremath{\noexpandarg
    \begin{bmatrix}
      \StrSubstitute{#1}{@}{{#2}} &
      \cdots &
      \StrSubstitute{#1}{@}{{#3}}
    \end{bmatrix}
  }
}
\newcommand{\rowmatfl}[2]{\ensuremath{\begin{bmatrix} #1 & \cdots & #2 \end{bmatrix}}}

\newcommand{\colmatfl}[2]{\ensuremath{\begin{bmatrix} #1 \\ \vdots \\ #2 \end{bmatrix}}}
\newcommand{\colmat}[3]{\ensuremath{\begin{bmatrix} #1_{#2} \\ \vdots \\ #1_{#3} \end{bmatrix}}}
\newcommand{\colmatfun}[3]{
  \ensuremath{\noexpandarg
    \begin{bmatrix}
      \StrSubstitute{#1}{@}{{#2}} \\
      \vdots \\
      \StrSubstitute{#1}{@}{{#3}}
    \end{bmatrix}
  }
}

\newcommand{\diagcornmat}[2]{
  \ensuremath{
    \begin{bmatrix}
      #1 & \cdots & 0 \\
      \vdots & \ddots & \vdots \\
      0 & \cdots & #2
    \end{bmatrix}
  }
}

\newcommand{\diagmat}[3]{
  \diagcornmat{#1_{#2}}{#1_{#3}}
}

\newcommand{\cornmatfun}[5]{
  \ensuremath{\noexpandarg
    \StrSubstitute{#1}{@r}{{#2}}[\cmftop]
    \StrSubstitute{#1}{@r}{{#3}}[\cmfbot]
    \expandarg
    \begin{bmatrix}
      \StrSubstitute{\cmftop}{@c}{{#4}}  & \cdots & \StrSubstitute{\cmftop}{@c}{{#5}} \\
    \vdots & \ddots & \vdots \\
    \StrSubstitute{\cmfbot}{@c}{{#4}} & \cdots & \StrSubstitute{\cmfbot}{@c}{{#5}}
  \end{bmatrix}
  }
}

\newcommand{\cornmatfunexp}[5]{
  \ensuremath{
    \StrSubstitute{#1}{@r}{{#2}}[\cmftop]
    \StrSubstitute{#1}{@r}{{#3}}[\cmfbot]
    \expandarg
    \begin{bmatrix}
      \StrSubstitute{\cmftop}{@c}{{#4}}  & \cdots & \StrSubstitute{\cmftop}{@c}{{#5}} \\
    \vdots & \ddots & \vdots \\
    \StrSubstitute{\cmfbot}{@c}{{#4}} & \cdots & \StrSubstitute{\cmfbot}{@c}{{#5}}
  \end{bmatrix}
  }
}

\newenvironment{remark}
               {\textbf{Remark.}}
               {}
               
\newenvironment{piecewise}
               {\left\{\begin{array}{r@{,\quad}l}}
               {\end{array}\right.}

\global\newlength\thewidth

%% file: sections/history/history_sec.tex
\section{Introduction}
\input{introduction}
\input{contributions}
\subsection{Related Work}
\label{sec:relatedwork}
\input{relatedwork}

%

%% file: sections/history/introduction.tex

Phase retrieval is the problem of solving a system of equations of the form \begin{equation} y = |A x_0|^2 + \eta, \label{eq:pr_bare} \end{equation} where $x_0 \in \C^d$ is the objective signal, $A \in \C^{D \times d}$ is a known measurement matrix, $\eta \in \R^D$ is an unknown perturbation vector, and $y \in \R^D$ is the vector of measurement data.  Here $|\cdot|^2$ acts componentwise so that for  $v \in \C^n$ we have $|v|^2_j = |v_j|^2$.  In phase retrieval, the goal is to recover an estimate of $x_0$ from knowledge of $y$ and $A$.  We sometimes rephrase the system \eqref{eq:pr_bare} as \begin{equation} y_j = | \langle a_j, x_0 \rangle |^2 + \eta_j,\label{eq:pr}\end{equation} where the $a_j^*$ stand for the rows of $A$ and are referred to as the measurement vectors.  The name \emph{phase retrieval} comes from viewing the $|\cdot|^2$ operation as erasing the phases of the complex-valued measurements $\langle a_j, x_0 \rangle$ and leaving only their magnitudes; solving for $x_0$ may be considered as a way of retrieving this phase information.  We immediately note that this problem contains an unavoidable phase ambiguity, in the sense that, for any solution $x$ and any $\theta \in [0, 2\pi)$, we will have that $\ee^{\ii \theta} x$ is also a solution, as $\abs{\inner{a_j, \eit x}}^2 = \abs{\inner{a_j, x}}^2$.

\subsection{Phase retrieval in imaging science}
  The phase retrieval problem appears in a multitude of imaging systems, since most optical sensors -- most significantly, charge-coupled devices and photographic film -- do not respond to the phase of an incoming light wave.  Rather they respond only to the number and energy of photons arriving at its surface, so they indicate only the intensity (absolute value squared), and not the phase, of the electromagnetic waves to which they are exposed.  This corresponds to our model in \eqref{eq:pr_bare} by imagining that the $i\th$ entry $a_i^* x$ of $Ax \in \C^D$ corresponds to the magnitude and phase of the light arriving at the $i\th$ pixel in an array of sensors.  
 Areas of optics that encounter this problem include astronomy \cite{fienup1987astronomy,walther1963question}, diffraction imaging \cite{
 rodenburg2008diffractive,shechtman2015phase}, laser pulse characterization \cite{bendory2017frog,
 sidorenko2017frog}, electron microscopy \cite{putkunz2012electron}, and x-ray crystallography \cite{bragg1915crystal_structure,marchesini2015coptych,dierolf2008ptych,
 hauptman1953monograph
 }.  Non-optical disciplines that can benefit from solutions to phase retrieval include speech recognition and audio processing \cite{balan2006signal,waldspurger2015cauchy,juang1993speechrec}, blind channel estimation \cite{strohmer2017wtf_deconv1
 }, and self-calibration \cite{strohmer2015self_calib}.

  The practice of these disciplines has produced many creative solutions to particular instances of the phase retrieval problem, and throughout the 20\th\ century the field largely evolved by the invention of \emph{ad hoc} solutions that resolved the data at hand.  
Notably, however, R.W.~Gerchberg and W.O.~Saxton in 1971 \cite{gerchberg1972practical} proposed an algorithm that can be applied to fairly general data, with remarkably minimal assumptions made on the structure of the object $x$ being detected.  This result inspired numerous variants (e.g., \cite{bauschke2003hybrid,bauschke2002phase,elser2003phase,fienup1978reconstruction,takajo1997numerical,takajo1998study,takajo1999further}), each of which empirically improved performance, but none of which produced a solid mathematical theory to explain why or when they would succeed.  Physicists, chemists, and biologists made  astounding scientific achievements in this fashion, but even with all this progress, the community remained largely in want of such a theoretical foundation that could offer reliable solutions in general settings until recent decades.
  
  There are three main questions about phase retrieval problems that the scientific community would wish to answer theoretically: first, in an ideal, noiseless case where $\eta = 0$, for what matrices $A \in \C^{D \times d}$ does the system of equations \eqref{eq:pr_bare} possess a unique solution (up to the known phase ambiguity)?  Second, given a case where a unique solution exists, is there an algorithm that can recover it?  Third, when a recovery process exists, is it stable so that in the presence of noise $\eta \neq 0$, the estimate $x$ does not differ much (or differs to a known degree, as a function of $\norm{\eta}$) from $x_0$?

This paper expands upon the theory of phase retrieval by studying a new class of matrices, that is of particular interest to ptychographic imaging, and an associated recovery algorithm that is proven to solve the system \eqref{eq:pr_bare} with guaranteed stability to noise and with known, competitive computational cost.
\if{\subsection{X-Ray Crystallography}
\label{sec:crystallography}
The history of phase retrieval cannot be told without making mention of x-ray crystallography, the field that first brought scientific interest to this problem and by many metrics its most decorated and fruitful application.  In x-ray crystallography, the goal is to gain an image of the positions of atoms within a molecule by illuminating a crystallized sample with x-rays.  The molecular structure is deduced from the pattern of the radiation diffracted by the sample.  A rough diagram of this setup is shown in \cref{fig:xray_cryst}.
\begin{figure}
  \centering\includegraphics[width=0.9\textwidth]{figs/crystal}
  \caption[Experimental setup for x-ray crystallography]
          {Experimental setup for x-ray crystallography.   {\small
              (Adapted from Cand\'{e}s, Li, and Soltanolkotabi, \emph{Phase Retrieval via Wirtinger Flow: Theory and Algorithms} \cite{candes2015wtf})}}
          \label{fig:xray_cryst}
\end{figure}
This seemingly simple technique has been indispensable for the study of chemistry, biology, and physics, having been used to confirm or identify the arrangements of atoms in a wide variety of important compounds.  Over a dozen discoveries made through x-ray crystallography -- or made in developing the technique -- have been recognized by Nobel Prizes in Physics, Chemistry, and Medicine or Physiology.  Indeed, the first Nobel Prize in Physics was awarded to Wilhelm R\"ontgen in 1901 for his discovery of x-rays.  The 1914 Prize in Physics was conferred upon Max von Laue for discovering the diffraction of x-rays by atomic crystals, and in 1915 William and Lawrence Bragg earned the same distinction for performing the first complete characterizations of atomic crystal structures \cite{galli2014nobel}.  Since the time of these highly esteemed pioneering discoveries, x-ray crystallography has been used to produce accurate molecular models of a number of drugs (e.g., \cite{cell2001antibios, rasmussen2007adrenergic, schindler2000kinase}), including penicillin in Dorothy Crowfoot Hodgkin's Nobel prize-winning work in 1963 \cite{hodgkin1963penicillin}.  It has elucidated several human biological compounds, including innumerable proteins \cite{kimber2003protein, varsani1993isomerase} and human DNA, for whose analysis in 1953 James Watson, Maurice Wilkins, and Francis Crick were awarded the Nobel prize in 1962, relying on the crystallographic images of Rosalind Franklin \cite{franklin1953nature,watson1962nobel_lecture,watson1953nature,wilkins1953nature}.  And this technology remains extremely relevant today, playing an active role in material sciences, where crystallography is being used to characterize the degradation of lithium-ion batteries \cite{hausbrand2015battery, andrej2018battery} and to study carbon nanostructures such as fullerenes \cite{lamb1990carbon, kroto1985fullerene}, whose analysis earned the 1996 Nobel Prize in Chemistry \cite{galli2014nobel}.

Perhaps the Nobel prize, among those awarded for developments related to crystallography, most cogent to the present work is that earned by Herbert A.~Hauptman and Jerome Karle in 1985 for their monograph, \emph{Solution of the Phase Problem, I.  The Centrosymmetric Crystal} \cite{hauptman1953monograph}, in which they presented a direct solution to the phase retrieval problem for periodic, discrete functions which are symmetric about the origin (i.e. centrosymmetric crystals).  Even though this was a special case that the literature of the time had focused on quite heavily, Hauptman and Karle's results eclipsed the somewhat heuristic methods put forth by several authors of the time \cite{harker1948phases,sayre1952shannon}, for the simple fact that they provided a recovery algorithm with a proof of success on a clearly defined set of problem instances.  This is precisely the class of solutions we attempt to expand, and with this thought, we proceed to a statement of the contributions and organization of this dissertation.
}\fi

%% file: sections/history/contributions.tex
Thus, we begin with a description of the application/setting which forms the subject of our analysis, along with a brief description of the phase retrieval strategy whose components we will study in more detail in later sections.  

\subsection{Local  Measurements and Ptychography}
\label{sec:locCorrMeas}
Consider the case where the vectors $a_j$ represent shifts of compactly-supported vectors $m_j, j = 1, \ldots, K$ for some $K \in \N$.  Using the notation $[n]_k:=\clopen{k, k + n}\subset \N$ and defining $[n]:=[n]_1$, we take $x_0, m_j \in \Cd$ with $\supp(m_j)\subset[\delta]\subset [d]$ for some $\delta \in \N$.  We also denote the space of Hermitian matrices in $\C^{k \times k}$ by $\H^k$.  Now we have measurements of the form \begin{equation} (y_\ell)_j = |\langle x_0, S^\ell m_j \rangle|^2, \quad (j, \ell) \in [K] \times P, \label{eq:shift_model} \end{equation} where $P \subset [d]_0$ is arbitrary and $S \in \Cdxd$ is the discrete circular shift operator, namely $(S x)_i = x_{i-1}$.  One can see that \eqref{eq:shift_model} represents the modulus squared of the correlation between $x_0$ and locally supported measurement vectors so we refer to the entries of $y$ as local correlation measurements. 

To see the connection to (a discretized version of) ptychography, consider $\gamma, x_0 \in \C^d$, denoting discretized versions of a known physical mask and unknown sample, respectively. In ptychographic imaging, small regions of a specimen
are illuminated one at a time, often with a physical mask between the specimen and the light source, and an intensity detector captures each of the resulting diffraction patterns. Thus each of the ptychographic measurements is a local measurement, which under certain assumptions (see, e.g., \cite{dierolf2008ptych,goodman2005introfourieroptics,  our_paper}) can be modeled by
\begin{equation}
    (y_\ell)_j = \left \vert \sum_{n=1}^d \gamma_n \, (x_0)_{n-\ell} \,
      \mathbbm e^{-\frac{2\pi \mathbbm i (j-1) (n-1)}d} 
      \right \vert^2, \quad (j,\ell) \in [d] \times  
      [d]_0, 
  \label{eq:1d_ptycho}
\end{equation}
where indexing is considered modulo-$d$. So, $(y_\ell)_j$ is a diffraction measurement corresponding to the $j^{th}$ Fourier mode of a circular $\ell$-shift of the specimen. Note that the use of circular shifts is for convenience only as one can zero-pad $x_0$ and $\gamma$ to obtain the same $(y_\ell)_j$. In practice, one may not need to use all the shifts $\ell \in [d]_0$ as a subset may suffice, { and we also consider this case in this paper}.  Now, defining $m_j \in \C^d$ by \begin{equation} (m_j)_n = \overline{ \gamma_n} \, \mathbbm e^{\frac{2\pi \mathbbm i (j-1) (n-1)}d} \label{eq:ptychm} \end{equation} and rearranging \eqref{eq:1d_ptycho}, we obtain
\begin{align}
    (y_\ell)_j &= \left \vert \sum_{n=1}^d (x_0)_{n-\ell} \, \overline{(m_j)_n} \right \vert^2 \label{eq:loco_measurements} 
        = \lvert \langle S^{-\ell} x_0, m_j \rangle \rvert^2 \notag     = \lvert \langle  x_0, S^{\ell} m_j \rangle \rvert^2 \notag. 
\end{align}
%
Thus \eqref{eq:loco_measurements} shows that ptychography (with $\ell$ ranging over any subset of $[d]_0$) represents a case of the general system seen in \eqref{eq:shift_model}.
Returning to \eqref{eq:shift_model} and following \cite{balan2006signal,candes2012phaselift,IVW2015_FastPhase}, the problem may be lifted to a linear system on the space of $\Cdxd$ matrices.  In particular, we observe that
\begin{align}
	(y_\ell)_j 
	&= \inner{x_0 x_0^*, S^\ell m_j m_j^* S^{\ell *}},
\end{align}
where the inner product above is the Hilbert-Schmidt inner product.  Restricting, for now, to the case $P = [d]_0$,  for every matrix $A \in \Span\{S^\ell m_j m_j^* S^{\ell *}\}_{\ell, j}$ we have $A_{ij} = 0$ whenever $|i - j| \mod d \ge \delta$.  Therefore, we introduce the family of operators $T_k : \Cdxd \to \Cdxd$ given by 
\begin{equation} T_k(A)_{ij} = \begin{piecewise} A_{ij} & |i - j| \mod d < k \\ 0 & \ow. \end{piecewise} \label{eq:T_delta} \end{equation}
 Note that $T_\delta$ is simply the orthogonal projection operator onto $T_\delta(\Cdxd),$ of which $\Span\{S^\ell m_j m_j^* S^{\ell *}\}_{\ell, j}$ is a subspace; therefore, \begin{equation} (y_\ell)_j = \inner{x_0 x_0^*, S^{\ell} m_j m_j^* S^{\ell *}} = \inner{T_{\delta}(x_0 x_0^*), S^{\ell} m_j m_j^* S^{\ell *}}, \quad (j, \ell) \in [K] \times P.  \label{eq:lifted_system} \end{equation}  For convenience, we set $D := K|P|$ to be the total number of measurements and define the map $\Ac : \Cdxd \to \C^D$ \begin{equation}\Ac(X)_{(\ell, j)} = \inner{X, S^{\ell} m_j m_j^* S^{\ell *}}.\label{eq:linear}\end{equation}  

With this in hand, we are prepared to consider our reconstruction strategy, which follows the outline laid out in \cite{IVW2015_FastPhase, our_paper}.  Namely, we will first consider the restriction $\Ac|_{T_{\delta}(\Cdxd)}$ of $\Ac$ to the domain $T_{\delta}(\Cdxd)$, the largest domain on which $\Ac$ may be injective.  Initially, the framework we consider consists of designing measurements $\{m_j\}$ (via the masks $\gamma$) such that $\left.\Ac\right|_{T_{\delta}(\Cdxd)}$ is invertible and then recovering an estimate of $x_0$ from \begin{equation}T_\delta(x_0 x_0^*) =: X_0. \label{equdef:X0} \end{equation}  This recovery process, in turn, is performed by deducing the magnitudes $\abs{x_0}$ and phases $\sgn(x_0):=\frac{x_0}{|x_0|}$ of $x_0$ separately.  This pseudo-algorithm is stated in \cref{alg:pr_basic}.

\begin{algorithm}
\renewcommand{\algorithmicrequire}{\textbf{Input:}}
\renewcommand{\algorithmicensure}{\textbf{Output:}}
\caption{Outline for our phase retrieval algorithm}
\label{alg:pr_basic}
\begin{algorithmic}[1]
    \REQUIRE Measurements $y\in \R^D$ as in \eqref{eq:shift_model}
    \ENSURE An estimate $x$ of $x_0$.
    \STATE Compute the matrix $X = \left.\Ac\right|_{T_\delta(\Cdxd)}^{-1}(y) \in T_{\delta}(\Hd)$.
    \STATE Compute the magnitudes $\abs{x}$ from $X$ (e.g., by methods described in 
    \cref{sec:blocky_block})
     \label{line:mag_rec}
    \STATE Compute the phases $\sgn(x)$ from $X$ (e.g., by methods described in \cref{sec:phase_est})
    \label{line:ang_sync}
    \STATE Return $x = \abs{x} \circ \sgn(x)$.
    \end{algorithmic}
\end{algorithm}

To give an example within the framework of \cref{alg:pr_basic}, similar to the algorithm studied in \cite{IVW2015_FastPhase}, one method of recovering $x_0$ from $X_0 = T_\delta(x_0 x_0^*)$ in the noiseless case would be to simply write $\abs{x}_i = \sqrt{(X_0)_{ii}}$.  To obtain the phases (up to a global shift), we could consider $\sgn(X_0)$ as a matrix of relative phases, in the sense that $\sgn(X_0)_{ij} = (x_0)_i \conj{(x_0)}_j$, allowing us to inductively set $\sgn(x)_1 = 1$ and $\sgn(x)_i = \sgn(x)_{i-1} \sgn(X_{i, i-1})$ for $i = 2,\ldots, d$. Since we will deal with the noisy scenario in this paper, we will strive to develop more sophisticated techniques than these. Indeed, having broken down our main model and recovery algorithm in this manner, we are prepared to chart out the structure and contributions of this paper, keeping in mind that we will generalize the framework of \cref{alg:pr_basic} to handle shifts $\ell$ that don't cover all $[d]_0$, and hence scenarios whereby $\Span\{S^\ell m_j m_j^* S^{\ell *}\}_{\ell, j}$ is a strict subspace of $T_\delta(\C^{d\times d})$. 

\subsection{Organization and Contributions}
\label{sec:organization}
From \cref{alg:pr_basic}, we can identify three main areas of study. The first  is the design of $\{m_j\}$ that permit invertible -- and well conditioned -- linear systems $\mathcal{A}$.  The second and third areas relate to the magnitude and phase recovery steps of lines \ref{line:mag_rec} and \ref{line:ang_sync}. We wish to propose provably efficient and robust algorithms for these sub-tasks, and then combine them to obtain a robust method.

This paper presents contributions in each of the three areas. Indeed, in the first part of the paper, we focus on the case where the full set of shifts is used, and we build on a paper by Iwen, Preskitt, Saab, and Viswanathan \cite{our_paper} that  in-turn improves upon the previous work by \citeauthor*{IVW2015_FastPhase} in \cite{IVW2015_FastPhase} concerning the framework in  \cref{sec:locCorrMeas}.  In \cref{ch:span_fam}, we derive a quickly calculable and exact expression for the condition number of the linear system $\mathcal{A}$, and we leverage this to expand our collection of known spanning families.  Specifically, we will discover that --- in line with ptychographic imaging --- setting $m_j = \gamma \circ f_j^d, j \in [2 \delta - 1]$, where $f_j^d$ is the $j\th$ Fourier vector in $\Cd$ and $\gamma \in \Rd$ has support $[\delta]$ produces an invertible system under a very mild condition.  We further prove that this condition holds for almost all $\gamma$; this result is particularly interesting considering that $\gamma\in \C^d$ has $\delta$ degrees of freedom, but must generate a spanning set for a subspace of dimension $d(2\delta-1)$ when all shifts are taken in \eqref{eq:T_delta}.

\begin{algorithm}
\renewcommand{\algorithmicrequire}{\textbf{Input:}}
\renewcommand{\algorithmicensure}{\textbf{Output:}}
\caption{Phase Retrieval from Local Ptychographic Measurements}
\label{alg:pty_pr}
\begin{algorithmic}[1]
    \REQUIRE A family of masks $\{m_j\}_{j = 1}^D$ of support $\delta$; $s, d, \dbar \in \N$ satisfying $d = \dbar s \ge 2 \delta - 1$.  A $(T_{\delta, s}, d)$-covering $\{J_i\}_{i \in [N]}$.  Measurements $\y = \Ac(\ux \ux^*) + n \in \R^{\dbar D}$, as in \eqref{eq:pty_meas_op} (see \cref{ch:ptychography}).
    \ENSURE $x \in \Cd$ with $x \approx \eit \ux$ for some $\theta \in [0, 2 \pi]$.
    \STATE Compute the matrix $X = \Ac|_{T_{\delta, s}(\Cdxd)}^{-1} y \in T_{\delta, s}(\Hd)$ as an estimate of $T_{\delta, s}(\ux \ux^*)$ (as in \cref{ch:ptychography}).
    \STATE Form the banded matrix of phases, $\tX = \sgn{(X)} \in T_{\delta, s}(\Hd)$, by normalizing the non-zero entries of $X$ (replacing any zero entries in the band with $1$'s).
    \STATE Compute $v$, the eigenvector of $\widetilde{X}$ corresponding to the largest eigenvalue and set $\widetilde{x} = \sgn(v)$ (as in \cref{sec:phase_est}).  
    \STATE Return $x = \BlkMag(X, \{J_i\}) \circ \tbx$ (see \cref{sec:blocky_block}).
\end{algorithmic}
\end{algorithm}

In the second part of the paper, we focus on increasing the match between our model  and the laboratory practices of ptychographers. To that end, we devise algorithms and derive theory that handle the practical setup where large shifts are used. The algorithmic framework we propose and analyze is summarized in \cref{alg:pty_pr}. 
Specifically, in \cref{ch:ptychography} we study the conditioning of the linear system arising from a set of shifts that is smaller than $[d]_0$.  In particular, we consider taking  shifts $\ell = s k, k \in [d / s]$, where $s$ is a fixed step size.  This leads to a new subspace $T_{\delta, s}(\Cdxd) \subset T_\delta(\Cdxd)$, for which we derive condition number estimates in the spirit of \cref{ch:span_fam}. In \cref{sec:blocky_block}, we propose and analyze a magnitude estimation step for \cref{alg:pty_pr} and prove that it is robust to noise. In \cref{sec:phase_est}, we extend the phase-estimation technique of \cite{our_paper} to the setting of large shifts, and show that this technique is robust. Finally, in \cref{sec:pty_std_rec}, we put our results together and prove that \cref{alg:pty_pr} comprises a stable phase retrieval method in the setting of large ptychographic shifts. This result is summarized in the theorem below, which is a slightly weaker but more streamlined version of \cref{thm:pty_rec}.

\begin{theorem}
Let $\Ac$ be the linear system arising from a set of measurements $$\{S^{s\ell}m_jm_j^*S^{-s\ell} \}_{(\ell,j)\in{[d/s]_0\times [D]}}$$
which spans $T_{\delta,s}(\C^{d\times d})$, where the vectors $m_j\in \C^d$ are of support $\delta$, and where $T_{\delta,s}$ is defined in \cref{eq:T_delta_s}. Let $\sigma_{min}^{-1}$ be the smallest singular value of $\Ac$ restricted to $T_{\delta,s}$. The error associated with recovering $x_0$ from the noisy measurements ~$y_{j,\ell}= |\langle S^{s\ell}m_j, \ux \rangle|^2 + n_{j,\ell}$, using \cref{alg:pty_pr} satisfies
\begin{equation}\begin{aligned} \mintheta \norm{x - \eit \ux}_2 &\le 
C \cdot  \sigma_{min} ^{-1} \cdot \left(  
\frac{{d}^{2} \cdot s  }{{\delta}^{5/2}   } \cdot \frac{ \|x_0\|_\infty  }{\min\limits_{j}|(x_0)_j|^2}  ~ + ~\frac{\delta^{-1/2}}{\min\limits_{j}|(x_0)_j|} \right) \cdot \|n\|_2.
\end{aligned} \label{eq:pty_rec}\end{equation}

\end{theorem}



%% file: sections/history/relatedwork.tex
The history of modern algorithmic phase retrieval begins in the 1970's with \cite{gerchberg1972practical} by \citeauthor*{gerchberg1972practical}, where the measurement data corresponded to knowing the magnitude of both the image $x_0$ and its Fourier transform.  This result was famously expanded upon by Fienup \cite{fienup1978reconstruction} later that decade, one significant improvement being that only the magnitude of the Fourier transform of $x_0$ must be known in the case of a signal $x_0$ belonging to some fixed convex set $\mathcal{C}$ (typically, $\mathcal{C}$ is the set of non-negative, real-valued signals restricted to a known domain).  Though these techniques work well in practice and have been popular for decades, they are notoriously difficult to analyze.  These are iterative methods that work by improving an initial guess until they stagnate.  In 2015, Marchesini et al.~proved that alternating projection schemes using generic measurements are guaranteed to converge to the correct solution {\em if provided with a sufficiently accurate initial guess} and algorithms for ptychography were explored in particular \cite{marchesini2015alternating}.  Waldspurger then proved that a spectral initialization  reaches this basin of attraction with high probability using a simple Gaussian suite of measurements \cite{waldspurger2018gerchsax}.  The application of alternating minimizations to sparse phase retrieval has received considerable attention, as well \cite{jagatap2017fast,eldar2017fienup}, although results in \cite{IVW2017_easy} suggest that virtually any phase retrieval method may easily be composed with compressed sensing techniques to target sparse signals.  However, despite this impressive body of work, no global recovery guarantees currently exist for alternating projection techniques using local measurements (i.e., finding a sufficiently accurate initial guess is not generally easy).

Other works have proved probabilistic recovery guarantees when provided with globally supported Gaussian measurements.  Methods for which such results exist vary in their approach, and include convex relaxations \cite{
candes2012phaselift,hassibi2018phasemax,waldspurger2015phasecut}, gradient descent strategies \cite{candes2015wtf
}, graph-theoretic \cite{alexeev2014phase,salanevich2015polarization} and frame-based approaches \cite{balan2009painless, bodmann2013stable,bodmann2017frames}, and variants on conventional alternating minimization ideas \cite{netrapalli2013phase,waldspurger2018gerchsax}.  The approach of non-convex optimization by gradient descent, named \emph{Wirtinger Flow} in its first application to phase retrieval \cite{candes2015wtf}, has enjoyed recent success in a variety of phase retrieval applications \cite{soltanolkotabi2018multiplexed} 
 as well as blind deconvolution \cite{strohmer2017wtf_deconv1} and low-rank matrix recovery \cite{soltanolkotabi2016procrustes}.

Several recovery algorithms achieve theoretical recovery guarantees while using at most $D = \bigO(d \log^4 d)$ masked Fourier coded diffraction pattern measurements, including both {\em PhaseLift} \cite{Candes2014WF,gross2015improved}, and {\em Wirtinger Flow} \cite{candes2015wtf}.  However, until recently, there has been no construction of these measurements that were not randomized, and -- to our knowledge -- the theory has not studied locally supported measurements of the type considered here.  Kueng, Gross, and others have tried to derandomize the constructions for \emph{PhaseLift} in particular by drawing the measurements from certain matrix groups \cite{kueng2015spherical,kueng2016clifford}, but the first completely deterministic, albeit non-local (and possibly non-physical), construction of a measurement system with provable, global recovery via PhaseLift appeared in \cite{kech2018explicit}.

Among the first treatments of local measurements are \cite{bendory2017stft,eldar2014sparse,jaganathan2016stft}, in which it is shown that STFT (short-time Fourier transform \cite{allen1977stft,portnoff1979stft}) measurements with specific properties can allow (sparse) phase retrieval in the noiseless setting, and several recovery methods have been proposed \cite{bendory2018stft,guo2018stft}.  Similarly, the phase retrieval approach from \cite{alexeev2014phase} was extended to STFT measurements in \cite{salanevich2015polarization} in order to produce recovery guarantees in the noiseless setting.  More recently, randomized robustness guarantees were developed for time-frequency measurements in \cite{salanevich2016polarization}.  However, no {\it deterministic} robust recovery guarantees have been proven in the noisy setting for any of these approaches.  Furthermore, none of the algorithms developed in these papers are demonstrated to be empirically competitive with standard alternating projection techniques for large signals when utilizing windowed Fourier and/or correlation-based measurements.  In \cite{IVW2015_FastPhase}, the authors first propose a deterministic measurement scheme  and prove the first deterministic robustness results in the recent literature, although these results treat a ``greedy'' recovery algorithm,  different from the one developed herein, and they obtain weaker recovery guarantees. Very recently, excellent work by Pearlmutter et al. \cite{perlmutter2019inverting} studied the related problem of inverting spectrogram measurements when the mask is locally supported and the signal is bandlimited.  Finally, \cite{melnyk2019phase} studied the ptychographic setup we consider in the second part of the paper, and proved an analogous (but slightly different) result to our \cref{thm:phase_stable} using the magnitude estimation techniques of \cite{our_paper}, rather than the more sophisticated techniques we consider herein.

In the midst of such an active and diverse field of research, the major contributions of our work are that it takes into account the local measurements that match the models for key applications such as ptychography.  In this setting, we have produced a provably fast and stable recovery algorithm for a deterministically stated class of measurement systems that corresponds well to ptychographic imaging.

%% file: sections/history/notation.tex
\label{sec:notation}
We take a moment to gather some of the notation that is used throughout the paper.  \Cref{tab:notation} displays some of the most commonly used objects.  We remark that, in this table and throughout this work, \emph{indices of a vector $x \in \Cn$ or matrix $A \in \Cmxn$ are always taken modulo the appropriate dimension.}  For example, $x_{n + 1} := x_1$ and $A_{00} = A_{mn}$.
%
Given matrices $V_j \in \C^{m_j \times n_j}$ for $j \in [n]$ \[\diag(V_j)_{j = 1}^n = \begin{bmatrix} V_1 & & \\ & \ddots & \\ & & V_n \end{bmatrix} \in \C^{\sum m_j \times \sum n_j}.\]  To conveniently switch between matrices and vectors of different sizes, $\mathcal{R}_d : \bigcup_{k = 1}^\infty \C^k \to \C^d$ is a resizing map, which truncates or zero-pads as appropriate. For $v \in \C^k$ and $i \in [d],$ \[\mathcal{R}_d(v)_i = \left\{\begin{array}{r@{,\quad}l} v_i & i \le k \\ 0 & \text{otherwise} \end{array}\right. \text{for}\ i \in [d].\]  
Similarly, $\mathcal{R}_{m \times n} : \bigcup_{k_1, k_2}^\infty \C^{k_1 \times k_2} \to \C^{m \times n}$ truncates or zero-pads matrices to size $m \times n$.  For $A \in \C^{m \times n}$, we define $\vec(A) \in \C^{m n}$ with  $\vec(A)_{(j - 1) m + i} = A_{i j}$ for $i, j \in [m] \times [n]$.  To invert $\vec$, we use $\mat_{(m, n)} : \C^{m n} \to \C^{m \times n},$ such that $\mat_{(m, n)}(v)_{ij} = v_{(j - 1) m + i}$.  
 
{\newcommand{\env}[1]{\texttt{#1}}\renewcommand{\thefootnote}{\fnsymbol{footnote}}
  \centering \renewcommand{\arraystretch}{1.5}
\input{notation_table}}
{\renewcommand{\thefootnote}{\fnsymbol{footnote}}\footnotetext[3]{We omit the subscript (or superscript) when it is obvious from context.}}

%% file: sections/history/notation_table.tex
\begin{longtable}
  {r|c|>{\setlength\hsize{1.22\hsize}}X|>{\setlength\hsize{.78\hsize}}X}
  \caption[Common operators, objects, and sets]{Common operators, objects, and sets.  Definitions assume $d, i, j, m, n,$ and $k$ are arbitrary elements of $\N$ unless otherwise stated.} \\
  \hline Parameters & Name and Type & Definition & Comments \\ \hline\hline
  \endfirsthead
  \caption[]{Common operators, objects, and sets, Continued}\\
  \hline
  Parameters & Name and Type & Definition & Comments \\ \hline\hline
  \endhead
  \endfoot
  \hline
  \endlastfoot
  & $\one_d \in \Cd$ & $(\one_d)_i = 1$ for $i \in [d]$ & $\one := \one_d$\footnotemark[3] \\
  $A \subset [d]$ & $\one_A^d \in \Cd$ & $(\one_A^d)_i = \chi_A(i)$ & $\one_A := \one_A^d$\footnotemark[3] \\
      $x \in \Cd, k \in [d]$ & $\circop_k(x) \in \C^{d \times k}$ & $\circop_k(x) = \begin{bmatrix} S^0 x & \cdots & S^{k - 1} x \end{bmatrix}$ & $\circop(x) = \circop_d(x)$ \\
  \parbox{2cm}{\raggedleft$\ell \in \Z,$\\$A \in \Cmxn$} & $\diag(A, \ell) \in \Cm$ & $\diag(A, \ell)_i = A_{i, i + \ell}$ & Notation overloaded with $\diag(\cdot)$. \\
  $x \in \Cd$ & $\diag(x) \in \Cdxd$ & $\diag(x) e_i = x_i e_i$ & Also written $D_x$ or $\diag(x_j)_{j=1}^d$. \\
  & $e_i^n \in \Rn$ & $e_i^n = I_n e_i$ & Usually  write $e_i$. \\
  & $F_d \in \Cdxd$ & $(F_d)_{ij} = \frac{1}{\sqrt{d}} \omega_d^{(i-1)(j-1)}$ &  $F_d$ is unitary.  $F := F_d$\footnotemark[3] \\
  & $f_j^d \in \Cd$ & $f_j^d = F_d e_j$ & $f_j := f_j^d$\footnotemark[3] \\
  & $\Hd \subset \Cdxd$ & $A \in \Hd$ iff $A = A^*$ & Hermitian matrices\\
  & $I_d \in \Rdxd$ & $I_d x = x$ & $I := I_d$\footnotemark[3] \\
  & $\abs{i - j} \mod d$ & $\min\{k \ge 0 : k \equiv \pm(i - j) \mod d\}$ & \\
  $n \in \Z$ & $[k]_n \subset \N$ & $\clopen{k, k + n} \cap \Z$ & We define $[k] = [k]_1$. \\  
  & $m \mod_k n \in [n]_k$ &  $p$  $\in[n]_k$ satisfying $p \equiv m \mod n$. & $m \mod n = m \mod_0 n$ \\
  & $\sym^d \subset \Rdxd$ & $\sym^d = \Rdxd \cap \Hd$ & Symmetric matrices \\
  & $S_d \in \Rdxd$ & $(S_d x)_i = x_{i-1}$ & $S := S_d$\footnotemark[3] \\
  & $R_d \in \Rdxd$ & $(R_d x)_i = x_{2-i}$ & $R := R_d$\footnotemark[3] \\
  $x, y \in \Cd$ & $x \circ y \in \Cd$ & $(x \circ y)_i = x_i y_i$ & Hadamard product \\
  $A \subset B$ & $\chi_A : B \to \{0, 1\}$ & $\chi_A(x) = \begin{piecewise} 1 & x \in A \\ 0 & \ow \end{piecewise}$ & \\
  & $\omega_d \in \C$ & $\omega_d = \ee^{\frac{2 \pi \ii}{d}}$ & $\omega := \omega_d$\footnotemark[3]
  \label{tab:notation}
\end{longtable}

%% file: sections/meas/meas_sec.tex
\label{ch:meas_sec}
\label{ch:meas}
\label{sec:meas_intro}
\input{meas_intro}
\subsection{Invertibility and condition numbers}
\label{sec:con_number}
\input{con_number}

\subsection{Explicit Examples of Spanning Families}
\label{sec:expl_span_fam}
\input{expl_span_fam}
\subsection{Inverting $\Ac$}
\label{sec:meas_expl_inv}
\input{meas_expl_inv}

%% file: sections/meas/meas_intro.tex
Herein we focus on the setup described in \Cref{sec:locCorrMeas}, while also accounting for noise as in \eqref{eq:pr}. Thus, by design, it is clear that the vectors $\{S^{\ell} m_j m_j^* S^{- \ell}\}$ are all contained in the subspace $T_\delta(\H^d)$ of $\H^d$ (defined in \eqref{eq:T_delta}), where $d$ is the ambient dimension (meaning $x_0, m_j \in \Cd$) and $\delta$ is the support size of the masks (so $\supp(m_j) \subset [\delta]$).  
%
%
%
Following the framework of \cref{alg:pr_basic} we now focus on the linear system \eqref{eq:linear}. In \cite{our_paper, IVW2015_FastPhase}, a total of  two examples of collections of vectors $\{m_j\}_{j \in [2 \delta - 1]} \subset \Cd$ such that this linear system was invertible on $T_\delta(\H^d)$ were given.  Considering that one of the major contributions of \cref{alg:pr_basic} is that it admits measurement models that intend to replicate laboratory conditions, our knowledge of which vectors are compatible with our algorithm and with the theory built for it is critical in promoting applicability.  In \Cref{sec:con_number}, we study the conditioning of the linear system in \eqref{eq:linear} as a function of the set of masks $\{m_j\}_{j = 1}^D$, resulting in \Cref{prop:meas_cond}.  This result also gives us a description of all sets of masks $\{m_j\}_{j = 1}^D$ (and $\gamma \in \Rd$ that generate such families) that are capable of spanning $T_\delta(\H^d)$, in the sense that we have a checkable condition that indicates whether the linear system of \eqref{eq:linear} is invertible.  
%
%
In \Cref{sec:expl_span_fam}, we provide a few examples of explicit $\gamma \in \Rd$ that are proven to satisfy the conditions to span $T_\delta(\H^d)$.  We consider the act of inverting $\Ac$ from a practical perspective in \Cref{sec:meas_expl_inv}.  We write its inverse explicitly and analyze the runtime of calculating $\Ac^{-1}(y)$ in \Cref{sec:inv_runtime}.  

\subsection{Preliminaries}
Before we begin these analyses, we introduce some definitions.  We say that $\{m_j\}_{j = 1}^D \subset \C^d$ is a \emph{local measurement system} or \emph{family of masks} of support $\delta$ if $1 \in \supp(m_j)$ and $\supp(m_j) \subset [\delta]$ for each $j$.
If we further have that each $m_j$ satisfies $m_j = \Rc_d(\sqrt{K} f_j^K) \circ \gamma$ for some $K \ge \max(\delta, D), \gamma \in \C^d$ satisfying $\supp(\gamma) = [\delta]$, then we call $\{m_j\}_{j = 1}^D$ a \emph{local Fourier measurement system} of support $\delta$ with mask $\gamma$ and modulation index $K$.  If $K = D = 2 \delta - 1$, then we simply refer to $\{m_j\}_{j = 1}^D$ as a \emph{local Fourier measurement system} of support $\delta$ with mask $\gamma$.  We add that, if we say that $\{m_j\}_{j = 1}^D$ is a local Fourier measurement system with support $\delta$ and mask $\gamma$, this implies an assertion that $\supp(\gamma) = [\delta]$. 

Given a local measurement system $\{m_j\}_{j = 1}^D$ in $\C^d$, the associated \emph{lifted measurement system} is the set $\mathcal{L}_{\{m_j\}} = \{S^{\ell} m_j m_j^* S^{- \ell}\}_{(\ell, j) \in [d]_0 \times [D]} \subset \C^{d \times d}$.  We then say that a family of masks $\{m_j\}_{j = 1}^D \subset \C^d$ of support $\delta$ is a \emph{spanning family} if $\Span_\R \Lc_{\{m_j\}} = T_\delta(\H^d)$.
The \emph{measurement operator} associated with a local measurement system is the operator
  \begin{gather}
    \mathcal{A} : T_\delta(\C^{d \times d}) \to \C^{[d]_0 \times [D]} \nonumber \\
    \mathcal{A}(X)_{(\ell, j)} = \langle S^{\ell} m_j m_j^* S^{-\ell}, X \rangle. \label{eq:meas_op}
  \end{gather}
  The \emph{canonical matrix representation} of $\Ac$ is the matrix $A \in \C^{d D \times d (2 \delta - 1)},$ defined by
  \begin{equation}
    \left(A \colmatfun{\diag(X, @)}{1 - \delta}{\delta - 1}\right)_{(j - 1) d + \ell} = \Ac(X)_{(\ell - 1, j)}.
    \label{eq:meas_mat}
  \end{equation}
  For convenience, we define the \emph{diagonal vectorization operator} $\Dc_I : \C^{d \times d} \to \C^{\abs{I} d}$ for any collection $(\ell_i)_{i = 1}^{\abs{I}} = I \subset [d]$  and $\Dc_k : \C^{d \times d} \to \C^{(2 k - 1) d}$ for any integer $k \le \frac{d + 1}{2}$ by
  \begin{align}
    \Dc_I(X) &= \colmatfun{\diag(X, \ell_@)}{1}{\abs{I}}
    \label{eq:diag_vec_set} \\
    \Dc_k(X) = \Dc_{[2 k - 1]_{1 - k}}(X) &= \colmatfun{\diag(X, @)}{1 - k}{k - 1},
    \label{eq:diag_vec}
  \end{align}
  so that \eqref{eq:meas_mat} becomes $A \Dc_{\delta} (X)_{(j - 1) d + \ell} = \Ac(X)_{(\ell - 1, j)}$.  We remark that, when $2 k - 1 \le d$, $\Dc_k$ is invertible on $T_k(\C^{d \times d})$, and for $v \in \C^{d (2 k - 1)}$, we use $\Dc_k^{-1}(v)$ or $\Dc_k^*(v)$ to represent the matrix in $T_k(\C^{d \times d})$ whose diagonals are given by the $2k - 1$ distinct $d$-length blocks of $v$.

%% file: sections/meas/con_number.tex
 We now calculate the singular values, and therefore the condition number, of the measurement operator $\Ac$ for an arbitrary local measurement system $\{m_j\}_{j = 1}^d$ in \cref{thm:meas_cond}, the main result of this section.  We remark that this is an important element in using the error bounds proven in \cite{our_paper},
 since they all unavoidably rely on the condition number $\kappa$ and singular values of the linear system solved in line 1 of \cref{alg:pr_basic}, and leaving this quantity unknown and unestimated renders these bounds far less useful.  We also remark that, in this section, as in \cref{sec:expl_span_fam}, we focus on calculations of the condition numbers $\kappa$, rather than the smallest singular value $\sigma_{\min}$ of the linear systems, since $\kappa$ is scale invariant.  Indeed, if $\{m_j\}_{j \in [D]}$ is a local measurement system with $\sigma_{\min}^{-1} = s,$ then $\{t m_j\}_{j \in [D]}$ has $\sigma_{\min}^{-1} = \frac{s}{t^2},$ so it would appear that simply making $t$ large could arbitrarily improve how well the estimate $X$ arising from line 1 of \cref{alg:pr_basic} approximates $T_\delta{(x_0x_0^*)}$; of course, this action would correspond to simply multiplying the observed measurements by the scalar $t^2$, and slips in its advantage by ignoring that $\norm{n}_2$ would also scale as $t^2$ in such a case.  This process clearly buys us nothing, so the $\sigma_{\min}^{-1}$ inequalities fail, in this way, to consider a sense of scale between $\norm{n}_2$ and $\norm{\Ac(x_0 x_0^*)}_2$.  Therefore, by referencing $\kappa$ and $\SNR = \frac{\norm{\Ac(x_0 x_0^*)}}{\norm{n}_2}$ 
 we have inequalities that more accurately 
 describe the relationship between the design of the linear system and the accuracy of the estimate produced by \cref{alg:pr_basic}.

We emphasize further that, perhaps equally or even more significantly, this result gives a description of all local measurement systems that are usable for phase retrieval in \cref{alg:pr_basic}, since we may simply check whether $\{m_j\}_{j \in [D]}$ leads to a system with any singular values of $0$.  This is an important addition to the framework of \cite{our_paper}, since previously we only possessed two examples of families of masks that produced invertible linear systems. Before stating the result, we introduce the operator $R$ which maps $x\in \C^d$ to $(x_{2-i})_{i=1}^{d}$ (where indexing is mod $d$ as always). We also introduce the  \emph{interleaving operators} $P^{(d, N)} : \C^{dN} \to \C^{dN}$ for any $d, N \in \N$, each of which is a permutation defined by \begin{equation} (P^{(d, N)} v)_{(i - 1)N + j} = v_{(j - 1)d + i}.\label{eq:interleave_def}\end{equation}  We can view this as beginning with $v \in C^{dN}$ written as $N$ blocks of $d$ entries, and interleaving them into $d$ blocks each of $N$ entries.

\begin{theorem} \label{thm:meas_cond} \label{prop:meas_cond}
  Given a family of masks $\{m_j\}_{j \in [D]}$ of support $\delta \le \frac{d + 1}{2}$, we define $g_m^j = \diag(m_j m_j^*, m),$ \[H = P^{(d, D)} \cornmatfun{R \conj{g}^@r_@c}{1}{D}{1 - \delta}{\delta - 1}
  ,\] and $M_j = \sqrt{d}\left(f_j^d \otimes I_D\right)^* H$.  Then the singular values of $\Ac$ as defined in \eqref{eq:meas_op} are $\{\sigma_i(M_j)\}_{(i, j) \in [D] \times [d]}$ and its condition number is \[\kappa(\mathcal{A}) = \dfrac{\max\limits_{i \in [d]} \sigma_{\max} (M_i)}{\min\limits_{i \in [d]} \sigma_{\min} (M_i)}.\]
\end{theorem}
Although \cref{thm:meas_cond} is satisfyingly general, perhaps the most useful result in this section is the strictly narrower \cref{prop:span_fam_cond} which greatly generalizes the mask construction in \cite{IVW2015_FastPhase}. 

%

\begin{proposition}
\label{prop:span_fam_cond}
\label{prop:gam_fam_cond}
  Let $\{m_j\}_{j = 1}^D \subset \C^d$ be a local Fourier measurement system with support $\delta$, mask $\gamma$, and modulation index $K$, where $D = 2\delta-1\leq d$.  Let $\Ac$ be the associated measurement operator as in \eqref{eq:meas_op}, with canonical matrix representation $A$ as in \eqref{eq:meas_mat}. 
\begin{itemize}
\item[]If $K = D,$ then the condition number of $\mathcal{A}$ is \begin{equation}\kappa(\mathcal{A}) = \dfrac{d^{-1/2} \norm{\gamma}_2^2}{\min\limits_{m \in [\delta]_0, j \in [d]} \lvert F_d^* (\gamma \circ S^{-m} \gamma)_j \rvert}.\label{eq:clean_cond}\end{equation}  

\item[] If $K > D$, the condition number is bounded by \begin{equation}\kappa(\mathcal{A}) \le \dfrac{d^{-1/2} \norm{\gamma}_2^2}{\min\limits_{m \in [\delta]_0, j \in [d]} \lvert F_d^* (\gamma \circ S^{-m} \gamma)_j \rvert} \kappa(\tF_K), \label{eq:messy_cond}\end{equation} 
where $\tF_K \in \C^{D \times D}$ is the $D \times D$ principal submatrix of $F_K$.
\end{itemize}
\end{proposition}

We recall that the design of local Fourier measurement systems is motivated by the application of ptychography -- in this type of laboratory setup, $\gamma$ can represent a mask or ``illumination function,'' describing the intensity of radiation applied to each segment of the sample -- so it is appropriate to the end user that our simplest and most conveniently applied result pertains to a realistic, broad class of local measurement systems.  In particular, we can now  determine  masks/illumination functions that are admissible for our phase retrieval algorithm: following almost immediately from \cref{prop:span_fam_cond}, we have sufficient conditions for a local Fourier measurement system to be a spanning family, which we state in \cref{cor:gam_fam_span}.\medskip

\begin{corollary}
  Let $\{m_j\}_{j \in D}$ be a local Fourier measurement system of support $\delta$ with mask $\gamma \in \Rd$ and modulation index $K$, where $D = 2 \delta - 1$.  Then $\{m_j\}_{j \in [D]}$ is a spanning family if $F_d^*(\gamma \circ S^{-m} \gamma)_j \neq 0$ for all $m \in [\delta]_0, j \in [d]$ and $K \ge D$.
  \label{cor:gam_fam_span}
\end{corollary}
\noindent\begin{remark} \noindent 
The condition in \cref{cor:gam_fam_span} for a local Fourier measurement system to be a spanning family is generic, in the sense that it fails to hold only on a subset of $\R^d$ with Lebesgue measure zero, except possibly when $\delta > d / 2$.  We consider that, whenever $m \neq d / 2$, the set of $\gamma \in \R^d$ giving at least one zero in $F_d^*(\gamma \circ S^{-m} \gamma)$ is a finite union of zero sets of non-trivial quadratic polynomials and hence a set of zero measure.  
Indeed, when $m \neq d / 2$, we have that \[F_d^*((e_1 + e_{m + 1}) \circ S^m(e_1 + e_{m + 1}))_k = f_k^{d*} e_{m + 1} = \omega^{-m(k-1)},\] so $\gamma \mapsto F_d^*(\gamma \circ S^m \gamma)_k$ is a non-zero, homogeneous quadratic polynomial and therefore has a zero locus of measure zero.
\end{remark}

The proofs of \cref{thm:meas_cond} and \cref{prop:span_fam_cond} require some preliminary work in defining and studying a few new operators pertaining to the structure of \eqref{eq:meas_op}.  These definitions and a number of results concerning them are contained in \cref{sec:interlemma}.  
\if{\subsection{Interleaving Operators and Circulant Structure {\color{red} Probably send to an appendix}}
\label{sec:interlemma}

To set the stage for the proof, we introduce a certain collection of permutation operators and study their interactions with circulant and block-circulant matrices.  The structure we identify here will be of much use to us in unraveling the linear systems we encounter in our model for phase retrieval with local correlation measurements.

For $\ell, N_1, N_2 \in \N, v \in \C^{\ell N_1}, k \in [\ell],$ and $H \in \C^{\ell N_1 \times N_2}$, we define the block circulant operator $\circop^{N_1}$ by
\begin{align*}
  \circop_k^{N_1}(v) &= \begin{bmatrix} v & S^{N_1} v & \cdots & S^{(k - 1)N_1} v \end{bmatrix} \\
  \circop_k^{N_1}(H) &= \begin{bmatrix} H & S^{N_1} H & \cdots & S^{(k - 1) N_1}H \end{bmatrix},
\end{align*}
where, as with $\circop(\cdot)$, when we omit the subscript we define $\circop^{N_1}(H) = \circop_\ell^{N_1}(H)$ and $\circop^{N_1}(v) = \circop_\ell^{N_1}(v)$.  We now proceed with the following lemmas; the first establishes the inverse of $P^{(d, N)}$.

\begin{lemma} \label{lem:interleave_inverse}
  For $d, N \in \N,$ we have \[(P^{(d, N)})^{-1} = P^{(d, N) *} = P^{(N, d)}.\]
\end{lemma}

\begin{proof}[Proof of \cref{lem:interleave_inverse}]
Simply take $v \in \C^{d N}$ and calculate, for $i \in [d], j \in [N]$,
  \begin{align*}
    (P^{(d, N)} P^{(N, d)} v)_{(i - 1) N + j} &= (P^{(d, N)} (P^{(N, d)} v))_{(i - 1) N + j} 
    = (P^{(N, d)} v)_{(j - 1) d + i} 
    = v_{(i - 1) N + j},
  \end{align*}
  with these equalities coming from the definition in \eqref{eq:interleave_def}.  
\end{proof}

We now observe some useful ways in which the interleaving operators commute with the construction of circulant matrices.

\begin{lemma}\label{lem:interleave}
  Suppose $V_i \in \C^{k \times n}, v_{ij} \in \C^k, w_j \in \C^{k N_1}$ for $i \in [N_1], j \in [N_2]$ and
  \begin{gather*}
    M_1 = \begin{bmatrix} \circop(V_1) \\ \vdots \\ \circop(V_{N_1}) \end{bmatrix},\quad
    M_2 = \begin{bmatrix} \circop^{N_1}(w_1) & \cdots & \circop^{N_1}(w_{N_2}) \end{bmatrix},\ \text{and} \\
    M_3 = \begin{bmatrix} \circop(v_{11}) & \cdots & \circop(v_{1 N_2}) \\ \vdots & \ddots & \vdots \\ \circop(v_{N_1 1}) & \cdots & \circop(v_{N_1 N_2}) \end{bmatrix}.\end{gather*}
  Then
  \begin{align}
    P^{(k, N_1)} M_1 &= \circop^{N_1}\left(P^{(k, N_1)} \begin{bmatrix} V_1 \\ \vdots \\ V_{N_1} \end{bmatrix}\right) \label{eq:M_1} \\
    M_2 P^{(k, N_2)*} &= \circop^{N_1}\left(\begin{bmatrix} w_1 & \cdots & w_{N_2} \end{bmatrix}\right) \label{eq:M_2} \\
    P^{(k, N_1)}M_3P^{(k, N_2)*} &= \circop^{N_1}\left(P^{(k, N_1)} \begin{bmatrix} v_{11} & \cdots & v_{1 N_2} \\ \vdots & \ddots & \vdots \\ v_{N_1 1} & \cdots & v_{N_1 N_2} \end{bmatrix}\right). \label{eq:M_3}
  \end{align}
\end{lemma}

\begin{proof}[Proof of lemma \ref{lem:interleave}]
  We index the matrices to check the equalities.  For \eqref{eq:M_1}, we take $(a, b, \ell, j) \in [d] \times [N_1] \times [k] \times [n]$ and have 
  \begin{align*}
    (P^{(k, N_1)} M_1)_{(a-1)N_1 + b, (\ell - 1) n + j} &= (M_1)_{(b - 1) k + a, (\ell - 1)n + j} 
    = \begin{bmatrix} S^{\ell - 1} V_1 \\ \vdots \\ S^{\ell - 1} V_{N_1} \end{bmatrix}_{(b - 1)k + a, j} 
    \\ &
    = (S^{\ell - 1}V_b)_{a, j} = (V_b)_{a + \ell - 1, j}
  \end{align*}
  and
  \begin{align*}
    \circop^{N_1}\left(P^{(k, N_1)} \begin{bmatrix} V_1 \\ \vdots \\ V_{N_1} \end{bmatrix}\right)_{(a-1)N_1 + b, (\ell - 1) n + j} &= \left(P^{(k, N_1)} \begin{bmatrix} V_1 \\ \vdots \\ V_{N_1} \end{bmatrix}\right)_{(a - 1)N_1 + b + (\ell-1)N_1, j} \\
    &= 
    (V_b)_{a + \ell - 1, j}
  \end{align*}
  For \eqref{eq:M_2}, we take $(a, b, j) \in [k] \times [N_2] \times [k N_1]$ and have
 \[(P^{(k, N_2)} M_2^*)_{(a - 1)N_2 + b, j} = (M_2)_{j, (b - 1) k + a} = (w_b)_{j + (a - 1)N_1}\]
  and
  \[\left(\circop^{N_1}\left(\begin{bmatrix} w_1 & \cdots & w_{N_2} \end{bmatrix}\right)\right)_{j, (a - 1)N_2 + b} = (S^{N_1(a - 1)} w_b)_j = (w_b)_{j + N_1(a - 1)},\]
and  \eqref{eq:M_3} follows immediately by combining \eqref{eq:M_1} and \eqref{eq:M_2}.
\end{proof}

\Cref{lem:interkron} introduces a few useful identities relating the interleaving operators to kronecker products.

\begin{lemma}\label{lem:interkron}
  For $v \in \C^N, V = \rowmat{V}{1}{\ell} \in \C^{N \times \ell}, A = \rowmat{A}{1}{m} \in \C^{d \times m}$, and $B_i \in \C^{m \times k}, i \in [\ell]$, we have
  \begin{align}
    P^{(d, N)} (v \kron A) &= A \kron v
    \label{eq:interkron_vec} \\
    P^{(d, N)} (V \kron A) &= \rowmat{A \kron V}{1}{\ell}
    \label{eq:interkron_mat} \\
    P^{(d, N)} (V \kron A) P^{(\ell, m)} &= A \kron V
    \label{eq:interkron_swap} \\
    (V \kron A) \diagmat{B}{1}{\ell} &= \rowmatfun{V_@ \kron A B_@}{1}{\ell}
    \label{eq:kron_diag}
  \end{align}
\end{lemma}

\begin{proof}[Proof of \cref{lem:interkron}]
  For \eqref{eq:interkron_vec}, we see that, for $i, j, k \in [d] \times [N] \times [m]$, we have
  \begin{align*}
    (P^{(d, N)} v \otimes A)_{(i - 1) N + j, k} &= (v \otimes A)_{(j - 1) d + i, k} \\
    &= v_j A_{i k}, \ \text{while} \\
    (A \kron v)_{(i - 1) N + j, k} &= A_{i k} v_j,
  \end{align*}
  and \eqref{eq:interkron_mat} follows by considering that $V \otimes A = \rowmatfun{V_@ \kron A}{1}{\ell}.$  To get \eqref{eq:interkron_swap}, we trace the positions of columns, considering that $(V \kron A) e_{(i - 1) m + j} = V_j \kron A_i$.  From \eqref{eq:interkron_mat}, we observe that $P^{(d, N)} (V \kron A) e_{(i - 1) m + j} = A_j \kron V_i,$ so \begin{align*}
    P^{(d, N)} (V \kron A) P^{(m, \ell)} e_{(j - 1) \ell + i} &= P^{(d, N)} (V \kron A) e_{(i - 1) m + j} \\
    &= A_j \kron V_i = (A \kron V) e_{(j - 1) \ell + i}.
  \end{align*}
 As for \eqref{eq:kron_diag}, we remark that \begin{gather*} (V \kron A) \diagmat{B}{1}{\ell}
    = (V \kron A) \rowmatfun{e^\ell_@ \kron B_@}{1}{\ell}  \\ = \rowmatfun{(V \kron A) (e_@^\ell \kron B_@)}{1}{\ell} = \rowmatfun{V_@ \kron A B_@}{1}{\ell},\end{gather*} as desired.
\end{proof}

The following lemma is a standard result (e.g., Theorem 13.26 in \cite{laub2004matrix}) regarding the kronecker product.

\begin{lemma}\label{lem:kronvec}
  We have $\vec(A B C) = (C^T \kron A) \vec(B)$ for any $A \in \C^{m \times n}, B \in \C^{n \times p}, C \in \C^{p \times k}.$  In particular, for $a, b \in \Cd, \vec a b^* = \conj{b} \kron a$, and \begin{equation}\label{id:kronsimp} \vec E_{jk} (\vec E_{j' k'})^* = E_{k k'} \kron E_{j j'}. \end{equation}
\end{lemma}
The next lemma covers the standard result concerning the diagonalization of circulant matrices, as well as a generalization to block-circulant matrices.
\begin{lemma}
  For any $v \in \Cd$, we have \begin{equation} \circop(v) = F_d \diag(\sqrt{d} F_d^* v) F_d^* = \sqrt{d} \sum_{j = 1}^d (f_j^{d *} v) f_j^d f_j^{d *} \label{eq:circ_dft_diag} \end{equation}
  Suppose $V \in C^{k N \times m}$, then $\circop^N(V)$ is block diagonalizable by \begin{equation} \circop^N(V) = \left(F_k \otimes I_N\right) \left(\diag(M_1, \ldots, M_k)\right) \left(F_k \otimes I_m\right)^*,  \label{eq:circ_dft_blk} \end{equation} where \begin{equation} \sqrt{k}\left(F_k \otimes I_N\right)^* V = \begin{bmatrix} M_1 \\ \vdots \\ M_k \end{bmatrix}, \quad \text{or} \quad M_j = \sqrt{k} (f_j^k \otimes I_N)^* V \label{eq:M_ell}\end{equation} \label{lem:circ_diag}
\end{lemma}

\begin{proof}[Proof of lemma \ref{lem:circ_diag}]
  The diagonalization in \eqref{eq:circ_dft_diag} is a standard result: see, e.g., Theorem 7 of \cite{gray2006circulant}.

  To prove \eqref{eq:circ_dft_blk}, we set $V_i$ to be the $k \times m$ blocks of $V$ such that $V^* = \rowmat{V^*}{1}{k}$ and begin by observing that, for $u \in \C^k$ and $W \in \C^{m \times p}$, the $\ell\th$ $k \times p$ block of $\circop^N(V)(u \otimes W)$ is given by \[\left(\circop^N(V)(u \otimes W)\right)_{[\ell]} = \sum_{i = 1}^k u_i (S^{N (i - 1)}V)_\ell W = \sum_{i = 1}^k u_i V_{\ell - i + 1} W.\]  Taking $u = f_j^k$ and $W = I_m$, this gives \begin{align*} \left(\circop^N(V)(f_j^k \otimes I_m)\right)_{[\ell]} &= \frac{1}{\sqrt{k}}\sum_{i = 1}^k \omega_k^{(j - 1) (i - 1)} V_{\ell - i + 1} I_m \\ &= \frac{1}{\sqrt{k}} \omega_k^{(j - 1) (\ell - 1)} \sum_{i = 1}^k \omega_k^{-(j - 1)(i - 1)} V_i \\ &= (f_j^k)_\ell \left(\sqrt{k} (f_j^k \otimes I_N)^* V \right) = (f_j^k)_\ell M_j. \end{align*}  This relation is equivalent to having \[\circop^N(V) (f_j^k \otimes I_m) = (f_j^k \otimes M_j) = (f_j^k \otimes I_N) M_j,\] which is the statement of the lemma.
\end{proof}

Lemma \ref{lem:circ_diag} immediately gives the following corollary regarding the conditioning of $\circop^N(V)$, with which we return to considering spanning families of masks.
\begin{corollary} \label{cor:circ_diag_condition}
  With notation as in lemma \ref{lem:circ_diag}, the condition number of $\circop^N(V)$ is \[\dfrac{\max\limits_{i \in [k]} \sigma_{\max} (M_i)}{\min\limits_{i \in [k]} \sigma_{\min} (M_i)}.\]
\end{corollary}
}\fi
We begin with the proof of \cref{thm:meas_cond}, of which \cref{prop:span_fam_cond} is a special case.

\begin{proof}[Proof of \cref{thm:meas_cond}]
  We consider the rows of the matrix $A$ representing the measurement operator $\mathcal{A}$, defined in \eqref{eq:meas_mat} and \eqref{eq:meas_op}.  We vectorize $X \in T_\delta(\C^{d \times d})$ by its diagonals with $\Dc_\delta$, as in \eqref{eq:diag_vec} and set $\chi_m = \diag(X, m), m = 1 - \delta, \ldots, \delta - 1$ and recall that $g_m^j = \diag(m_j m_j^*, m)$. Then
  \begin{align*}
    \mathcal{A}(X)_{(\ell, j)} &= \langle S^{\ell} m_j m_j^* S^{-\ell}, X \rangle 
    = \sum_{m = 1 - \delta}^{\delta - 1} \langle S^{\ell} g_m^j, \chi_m \rangle,
  \end{align*}
  so that the definition of $A$ in \eqref{eq:meas_mat} immediately yields its $(j - 1) d + \ell\th$ row as $\rowmatfun{g_@^{j*} S^{1 - \ell}}{1 - \delta}{\delta - 1}$.  Transposing this expression, we see that the $(j - 1)d + 1\st$ through $(j - 1) d + d\th$ rows of $A$ compose $\rowmatfun{\circop(g^j_@)^*}{1 - \delta}{\delta - 1}$. 
  Together with $\circop(v)^* = \circop(R \conj{v})$,  $A$ is the block matrix given by \[A = \cornmatfun{\circop(g_@c^@r)^*}{1}{D}{1 - \delta}{\delta - 1} = \cornmatfunexp{\circop(R \conj{g}_@c^@r)}{1}{D}{1 - \delta}{\delta - 1},\] which may be transformed, by \eqref{eq:M_3} of \cref{lem:interleave}, to
  \begin{equation}
    P^{(d, D)} A P^{(d, 2 \delta - 1)*} = \circop^D\left(P^{(d, D)} \cornmatfun{R \conj{g}_@c^@r}{1}{D}{1 - \delta}{\delta - 1} \right) = \circop^D(H).
    \label{eq:interleaved_meas}
  \end{equation}
%
Quoting \cref{cor:circ_diag_condition} establishes the theorem.
\end{proof}

We are now able to prove \cref{prop:span_fam_cond}.

\begin{proof}[Proof of \cref{prop:span_fam_cond}]
  We begin by remarking that, for any $j, m \in [d]$, the Cauchy-Schwarz inequality gives us
  \begin{align*}
    f_j^{d*} (\gamma \circ S^{-m} \gamma) &= \left(D_{f_j^{d*}} \gamma\right) \cdot \left(S^{-m} \gamma\right) 
    \le \frac{1}{\sqrt{d}} \norm{\gamma}_2^2.
  \end{align*}
  Observing that $f_1^{d*}(\gamma \circ \gamma) = \frac{1}{\sqrt{d}} \norm{\gamma}_2^2$, we have that \begin{equation} \max_{(j, m) \in [d] \times [\delta]_0} F_d^* (\gamma \circ S^{-m} \gamma)_j = \frac{1}{\sqrt{d}} \norm{\gamma}_2^2. \label{eq:max_sv} \end{equation}
Recall that $D = 2 \delta - 1 \le d$ and set $\tF_K \in \C^{2 \delta - 1 \times 2 \delta - 1}, (\tF_K)_{ij} = \frac{1}{\sqrt{K}}\omega_K^{(i-1)(j-\delta)}$ to be the principal submatrix of $\diag(\sqrt{K} f^K_{2 - \delta}) F_K$, let $v_j=\Rc_d(\sqrt{K} f_j^K)$.  For a local Fourier measurement system (i.e., $K=2\delta-1$), we have \begin{equation} \begin{aligned} g_m^j &= \diag(m_j m_j^*, m) = \diag((\gamma \circ v_j) (\gamma \circ v_j)^*, m) \\ &= \diag(D_{v_j} \gamma \gamma^* D_{\overline{v_j}}, m) = \omega_K^{-m(j - 1)} \diag(\gamma \gamma^*, m), \end{aligned} \label{eq:gam_diag}\end{equation} so, setting $g_m = \diag(\gamma \gamma^*, m),$ we have $g_m^j = \diag(m_j m_j^*, m) = \omega_K^{-m(j - 1)} g_m$.  Therefore, we label the $2 \delta - 1 \times 2 \delta - 1$ blocks of $H$ (where $H$ is as in \cref{thm:meas_cond}) by $H^* = \begin{bmatrix} H_1^* & \cdots & H_d^* \end{bmatrix}$, so that \[(H_\ell)_{ij} = (R \conj{g}_{j - \delta}^i)_\ell = \omega_K^{(i - 1)(j - \delta)}(R g_{j - \delta})_\ell\] and $M_\ell = \sqrt{d} (f_j^d \otimes I_D)^* H = \sum_{k = 1}^d \omega_d^{-(\ell - 1)(k - 1)} H_k,$ giving \begin{align*} (M_\ell)_{ij} &= \sum_{k = 1}^d \omega_d^{-(\ell - 1)(k - 1)} (H_k)_{ij} = \omega_K^{(i - 1)(j - \delta)} \sum_{k = 1}^d \omega_d^{-(\ell - 1)(k - 1)} (Rg_{j - \delta})_k \\ &= \sqrt{d} \omega_K^{(i - 1)(j - \delta)} (\conj{F}_d^* g_{j - \delta})_\ell. \end{align*}  In other words, \begin{equation} M_\ell = \sqrt{d K} \tF_K \diag(f_{2 - \ell}^{d*} g_{1 - \delta},\, \ldots\, , f_{2 - \ell}^{d*} g_{\delta - 1}). \label{eq:block_diag_components} \end{equation}  If $K = 2 \delta - 1$, then $\tF_K$ is unitary, and the singular values of $M_\ell$ are $\{\sqrt{d K} \abs{f_\ell^{d *} g_j}\}_{j = 1 - \delta}^{\delta - 1}$ (since $\abs{f_{2 - \ell}^{d*} g_j} = \abs{\conj{f}_{2 - \ell}^{d*} g_j} = \abs{f_\ell^{d*} g_j}$).  Recognizing that $S^j g_j = g_{-j}$, then \cref{prop:meas_cond} and \eqref{eq:max_sv} take us to \eqref{eq:clean_cond}.

  If $D = 2 \delta - 1 < K$, then the argument remains unchanged, except that the singular values of $M_\ell$, instead of being known explicitly, are bounded above and below by $\max\limits_{|j| < \delta} |f_\ell^{d *} g_j| \sigma_{\max}(\tF_K)$ and $\min\limits_{|j| < \delta} |f_\ell^{d *} g_j| \sigma_{\min}(\tF_K)$ respectively, which gives  \eqref{eq:messy_cond}.
%
\end{proof}

%% file: sections/meas/expl_span_fam.tex
\label{sec:reanalyze}

In this section, we analyze three explicit examples of masks $\gamma \in \Rd$ and their corresponding local Fourier measurement systems, and prove under what conditions these constitute spanning families.  The goal is to constructively provide examples of spanning families that are well-conditioned, and which are scalable in the sense that they may be used for any choice of $d$ and $\delta$.  Specifically, we analyze Example 1 of \cite{IVW2015_FastPhase} (also known as ``exponential masks,'' as we take $\gamma_i = C a^i$ for some $C, a \in \R$) with the new results above, and find improvements on the bounds of its condition number, which scales roughly like $\kappa \approx \delta^2$.  Then, a new set of masks is proposed and studied.  These masks, referred to as ``near-flat masks,'' are constructed by taking $\gamma = a e_1 + \one_{[\delta]} \in \Rd$, and we provide a choice of $a$ that achieves a condition number that is asymptotically linear in $\delta$ -- a notable improvement over the conditioning of the exponential masks.  Finally,  we note the somewhat curious case of a constant mask, $\gamma = \one_{[\delta]}$.  Here, $\gamma$ produces a spanning local Fourier measurement system -- with poor conditioning -- when $d$'s prime divisors are each greater than $\delta$.

\subsection*{Example 1: Exponential Masks}
We first briefly review an example initially proposed in \cite{IVW2015_FastPhase}; in fact, this family of masks is the first spanning family to have been studied in the relevant literature.
Here, we will take $d \in \N$ to be the ambient dimension and $\delta \le \frac{d + 1}{2}$.  Then, we let $\{m_j(a)\}_{j = 1}^{2 \delta - 1}$ be the local Fourier measurement system with mask $\gamma(a) \in \Rd$ defined by $\gamma(a)_i = a^{i - 1},$ for some $0 < a \in \R, a \neq 1$.  The authors in \cite{IVW2015_FastPhase} show that the measurement operator $\Ac$ associated with this family has a condition number bounded by $\kappa \le \max\{144 \ee^2, (\frac{3 \ee (\delta - 1)}{2})^2\}$ when $a$ is chosen to be $\max\{4, \frac{\delta - 1}{2}\}$.  This bound was slightly strengthened in \cite{BP_thesis} through a different choice of $a$, but it was shown to be optimal up to a constant in the sense that, for optimal $a(\delta)$, $\kappa$ still grows as $\bigO(\delta^2)$.

\subsection*{Example 2: Near-Flat Masks}
\label{sec:nearflat_mask}
We now analyze masks of the form $\gamma = a e_1 + \one_{[\delta]}$ in \cref{prop:spike_one}, which, for certain choices of $a(\delta)$, are shown in \cref{prop:spike_one} to produce linear systems with a condition number that is asymptotically less than that of the exponential masks by a factor of $\delta$.

\begin{proposition}
  Let $\{m_j\}_{j \in [D]} \subset \Cd$ be the local Fourier measurement system of support $\delta \le \frac{d + 1}{2}$ with mask $\gamma \in \Rd$ given by $\gamma = a e_1 + \one_{[\delta]}$ where $a > \delta - 1$.  Then this is a spanning family with condition number bounded above by \begin{equation} \kappa \le \dfrac{a^2 + 2 a + \delta}{a - \delta + 1}. \label{eq:spike_cond} \end{equation}  If we choose $a = 2 \delta - 1$, we have $\kappa \le 4 \delta + 1$.
  \label{prop:spike_one}
\end{proposition}

\begin{proof}[Proof of \cref{prop:spike_one}]
  We calculate the condition number directly.  We immediately have $\norm{\gamma}_2^2 = (a + 1)^2 + (\delta - 1) = a^2 + 2a + \delta$, which is the numerator of \eqref{eq:spike_cond}, so it remains only to provide a lower bound on $\sqrt{d} f_j^{d*}(\gamma \circ S^{-m} \gamma)$.  To achieve this, we remark that, for $m \ge 1$, \[\sqrt{d} f_j^{d*} (\gamma \circ S^{-m} \gamma) = a + \sum_{i = 1}^{\delta - m} \omega_d^{(j - 1)(i - 1)} = \begin{piecewise} a + \delta - m & j = 1 \\ a + \dfrac{1 - \omega_d^{(j - 1)(\delta - m)}}{1 - \omega_d^{j - 1}} & \ow \end{piecewise}.\]  Clearly, this expression has its maximum absolute value when $j = 1$, as $\abs{a + \sum_{i = 1}^{\delta - m} \omega_d^{(j - 1) (i - 1)}} \le a + \sum_{i = 1}^{\delta - m} \abs{\omega_d^{(j - 1)(i - 1)}} = a + \delta - m$, so we consider that, for $j \neq 1$, we have \begin{equation}\abs*{\sqrt{d} f_j^{d*} (\gamma \circ S^{-m} \gamma)} \ge \abs*{\Re\left(a + \dfrac{1 - \omega_d^{(j - 1)(\delta - m)}}{1 - \omega_d^{j - 1}}\right)}.\label{eq:flat_ge_real}\end{equation}  We then reduce the term $\Re\left(\frac{1 - \omega_d^{(j - 1)(\delta - m)}}{1 - \omega_d^{j - 1}}\right)$ by setting $\eit := \omega_d^{j - 1}$ and $k := \delta - m$ and finding
  \begin{equation}
    \begin{aligned}
      \Re\left(\dfrac{1 - \ee^{\ii k \theta}}{1 - \eit}\right) &= \Re\left(\dfrac{(1 - \ee^{\ii k \theta})(1 - \ee^{-\theta \ii})}{2 - 2 \cos \theta}\right) \\
      &= \dfrac{(1 - \cos \theta) + \cos(k - 1) \theta - \cos k \theta}{2 - 2 \cos \theta} \\
      &= \frac{1}{2} + \dfrac{\sin(k - \frac{1}{2}) \theta \sin \frac{\theta}{2}}{2 \sin^2 \frac{\theta}{2}},
    \end{aligned}
    \label{eq:flat_real}
  \end{equation}
  where the third line comes from $\sin^2 x = (1 - \cos(2x)) / 2$ and $\cos a - \cos b = \frac{1}{2} \sin\left(\frac{a + b}{2}\right) \sin\left(\frac{b - a}{2}\right)$.  Using $\abs{\sin n\theta} \le n \abs{\sin \theta}$, this gives us that $\Re\left(\frac{1 - \ee^{\ii k \theta}}{1 - \eit}\right) \ge \frac{1}{2} - \frac{2k - 1}{2} = -k$, and hence \begin{equation} \abs*{\sqrt{d} f_j^{d*}(\gamma \circ S^{-m} \gamma)} \ge a - \delta + m \ge a - \delta + 1 \label{eq:nearflat_1} \end{equation} for all $1 \le m < \delta$.  For $m = 0$, a similar calculation gives that \begin{equation} \abs*{\sqrt{d} f_j^{d*}(\gamma \circ \gamma)} \ge a^2 + 2a - \delta, \label{eq:nearflat_2}\end{equation} and we notice that this bound is \emph{greater} than the bound for the case $1 \le m$, stated in \eqref{eq:nearflat_1} whenever $a \ge \frac{1 + \sqrt{5}}{2}$.  This means \eqref{eq:nearflat_1} is tighter than \eqref{eq:nearflat_2} whenever $a > \delta - 1$, which is necessary for these bounds to be positive and meaningful.  Therefore, we may restrict to choices of $a > \delta - 1$ and use the bound $\abs*{\sqrt{d} f_j^{d*}(\gamma \circ S^{-m} \gamma)} \ge a - \delta + 1$.  This completes the proof.
\end{proof}

This result is a new contribution to the collection of well-conditioned local measurement systems.  Example 2 of \cite{our_paper} constituted a local measurement system with a condition number that scales as $\bigO(\delta)$, but it was a somewhat cumbersome construction.  Each of its members $m_j$ was quite sparse, having either 1 or 2 nonzero entries, and did not correspond to a local Fourier measurement system (which corresponds to the diffraction process that we expect to govern our measurement apparatus).    By contrast, the example  in \cref{prop:spike_one}, while not necessarily simple to achieve in the lab, at least has the merit of accomodating the Fresnel diffraction model which motivates the study of local Fourier measurement systems, and its conditioning asymptotically equals that of the sparse construction that was previously the most well-conditioned measurement system known.



\subsection*{Example 3: Constant Masks}
\label{sec:const_mask}
After these two examples, we also remark that the simplest type of mask -- a constant mask, where $\gamma = \one_{[\delta]}$ -- can actually produce a spanning family, albeit a badly conditioned one. Moreover, the conditions required of $\delta$ and $d$ to admit this are upsettingly number theoretical, so we present this result in \cref{prop:constant_gamma} as a negative, though relevant, result.

\begin{proposition}
  Take $d \in \N$ and $\delta \le d$.  Then, with $D = \min(d, 2 \delta - 1)$, the local Fourier measurement system $\{m_j\}_{j = 1}^D$ of support $\delta$ and mask $\gamma = \one_{[\delta]}$ is a spanning family if and only if $d$ is strictly $\delta$-rough, in the sense that $k \mid d \implies k > \delta$.  In this event, and if we additionally take $d > 4$, the condition number of $\Ac$ is bounded by $\kappa \le \delta d^2 / 8$.
  \label{prop:constant_gamma}
\end{proposition}

\begin{proof}[Proof of \cref{prop:constant_gamma}]
  We begin by remarking that $\gamma \circ S^{-(\delta - k)} \gamma = \one_{[k]}$, such that \begin{equation} \sqrt{d} f_j^{d*} (\gamma \circ S^{-(\delta - k)} \gamma) = \sum_{i = 1}^k \omega_d^{(j - 1) (i - 1)} = \begin{piecewise} k & j = 1 \\ \dfrac{1 - \omega_d^{(j - 1) k}}{1 - \omega_d^{(j - 1)}} & \ow \end{piecewise}, \label{eq:const_fj}\end{equation} and hence $f_j^{d*}(\gamma \circ S^{- (\delta - k)} \gamma) = 0$ iff $(j - 1) k = n d$ for some positive integer $n$.  By \cref{cor:gam_fam_span}, this means that $\gamma$ produces a spanning family iff there does not exist a pair $(j, k) \in [d] \times [\delta]_0$ such that $j k = n d$ for some positive integer $n$.  This condition occurs iff there is no pair $(j, k) \in [d] \times [\delta]$ satisfying $j k = d$, which is to say that $\gamma$ produces a spanning family iff $d$ is strictly $\delta$-rough.
  To get the condition number, consider that $\norm{\gamma}_2^2 = \delta$.  It suffices, then, to get a lower bound on $\sqrt{d} \abs*{f_j^{d*} \one_{[k]}}$ for all $k \in [\delta]$.   Trivially following from \eqref{eq:const_fj}, we may write \[\sqrt{d} \abs*{f_j^{d*} \one_{[k]}} \ge \dfrac{\abs{1 - \omega_d}}{2} \ge \dfrac{\abs{\Re(1 - \omega_d)}}{2} = \dfrac{1 - \cos\left(\frac{2 \pi}{d}\right)}{2}.\]  For $d > 4$, we use $1 - \cos(x) \ge (2 x / \pi)^2$ to get that $(1 - \cos(\frac{2 \pi}{d})) / 2 \ge 8 / d^2$, which completes the proof.  The only possible case when $d \le 4$ is $d = 3, \delta = 2$, and we can find by exhaustive calculation that $\kappa = \frac{2}{2 - \sqrt{3}}$.
\end{proof}

This result shows that, while constant masks can produce spanning families in some circumstances, the condition number of the resulting linear system is remarkably unstable as a function of the parameters of the discretization, $d$ and $\delta$.  At the very least, we have that if $(\delta, d)$ admits a constant spanning local Fourier measurement system, then $(\delta, d + 1)$ will not.  Since $d$ is intended to represent the number of pixels in the sensor array, this is a prohibitively specific requirement to be made of the discretization of the phase retrieval problem, so we emphasize that the result of \cref{prop:constant_gamma} is of primarily mathematical interest.

%% file: sections/meas/meas_expl_inv.tex
In this section, we use the results of \cref{sec:con_number} to explicitly state the inverse of the measurement operator $\Ac$ and remark that it can be used to easily deduce the computational complexity of calculating its inverse (initially calculated in \cite{IVW2015_FastPhase}, see also \cite{our_paper}). 
%
\label{sec:inv_runtime}
Fix a local measurement system $\{m_j\}_{j = 1}^D$ with support $\delta$, with associated measurement operator $\Ac$ and canonical matrix representation $A$, as in \cref{eq:meas_op,eq:meas_mat}.  Then, \eqref{eq:interleaved_meas} and \cref{lem:circ_diag} give \[A = P^{(D, d)} (F_d \otimes I_D) \diag(M_\ell)_{\ell = 1}^d (F_d \otimes I_D)^* P^{(d, D)},\]  where we recall $M_\ell$ from \eqref{eq:M_ell}.  In the case where $\{m_j\}_{j = 1}^D$ is a local Fourier measurement system with mask $\gamma$ and modulation index $K$, we define $Z \in \C^{D \times d}$ by \begin{equation} Z_{m \ell} = \sqrt{d K} f_{2 - \ell}^{d*} g_{m - \delta}. \label{eq:four_shift_mat}\end{equation}  Setting $z_\ell = Z e_\ell$, we have \begin{equation} M_\ell = \tF_K D_{z_\ell} \ \text{and} \ \diag(M_\ell)_{\ell = 1}^d = (I_d \otimes \tF_K) \diag(\vec(Z)), \label{eq:M_rearr} \end{equation} by which we further reduce $A$ to \[A = P^{(D, d)} (F_d \otimes I_D) (I_d \otimes \tF_K) \diag(\vec(Z)) (F_d \otimes I_D)^* P^{(d, D)}.\]  This reasoning immediately produces the inverse of $A$, which we state in \cref{prop:A_inv_gam}. 
\begin{proposition} \label{prop:A_inv_gam}
  Let $A \in \C^{d D \times d D}$ be the canonical representation of the measurement operator $\Ac$ associated with a local Fourier measurement system $\{m_j\}_{j = 1}^d$ of support $\delta \le \frac{d + 1}{2}$ with mask $\gamma \in \R^d$.  Defining $Z$ as in \eqref{eq:four_shift_mat}, we have
  \begin{equation} A^{-1} = P^{(D, d)} (F_d \otimes I_D) (\diag(\vec(Z)))^{-1} (I_d \otimes \tF_K^*) (F_d \otimes I_D)^* P^{(d, D)}. \label{eq:A_inv_gam} \end{equation} If $\{m_j\}_{j = 1}^D$ is a general local measurement system, and $\Ac$ is invertible, then its inverse is given by \[A^{-1} = P^{(D, d)} (F_d \otimes I_D) \diag(M_\ell^{-1})_{\ell = 1}^d (F_d \otimes I_D)^* P^{(d, D)}.\]
\end{proposition}
\noindent This formulation makes it straightforward to deduce the computational complexity of inverting $\Ac$. In the case of local Fourier measurement systems, the dominant cost is that of computing $D=\bigO(\delta)$ Fourier transforms of size $d$, so  the cost of inverting $A$ comes out to $\bigO(\delta\, d \log d)$, as in 
 \cite{IVW2015_FastPhase}.


\if{\subsection{Preliminaries in Probability}
Prior to discussing the variance of the images of vectors under $A^{-1}$ or $\Ac^{-1}$, we review some preliminaries regarding the real and complex multivariate Gaussian distributions.  These results and the notation with which we express them may be found in many standard texts, for example \cite{tong1990normal} for the real case and section 7.8.1 of \cite{gallager2008digital} for the complex.  For $\mu \in \Rn$ and $0 \prec \Sigma \in \sym^n, \Nc(\mu, \Sigma)$ refers to the multivariate normal distribution on $\Rn$ with mean $\mu$ and covariance matrix $\Sigma = \mathbb{E}_x[(x - \mu) (x - \mu)^T]$, and is determined by its probability density function \[f_{\Nc}(x; \mu, \Sigma) = \frac{1}{(2 \pi)^{n / 2} \det(\Sigma)^{1 / 2}} \exp\left(- \dfrac{(x - \mu)^T \Sigma^{-1} (x - \mu)}{2}\right).\]  For $\nu \in \Cn$ and $0 \prec \Xi \in \H^n, \CN(\nu, \Xi)$ is the (circularly symmetric about $\nu$) multivariate complex normal distribution on $\Cn$ with mean $\nu$, covariance $\Xi = \mathbb{E}_z[(z - \nu) (z - \nu)^*]$, and density function \[f_{\CN}(z; \nu, \Xi) = \dfrac{1}{\pi^n \det(\Xi)} \exp(- (z - \nu)^* \Xi^{-1} (z - \nu))\]  We remark that circular symmetry is defined by having that the real and imaginary parts of $z \sim \CN(0, \Xi)$ be i.i.d., and is ensured by tacitly requiring, as we shall throughout this dissertation, that $\mathbb{E}[(z - \nu) (z - \nu)^T] = 0$ (see Theorem 7.8.1 of \cite{gallager2008digital}).

We now relate a standard result from the literature concerning linear transformations of Gaussian random vectors.
\begin{proposition}[Theorem 3.3.3 in \cite{tong1990normal} and Section 7.8.1 of \cite{gallager2008digital}] \label{prop:gausaff}
  Suppose $x \sim \Nc(\mu, \Sigma)$.  Then for $A \in \Rmxn, b \in \Rm,$ \begin{equation} Ax + b \sim \Nc(A \mu + b, A \Sigma A^T). \label{eq:real_gausaff} \end{equation}

  Suppose $z \sim \CN(\nu, \Xi)$.  Then for $B \in \Cmxn, c \in \Cm,$ \begin{equation} B z + c \sim \CN(B \nu + c, B \Xi B^*). \label{eq:cpx_gausaff}\end{equation}
\end{proposition}
We remark that the result regarding complex Gaussian vectors implies that linear transformations preserve the property of circular symmetry.

\subsection{Distribution of variance}
\label{sec:dist_var}
In the interest of studying the propagation of error through our recovery algorithm, in this section, we describe the probability distribution of the noise on each entry of $\Ac^{-1}(\Ac(x x^*) + \eta)$ as a function of the noise vector's distribution.  To keep things tractable, we assume that the entries of $\eta$ are identically and independently distributed; specifically, we will assume that $\eta_{(\ell, j)} \iid \Nc(0, \sigma^2)$ for some $\sigma \ge 0$.

Before we begin, we remark that the distribution of $\Ac^{-1}(\eta)_{i, i + m}$ will depend only on $m$.  This follows intuitively by noting that $\Ac$ commutes nicely with ``diagonal shifts'' in the sense that \[\Ac(S^k X S^{-k})_{(\ell, j)} = \inner{S^{\ell} m_j m_j^* S^{-\ell}, S^k X S^{-k}} = \inner{S^{\ell - k} m_j m_j^* S^{-(\ell - k)}, X} = \Ac(X)_{(\ell - k, j)},\] so that if $y_1, y_2 \in \R^{[d] \times [D]}$ satisfy $(y_1)_{\ell, j} = (y_2)_{\ell - k, j}$ for some $k \in \N$, we will have $\Ac^{-1}(y_1) = S^k \Ac^{-1}(y_2) S^{-k}$.  In particular, if the entries of $y_1$ are identically, independently distributed random variables and $y_2$ is \emph{defined} by $(y_2)_{\ell - k, j} = (y_1)_{\ell, j}$, then $y_1$ and $y_2$ are drawn from the same distribution on $\R^{[d] \times [D]}$.  This means that $\Ac^{-1}(y_1)$ and $\Ac^{-1}(y_2)$ are identically distributed, but $\Ac^{-1}(y_1) = S^k \Ac^{-1}(y_2) S^{-k}$, so the distributions of $\Ac^{-1}(y_1)$ and $S^k \Ac^{-1}(y_1) S^{-k}$ are identical.  In particular, the distribution of the image of i.i.d.~noise under $\Ac^{-1}$ is invariant under such diagonal shifts, so $\Ac^{-1}(\eta)_{i, i + m}$ is distributed identically to (though not necessarily independently from!) $\Ac^{-1}(\eta)_{1, 1 + m}$, and this conclusion holds for all $i$ and $m$.

To make this precise, and to discover the distribution of $\Ac^{-1}(\eta)_{1, 1 + m}$ exactly, we will present another means by which $\Ac^{-1}(y)$ may be calculated from $y$.  With this in mind, we remark that \eqref{eq:interkron_mat} of \cref{lem:interkron} gives that \[P^{(D, d)} (F_d \kron I_D) = \rowmat{I_D \kron f^d}{1}{d},\] so $A$ may be transformed by \eqref{eq:kron_diag} to give \[A = \rowmatfun{M_@ \kron f_@^d}{1}{d} \colmat{I_D \kron f^d}{1}{d} = \sum_{j = 1}^d M_j \kron f_j^d f_j^{d*}.\]  Setting $X = \rowmat{\chi}{1 - \delta}{\delta - 1} \in \C^{d \times 2 \delta - 1}$, \cref{lem:kronvec} gives us \[A \colmat{\chi}{1 - \delta}{\delta - 1} = \vec\left(\sum_{j = 1}^d f_j^d f_j^{d*} X M_j^T\right).\]  Recalling \eqref{eq:M_rearr} along with $\tF_D^{-T} = (\tF_D^T)^* = \conj{\tF}_D$, we have \[f_j^{d*} \mat_{(d, D)}\left(A \colmat{\chi}{1 - \delta}{\delta - 1}\right) \conj{\tF}_D = f_j^{d*} X D_{z_j},\]
 so that, for $\ell \in [2 \delta - 1]$, and recalling $\conj{\tF}_D e_\ell = \conj{f}^D_{\ell + 1 - \delta} = f_{\delta + 1 - \ell}^D$, we have \[f_j^{d *} \mat_{(d, D)}(A \vec(X)) f_{\delta + 1 - \ell}^D = f_j^{d*} \mat_{(d, D)}(A \vec(X)) \conj{\tF}_D e_\ell = f_j^{d*} \chi_{\ell - \delta} Z_{\ell j}.\] In this way, from $A \vec(X)$ we may recover \[b_{\ell} := F_d^* \mat_{(d, D)}(A \vec(X)) f_{\delta + 1 - \ell}^D = \vec(f_j^{d*} Z_{\ell j} \chi_{\ell - \delta})_{j = 1}^d = D_{Z^T e_{\ell}} F_d^* \chi_{\ell - \delta}\] for each $\ell$, from which $\chi_{\ell - \delta}$ is determined by taking $\chi_{\ell - \delta} = F_d D_{Z^T e_{\ell}}^{-1} b_\ell$.
In other words, for $y \in \R^{d D}$, the $m\th$ diagonal of $\Dc_\delta^{-1}(A^{-1} y)$, for $m = 1 - \delta, \ldots, \delta - 1$, is given by \begin{equation} \chi_m = F_d D^{-1}_{Z^T e_{m + \delta}} F_d^* \mat_{(d, D)}(y) f_{1 - m}^D, \label{eq:chi_m} \end{equation} and we use this expression to deduce the distribution of noise on the $m\th$ diagonal when $y$ is a random variable.  For instance, if we consider that \[\Dc_\delta^{-1} A^{-1}(A \Dc_\delta(T_\delta(x x^*)) + \eta) = T_\delta(x x^*) + \Dc_\delta^{-1} (A^{-1} \eta),\] where $\eta_j \iid \Nc(0, \sigma^2)$ for $j \in [d (2 \delta - 1)]$, knowing the distribution of $A^{-1} \eta$ will tell us the distribution of noise on our recovered estimate of $T_\delta(x x^*)$.  We deduce this result directly from \eqref{eq:chi_m}, and summarize it in \cref{prop:var_dist}.  We remark that, in the statement and proof of this proposition, exponentiation of vectors is entry-wise.
\begin{proposition}
  Suppose that $\{m_j\}_{j = 1}^{2 \delta - 1}$ is a local Fourier measurement system with support $\delta \le \frac{d + 1}{2}$, mask $\gamma$, and modulation index $K = D = 2 \delta - 1$, with associated measurement operator and representation matrix $\Ac$ and $A$.  Suppose further that $\eta \in \R^{d (2 \delta - 1)}$ has entries that are i.i.d.~Gaussian random variables, namely $\eta_j \iid \Nc(0, \sigma^2)$ for $j \in [d (2 \delta - 1)]$ and some $\sigma \ge 0$.  Then, setting $N = \Dc_\delta^{-1} A^{-1} \eta \in T_\delta(\Cdxd)$, then $N$ is Hermitian, its $0, \ldots, \delta - 1\st$ diagonals are distributed  independently from one another, and the $m\th$ diagonal of $N$ is distributed by \begin{align} \diag(N, m) &\sim \CN(0, \sigma^2 \circop(s_m)), m \in [\delta - 1], \ \text{and} \nonumber \\ \diag(N, 0) &\sim \Nc(0, \sigma^2 \circop(s_0)), \ \text{where} \nonumber \\ s_m &= \frac{1}{D d^{3 / 2}} F_d^* \abs*{F_d^* g_m}^{-2}.\label{eq:s_m} \end{align}
  \label{prop:var_dist}
\end{proposition}

\begin{proof}[Proof of \cref{prop:var_dist}]
  To establish that $N$ is Hermitian, we remark that $\Cdxd = \H^d \oplus \Skew(d)$, where \[\Skew(d) = \{B \in \Cdxd : B^* = -B\}\] is the set of skew-Hermitian matrices and where $\oplus$ represents the direct product.  In other words, given any $M \in \Cdxd$, we have a unique $H \in \H^d, B \in \Skew(d)$ such that $M = H + B$.  We now consider that, if $H \in \H^d$ and $B \in \Skew(d)$, we have that
  \begin{align*}
    \Tr(H^* B) &= \Tr(H B) = \Tr(B H) \\
    &= -\Tr(B^* H) = - \conj{\Tr(H^* B)},
  \end{align*}
  meaning that $\Tr(H^* B) \in \ii \R$.  Additionally, we remark that the Hermitian matrices are a real Hilbert space, so for another $H' \in \H^d$, we have $\Re \inner{H', M} = \inner{H', H}$ and $\Im \inner{H', M} = \inner{H', B} / \ii$.  Therefore, since all the measurement matrices $S^\ell m_j m_j^* S^{- \ell}$ appearing in $\Ac$ are Hermitian, given $M \in \Cdxd$ decomposed into its Hermitian and skew-Hermitian parts $M = H + B$, we have that $\Re \Ac(M) = \Ac(H)$ and $\Im \Ac(M) = \Ac(B) / \ii$.  In particular, $\Ac^{-1}(\R^{[d]_0 \times [D]}) \subset \H^d$, and $N$, being the inverse image of a real vector $\eta$, will be Hermitian.

  For convenience, throughout the remainder of the proof we will set $\chi_m = \diag(N, m)$ for $m = 1 - \delta, \ldots, \delta - 1$.  To prove the independence of $\{\chi_m\}_{m \in [\delta]_0}$, we consider that \eqref{eq:chi_m} tells us that \[\chi_m = F_d D^{-1}_{Z^T e_{m + \delta}} F_d^* \left(\mat_{(d, D)}(\eta) f^D_{1 - m}\right),\] and we focus on the term $\mat_{(d, D)}(\eta) f_{1 - m}^D$.  Considering that $\eta \sim \Nc(0, \sigma^2 I_{d D})$, we have that the rows $\rowmat{r}{1}{d}^T$ of $\mat_{(d, D)}(\eta)$ are distributed according to $r_i = \mat_{(d, D)}(\eta)^T e_i \sim \Nc(0, \sigma^2 I_D)$.  At this point, it would be convenient if we could merely cite \cref{prop:gausaff} to establish the distribution of $r_i^T f^D_{1 -m }$ or indeed $\mat_{(d, D)}(\eta) f_{1 - m}^D$, but we remark that we don't have a result for the image of a real Gaussian vector under a complex linear transformation.  Therefore, we consider the real and imaginary parts of $\mat_{(d, D)}(\eta) f_{1 - m}^D$ separately, setting $v_m = \mat_{(d, D)}(\eta) f_{1 - m}^D$ and seeing that \[(v_m)_i = r_i^T \Re(f_{1 - m}^D) + \ii\, r_i^T \Im(f_{1 - m}^D) = \Re(f_{1 - m}^D)^T r_i + \ii\, \Im(f_{1 - m}^D)^T r_i.\]  Since $\{f_1^D\} \cup \{\sqrt{2} \Re(f_{1 - m}^D)\}_{m \in [\delta - 1]} \cup \{\sqrt{2} \Im(f_{1 - m}^D)\}_{m \in [\delta - 1]}$ is an orthonormal basis for $\R^D$, then compiling these vectors into a matrix $Q$, we have from \cref{prop:gausaff} that $Q^T r_i \sim \Nc(0, \sigma^2 Q^T Q) = \Nc(0, \sigma^2 I_D)$, meaning the real and imaginary components of $v_1, \ldots, v_{\delta - 1}$ and $v_0$ are all independent, with $v_0 \sim \Nc(0, \sigma^2 I_d), \Re(v_m) \sim \Nc(0, \frac{\sigma^2}{2} I_d),$ and $\Im(v_m) \sim \Nc(0, \frac{\sigma^2}{2} I_d)$ for $m \in [\delta - 1]$.  Therefore, $v_m \iid \CN(0, \sigma^2 I_d)$ and $v_0 \sim \Nc(0, \sigma^2 I_d)$, independently from the other $v_m$.  Since the $v_m$ are independent, clearly their images under non-singular, fixed (independently of the random process) linear transformations will also be independent.  In particular, the diagonals $\chi_m = F_d D^{-1}_{Z^T e_{m + \delta}} F_d^* v_m$ will be independent of one another for $m \in [\delta]_0$.

  To get the actual distribution of the $\chi_m$ for $m \in [\delta - 1]$, we simply quote \cref{prop:gausaff} again to get that \[\chi_m \sim \CN(0, \sigma^2 F_d D^{-1}_{Z^T e_{m + \delta}} F_d^* I_d F_d D^{-1*}_{Z^T e_{m + \delta}} F_d^*) = \CN(0, \sigma^2 F_d D^{-2}_{\abs{Z^T e_{m + \delta}}} F_d^*).\] The covariance matrix becomes, by recalling \eqref{eq:circ_dft_diag} and the definition of $Z$ in \eqref{eq:four_shift_mat}, \[\sigma^2 F_d D^{-2}_{\abs{Z^T e_{m + \delta}}} F_d^* = \circop\left(d^{-1/2} F_d^* \abs{Z^T e_{m + \delta}}^{-2}\right) = \circop(s_m).\]  The only distinction for $\chi_0$ is that the distributions are all $\Nc$ instead of $\CN$; the same calculations as above give $\chi_0 \sim \Nc(0, \sigma^2 \circop(s_0))$.  To verify that $s_0 \in \Rd$, we consider that $g_0 = \gamma \circ \gamma$, so $\tilde{g} := F_d^* g_0 = \frac{1}{\sqrt{d}} (F_d^* \gamma) * (F_d^* \gamma) $ satisfies $\tilde{g}_i = \tilde{g}_{2 - i}$ (equivalently, $R \tilde{g} = \tilde{g}$).  Similarly, $R(\tilde{g}^{-2}) = \tilde{g}^{-2}$, which guarantees that $s_0 = \frac{1}{D \sqrt{d}} F_d^* \tilde{g}^{-2} \in \R^d$.
\end{proof}
}\fi
\if{\subsection{Slight Improvement to Magnitude Estimation Bounds}
{\color{red} How best to incorporate this? Or should we get rid of it?}
One useful corollary of the calculations in \cref{sec:dist_var} allows us to sharpen our bounds on the error of the magnitude estimation step of \cref{alg:phaseRetrieval1} studied in \cref{sec:recov_guar}.  In particular, because \eqref{eq:chi_m} shows us how to calculate a specific diagonal $\chi_m$ from the measurement data $y$, we can get a bound on $\norm{\diag(X - X_0)}_1$ that is sharper than the $\norm{X - X_0}_F$ that we used in the proof of \cref{lem:diag_mag_diff}.  We briefly state and prove this improvement here, and cite it later in the finalized, sharpest bounds presented in \cref{sec:ang_sync_improve,sec:2d_recov}.

To describe precisely what we mean, we define the diagonal projection operators $P_{\diag(\Cmxn, \ell)} : \Cmxn \to \Cm$ by $P_{\diag(\Cmxn, \ell)}(A) = \diag(A, \ell)$.  With this in hand, we consider that \eqref{eq:chi_m} may be rewritten as \begin{equation} (P_{\diag(\Cdxd, m)} \circ A^{-1}) y = F_d D^{-1}_{Z^T e_{m + \delta}} F_d^* \mat_{(d, D)}(y) f_{1 - m}^D, \label{eq:chi_mop} \end{equation} where all terms are defined as in \cref{sec:dist_var}.  We may rewrite this, using \cref{lem:kronvec}, as \[(P_{\diag(\Cdxd, m)} \circ A^{-1}) y = (f_{1 - m}^{D T} \kron F_d D^{-1}_{Z^T e_{m + \delta}} F_d^*) y.\]  The singular values of $f_{1 - m}^{D T} \kron F_d D^{-1}_{Z^T e_{m + \delta}} F_d^*$ are the singular values of $F_d D^{-1}_{Z^T e_{m + \delta}} F_d^*$, which in turn are simply the entries of $Z^T e_{m + \delta}$.  When we consider the main diagonal $\chi_0$, we notice from \eqref{eq:four_shift_mat} that this gives $\sigma = \sqrt{d K} \abs{F_d^* (\gamma \circ \gamma)}$.  We state this result, along with specific bounds on $\sigma_{\min}$ and $\kappa$ for $\gamma_{\exp}$ and $\gamma_{\rm flat}$ studied in \cref{sec:expl_span_fam}, in \cref{prop:mag_est_imp}.

\begin{proposition} \label{prop:mag_est_imp}
  Consider a local Fourier measurement system with support $\delta$ with $D := 2 \delta - 1 \le d$, mask $\gamma$, and matrix representation $A \in \C^{d D \times d D}$.  Set $Y = (P_{\diag(\Cdxd, 0)} \circ A^{-1})$.  Then 
  \begin{equation}
    \begin{gathered}
      \sigma_{\max}(Y) = \sqrt{D} \norm{\gamma}_2^2 \\
      \sigma_{\min}(Y) = \sqrt{d D} \min_{j \in [d]} \abs{f_j^{d*} (\gamma \circ \gamma)}.
    \end{gathered}
    \label{eq:magimp}
  \end{equation}
  When $\gamma = \gamma_{\exp}$ or $\gamma_{\rm flat}$ with the recommended parameter $a_\delta$ from \cref{prop:exp_kappa,prop:nearflat_more}, then we have, for $\delta \ge 5$,
  \begin{equation}
    \begin{gathered}
      \sigma_{\min}(Y_{\exp}) \ge 20 \sqrt{D}, \ \kappa(Y_{\exp}) \le \delta / 2 \\
      \sigma_{\min}(Y_{\rm flat})  \ge \frac{\sqrt{D}}{6} (\delta - 1)^2, \ \kappa(Y_{\rm flat})  \le 1 + \frac{18}{\delta - 1}
    \end{gathered}
    \label{eq:mask_magimp}
  \end{equation}
\end{proposition}

\begin{proof}[Proof of \cref{prop:mag_est_imp}]
  The main statement, \eqref{eq:magimp}, is proven in the remarks preceding the proposition.  The inequalities of \eqref{eq:mask_magimp} are proven by straightforward calculation.

  To get $\sigma_{\min}(Y_{\exp}),$ we write that, for any $\omega \in \Sbb^1$,
  \begin{align*}
    \frac{1}{\sqrt{D}} \sigma_{\min}(Y_{\exp}) &\ge \abs*{\sum_{i = 1}^\delta a_\delta^{2(i - 1)} \omega^{i - 1}} = \abs*{\dfrac{1 - a_\delta^{2 \delta} \omega^\delta}{1 - a_\delta^2 \omega}} \\
    &\ge \dfrac{a_\delta^{2 \delta} - 1}{a_\delta^2 + 1} = \dfrac{\left(1 + \frac{4}{\delta - 2}\right)^\delta - 1}{2 + \frac{4}{\delta - 2}},
  \end{align*}
  which may be shown numerically to exceed $20$.  For $\kappa(Y_{\exp})$, we use $\frac{1}{\sqrt{D}} \sigma_{\min}(Y_{\exp}) \ge \frac{a_\delta^{2 \delta} - 1}{a_\delta^2 + 1}$ and $\frac{1}{\sqrt{D}} \sigma_{\max}(Y_{\exp}) = \frac{a_\delta^{2 \delta} - 1}{a_\delta^2 - 1}$ to get \[\kappa(Y_{\exp}) \le \dfrac{a_\delta^2 + 1}{a_\delta^2 - 1} = \dfrac{2 + \frac{4}{\delta - 2}}{\frac{4}{\delta - 2}} = \delta / 2.\]

  To get $\sigma_{\min}(Y_{\rm flat})$, we remark that 
  \begin{align*}
    \frac{1}{\sqrt{D}} \sigma_{\min}(Y_{\rm flat}) &\ge (a_\delta + 1)^2 - (\delta - 1) = 4 \abs{C_0}^2 \delta^2 + (4 \abs{C_0} - 1) \delta + 2 \\
    &\ge \frac{1}{6} \delta^2 - \frac{1}{6} \delta + 2 \ge \frac{1}{6}(\delta - 1)^2,
  \end{align*}
  where we have merely used numerical bounds on the constant $\abs{C_0} \approx 0.434467$ discussed in the proof of \cref{prop:nearflat_more}.  For $\kappa(Y_{\rm flat})$, we use that $\frac{1}{\sqrt{D}} \sigma_{\max}(Y_{\rm flat}) = 4 C_0^2 \delta^2 + (4 C_0 + 1) \delta + 2 = \frac{1}{\sqrt{D}} \sigma_{\min}(Y_{\rm flat}) + 2 (\delta - 1)$, such that \[\kappa(Y_{\rm flat}) = 1 + \dfrac{2 \delta - 1}{\frac{1}{\sqrt{D}} \sigma_{\min}(Y_{\rm flat})} \le 1 + \dfrac{18}{\delta - 1},\] where we have used that $2 \delta - 1 \le 3 \delta - 3$.  
\end{proof}
}\fi

%% file: sections/ptychography/ptychography_sec.tex
\input{ptych_intro}
\label{sec:ptych_intro}

\subsection{Conditioning of $\Ac$ for Ptychography}
\input{con_number_ptych}
\label{sec:con_number_ptych}

\section{Analysis of the Recovery Algorithm}
\input{ptych_recovery}
\label{sec:ptych_recov}
\label{sec:ptych_recovery}

\input{ptych_num}
\label{sec:ptych_num}

%% file: sections/ptychography/ptych_intro.tex
In our model for the ptychographic setup of \eqref{eq:shift_model}, we have so far assumed that measurements are taken corresponding to all shifts $\ell \in [d]_0$; in the notation of \eqref{eq:shift_model}, this is equivalent to taking $P = [d]_0$.  Herein we present a useful generalization to the case where $P = s[d / s]_0$, and $s \in \N$ divides $d$.  

The motivation for studying this case is that, unfortunately, in practice taking $P = [d]_0$ is usually an impossibility, since in many cases an illumination of the sample can cause damage to the sample \cite{starodub2008damage}, and applying the illumination beam (which can be highly irradiative) repeatedly at a single point can destroy it.  In ptychography as it is usually performed in the lab, the beam is shifted by a far larger distance than the width of a single pixel\footnote{The difficulty of moving the illumination apparatus at a scale equal to the desired optical resolution is another reason taking $P = [d]_0$ is a cumbersome assumption.} -- instead of overlapping on $\delta - 1$ of $\delta$ pixels, adjacent illumination regions will typically overlap on a percentage of their support, on the order of $50\%$ (or even less sometimes) \cite{marchesini2015coptych,shapiro2014nanometer}.  Considering the risks to the sample and the costs of operating the measurement equipment, there are strong incentives to reduce the number of illuminations applied to an object, so our theory ought to address a model that reflects this.

\subsection{Measurement Operator and Its Domain}
\label{sec:pty_meas_op}
\begin{figure}
  \centering
  \begin{subfigure}[b]{0.4\textwidth}
    \begin{tikzpicture}[ampersand replacement=\&,baseline=-\the\dimexpr\fontdimen22\textfont2\relax]
    \matrix (m)[matrix of math nodes,left delimiter=(,right delimiter=)]
            {
              * \& * \& * \&   \&   \&   \& * \& * \\
              * \& * \& * \& * \&   \&   \&   \& * \\
              * \& * \& * \& * \& * \&   \&   \&   \\
                \& * \& * \& * \& * \& * \&   \&   \\
                \&   \& * \& * \& * \& * \& * \&   \\
                \&   \&   \& * \& * \& * \& * \& * \\
              * \&   \&   \&   \& * \& * \& * \& * \\
              * \& * \&   \&   \&   \& * \& * \& * \\
            };

            \begin{pgfonlayer}{myback}
              \fhighlight[blue!30]{m-1-1}{m-3-3}
              \fhighlight[blue!30]{m-2-2}{m-4-4}
              \fhighlight[blue!30]{m-3-3}{m-5-5}
              \fhighlight[blue!30]{m-4-4}{m-6-6}
              \fhighlight[blue!30]{m-5-5}{m-7-7}
              \fhighlight[blue!30]{m-6-6}{m-8-8}
              \fhighlight[blue!30]{m-7-1}{m-8-1}
              \fhighlight[blue!30]{m-1-7}{m-1-8}
              \fhighlight[blue!30]{m-8-8}{m-8-8}
              \fhighlight[blue!30]{m-8-1}{m-8-2}
              \fhighlight[blue!30]{m-1-8}{m-2-8}
              \fhighlight[blue!30]{m-1-1}{m-2-2}
            \end{pgfonlayer}
  \end{tikzpicture}
    \caption{$T_3(\C^{8 \times 8})$}
  \end{subfigure}
  \begin{subfigure}[b]{0.4\textwidth}
    \begin{tikzpicture}[ampersand replacement=\&,baseline=-\the\dimexpr\fontdimen22\textfont2\relax]
    \matrix (m)[matrix of math nodes,left delimiter=(,right delimiter=)]
            {
              * \& * \& * \&   \&   \&   \& * \& * \\
              * \& * \& * \& * \&   \&   \&   \& * \\
              * \& * \& * \& * \& * \&   \&   \&   \\
                \& * \& * \& * \& * \& * \&   \&   \\
                \&   \& * \& * \& * \& * \& * \&   \\
                \&   \&   \& * \& * \& * \& * \& * \\
              * \&   \&   \&   \& * \& * \& * \& * \\
              * \& * \&   \&   \&   \& * \& * \& * \\
            };

            \begin{pgfonlayer}{myback}
              \fhighlight[blue!30]{m-1-1}{m-3-3}
              \fhighlight[blue!30]{m-3-3}{m-5-5}
              \fhighlight[blue!30]{m-5-5}{m-7-7}
              \fhighlight[blue!30]{m-7-1}{m-8-1}
              \fhighlight[blue!30]{m-1-7}{m-1-8}
              \fhighlight[blue!30]{m-7-7}{m-8-8}
            \end{pgfonlayer}
    \end{tikzpicture}
    \caption{$T_{3, 2}(\C^{8 \times 8})$}
  \end{subfigure}
  \caption{The support of the subspaces $T_\delta(\C^{d \times d})$ vs. $T_{\delta, s}(\C^{d \times d})$ for $d = 8, \delta = 3, s = 2$}
  \label{fig:T_delta_s}  
\end{figure}
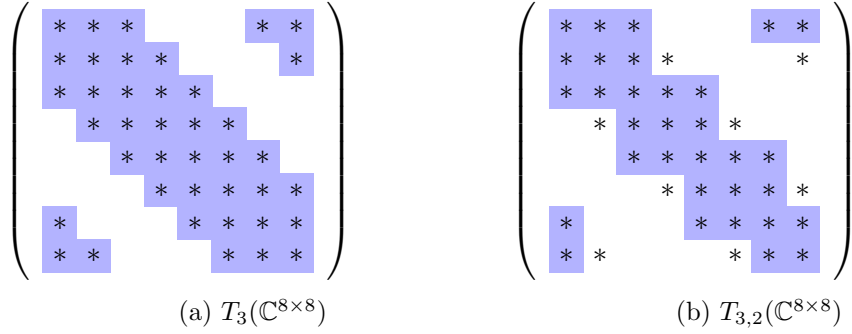

Towards this model, instead of using all shifts in our lifted measurement system, we fix a shift size $s \in \N$, $s<\delta$, where $d = \dbar s$ with $\dbar \in \N$ and use $S^{s \ell} m_j m_j^* S^{-s \ell}$ for $\ell \in [\dbar]_0$.  We introduce the following generalization of the lifted measurement system: given a family of masks of support $\delta$, $\{m_j\}_{j \in [D]} \subset \C^d$, and $s, \dbar \in \N$ with $\dbar = d / s$, the associated \emph{lifted measurement system of shift $s$} is \begin{equation} \Lc_{\{m_j\}}^s := \{S^{s \ell} m_j m_j^* S^{-s \ell}\}_{(\ell, j) \in [\dbar]_0 \times [D]} \subset \C^{d \times d}. \label{eq:lift_pty_sys}\end{equation}
This leads to an obvious redefinition of the measurement operator, now $\Ac : \Cdxd \to \R^{[\dbar]_0 \times [D]}$: \begin{equation} \mathcal{A}(X)_{(\ell, j)} = \inner{S^{s \ell} m_j m_j^* S^{-s \ell}, X}, \quad (\ell, j) \in [\dbar]_0 \times [D]. \label{eq:pty_meas_op} \end{equation}  This will also force us to reconsider the subspace of $\Cdxd$ with which we are working in the domain of $\Ac$, since it is clearly impossible, by inspection of \cref{fig:T_delta_s}, for $\Lc^s$ to span $T_\delta(\Cdxd)$ with a shift size $s > 1$.  In an effort to define the subspace analogous to $T_\delta(\Cdxd)$ in the ptychographic case, we let $\Jc_{\delta, s} = \bigcup_{\ell \in [\dbar]_0}\supp(S^{s \ell} \one_{[\delta]} \one_{[\delta]}^* S^{-s \ell})$ be the set of indices ``reached'' by this system, and we let \begin{equation} T_{\delta, s}(X) = \left\{\begin{array}{r@{,\quad}l} X_{ij} & (i, j) \in \Jc_{\delta, s} \\ 0 & \text{otherwise}\end{array}\right.\label{eq:T_delta_s}\end{equation} be the projection onto the associated subspace of $\Cdxd$.  $T_{\delta, s}$ is visualized in \cref{fig:T_delta_s}.

%% file: sections/ptychography/con_number_ptych.tex
With this setup in hand, we begin our analysis of the linear system $\Ac(X) = y$ with a number of lemmas that unravel the structure of this operator.  Our goal will be to proceed similarly to \cref{sec:con_number} by rewriting $\Ac$ as a product of a block-circulant matrix with certain permutations, at which point we will be able to cite \cref{cor:circ_diag_condition}, which renders a convenient expression for the condition number.
%
\label{sec:pty_new_op}
In service of this strategy, in this section we introduce a few new operators that are useful in the analysis of $\Ac$.  For $N \in \N$, we define $\Tc_N : \bigcup_{\ell \in \N} \C^{\ell N \times m} \to \bigcup_{\ell \in \N} \C^{\ell m \times N},$ the blockwise transpose operator, defined by \begin{equation} \Tc_N\left(\begin{bmatrix} V_1 \\ \vdots \\ V_\ell \end{bmatrix}\right) = \begin{bmatrix} V_1^* \\ \vdots \\ V_{\ell}^* \end{bmatrix}\label{eq:blk_trans_op} \end{equation} for $V_1, \ldots, V_\ell \in \C^{N \times m}$.
We also define, for $(k_j)_{j = 1}^n$ and permutation $P \in \{0, 1\}^{n \times n}$, the \emph{blockwise permutation operator} $\Pc(P, (k_j)) : \C^{K \times K} \to \C^{K \times K}$, where $K = \sum_{j = 1}^n k_j$.  Our intention will be to permute the blocks of a block vector $\rowmat{v^T}{1}{n}^T$, where $v_j \in \C^{k_j}$.  In order to specify $\Pc(P, (k_j))$ precisely, we permit an overloading of notation on permutations: namely, if $P \in \{0, 1\}^{m \times m}$ is a permutation, then we identify $P$ with the mapping $\pi : [m] \to [m]$ where $\pi(i) = j$ whenever $P e_i = e_j$.  In particular, if we write $P(i)$, we mean ``$j$ such that $P e_i = e_j$.''  With this in mind, $\Pc(P, (k_j))$ is defined, for $v_j \in \C^{k_j}$, by \begin{equation} \Pc(P, (k_j)) \colmatfl{v_1}{v_n} = \colmatfl{v_{P(1)}}{v_{P(n)}}. \label{eq:block_perm} \end{equation}  

\if{
\subsection{Lemmas on Block Circulant Structure}
\label{sec:ptych_lem}
We begin with \cref{lem:circ_transpose}, which describes the transposes of block circulant matrices.  For this lemma and the remainder of this section, the reader is advised to recall the definitions of $R_k$ and $P^{(d, N)}$ from \cref{sec:notation} and \eqref{eq:interleave_def}, as well as $\Pc(P, \{k_i\})$ and $\Tc_N$ from \cref{sec:pty_new_op}.
\begin{lemma}
\label{lem:circ_transpose} Given $k, N, m \in \N$ and $V \in \C^{kN \times m}$, we have \[\circop^N(V)^* = \circop^m\left( (R_k \otimes I_m) \Tc_N(V) \right).\] 
\end{lemma}
\begin{proof}[Proof of \cref{lem:circ_transpose}]
  Suppose $V_i$ are the $N \times m$ blocks of $V$, such that $V = \left[V_1^T \cdots V_k^T\right]^T$.  Indexing blockwise, we have $\circop^N(V)_{[ij]} = V_{i - j + 1}$, so that $\circop^N(V)^*_{[ij]} = V_{j - i + 1}^*$.  In other words,
  \[
  \circop^N(V)^* = \begin{bmatrix} V_1^* & V_2^* & \cdots & V_N^* \\ V_N^* & V_1^* & \cdots & V_{N - 1}^* \\ \vdots & & \ddots & \vdots \\ V_2^* & V_3^* & \cdots & V_1^* \end{bmatrix} = \circop^m((R_k \otimes I_m) \Tc_N(V) )
 . \]
\end{proof}

\Cref{lem:block_circ_right,lem:diag_kron_perm} provide identities for a few block matrix structures that will be of interest.
\begin{lemma}\label{lem:block_circ_right}
  Given $N_1, N_2, k, m \in \N$ and $V_i \in \C^{k N_1 \times m}$ for $i \in [N_2]$, we have \[\begin{bmatrix} \circop^{N_1}(V_1) & \cdots & \circop^{N_1}(V_{N_2}) \end{bmatrix} (P^{(k, N_2)} \otimes I_m)^* = \circop^{N_1}\left(\begin{bmatrix} V_1 & \cdots & V_{N_2}\end{bmatrix}\right).\]
\end{lemma}

\begin{proof}[Proof of \cref{lem:block_circ_right}]
  We quote \eqref{eq:M_2} from \cref{lem:interleave} and consider that $P^{(k, N_2)} \otimes I_m$ is a permutation that changes the blockwise indices of $m \times p$ blocks (or, acting from the right, $p \times m$ blocks) exactly the way that $P^{(k, N_2)}$ changes the indices of a vector.
\end{proof}

\begin{lemma} \label{lem:diag_kron_perm}
  Given $k, n \in \N$ and $V_j \in \C^{m_j \times n_j}$ for $j \in [n]$ and setting $M = \sum_{j = 1}^n m_j, N = \sum_{j = 1}^N n_j,$ and $D = \diag(I_k \kron V_j)_{j = 1}^n \in \C^{k M \times k N}$, we have \[D = \diagcornmat{I_k \kron V_1}{I_k \kron V_n} = P_1 (I_k \otimes \diag(V_j)_{j = 1}^n) P_2^*\] where $P_1 = \Pc(P^{(n, k)}, (m_{j \mod_1 n})_{j = 1}^{kn})$ and $P_2 = \Pc(P^{(n, k)}, (n_{j \mod_1 n})_{j = 1}^{kn})$.
\end{lemma}

\begin{proof}[Proof of \cref{lem:diag_kron_perm}]
  We immediately reduce to the case $m_j = n_j = 1$ for all $j$ by observing that $P_1$ and $P_2$ will act on blockwise indices precisely as $P^{(n, k)}$ acts on individual indices.  Here, we replace $V_j$ with $v_j \in \C$, and note that $\diag(V_j)_{j = 1}^n = \diag(v)$.  Hence, we need only remark that \[\left(\diag(I_k \otimes v_\ell)_{\ell = 1}^n\right)_{((i_1 - 1)k + i_2) ((j_1 - 1)k + j_2)} = \left\{\begin{array}{r@{,\quad}l} v_{i_1} & i_1 = j_1 \ \text{and} \ i_2 = j_2 \\ 0 & \text{otherwise} \end{array}\right.,\] while 
  \begin{align*} (P^{(n, k)} (I_k \otimes \diag(v)) P^{(n, k)*})_{((i_1 - 1)k + i_2) ((j_1 - 1)k + j_2)} 
 & = (I_k \otimes \diag(v))_{((i_2 - 1)n + i_1) ((j_2 - 1)n + j_1)} 
  \\ &
  = \left\{\begin{array}{r@{,\quad}l} v_{i_1} & i_1 = j_1 \ \text{and} \ i_2 = j_2 \\ 0 & \text{otherwise} \end{array}\right..\end{align*}
\end{proof}
}\fi
\subsection{Matrix Representation and Conditioning of $\Ac$}
\label{sec:pty_zero_col}
To begin the discussion of the matrix representation of $\Ac$, we refresh our notation: $d, \delta, \dbar, s \in \N$ satisfy $d = \dbar s$ and $2 \delta - 1 \le d$.  We have $D \in \N$ (arbitrary for now) measurement vectors $\{m_j\}_{j \in [D]} \in \Cd$ satisfying $1 \in \supp(m_j) \subset [\delta]$, and we set $g_m^k = \diag(m_k m_k^*, m)$ for all $1 - \delta \le m \le \delta - 1$ and all $k \in [D]$.  The expressions $\Lc_{\{m_j\}}^s, T_{\delta, s},$ and $\Ac$ are defined in \cref{eq:lift_pty_sys,eq:T_delta_s,eq:pty_meas_op}.

We now consider the question of when $\Span \Lc_{\{m_j\}}^s = T_{\delta, s}$ and what the condition number of $\Ac$ will be.  As in \eqref{eq:diag_vec}, we vectorize $X$ by its diagonals\footnote{Notice that this will force $A$, the matrix representing $\Ac$, to be singular.  We expand on this later.} with $\Dc_\delta(X) \in \C^{d (2 \delta - 1)}$ and write $A \in \C^{\dbar D \times (2 \delta - 1) d}$ such that \begin{align*} \left(A \Dc_\delta(X) \right)_{(j-1) \dbar + \ell} &=  \left(A \colmatfl{\diag(X, 1 - \delta)}{\diag(X, \delta - 1)}\right)_{(j - 1) \dbar + \ell} = \Ac(X)_{(\ell - 1, j)} \\[.5em] &= \inner*{S^{s(\ell - 1)} m_j m_j^* S^{s (\ell - 1)}, X} = \sum_{m = 1 - \delta}^{\delta - 1} S^{s (\ell - 1)} g_m^{j*} \diag(X, m), 
\end{align*}
which gives the $(j - 1) \dbar + \ell\th$ row of $A$ as \[\begin{bmatrix} S^{s (\ell - 1)} g_{1 - \delta}^j \\ \vdots \\ S^{s (\ell - 1)} g_{\delta - 1}^j \end{bmatrix}^*\] so that, by \cref{lem:circ_transpose}, we have\footnote{For reference, we remark that $\circop^s(g_m^k) \in \C^{d \times \dbar}$ and $\circop(R_{\dbar} \Tc_s g_m^k) \in \C^{\dbar \times d}$.} \begin{equation} A = \begin{bmatrix} \circop^s(g_{1 - \delta}^1) & \cdots & \circop^s(g_{1 - \delta}^D) \\ \vdots & \ddots & \vdots \\ \circop^s(g_{\delta - 1}^1) & \cdots & \circop^s(g_{\delta - 1}^D) \end{bmatrix}^* = \begin{bmatrix} \circop(R_{\dbar} \Tc_s g_{1 - \delta}^1) & \cdots & \circop(R_{\dbar} \Tc_s g_{\delta - 1}^1) \\ \vdots & \ddots & \vdots \\ \circop(R_{\dbar} \Tc_s g_{1 - \delta}^D) & \cdots & \circop(R_{\dbar} \Tc_s g_{\delta - 1}^D) \end{bmatrix}. \label{eq:A_block_ptych} \end{equation}  However, because $T_{\delta, s} \subsetneq T_\delta$ when $s > 1$, this operator can never be invertible.  In fact, when $s > 1$, $A$ has several \emph{completely zero} columns, corresponding  to the coordinates of entries in $T_\delta / T_{\delta, s}$,\footnote{By $V / W$, where $W \subset V$, we mean $V \cap W^\perp$.} in the sense that $\Dc_\delta(T_\delta / T_{\delta, s}) \subset \Nul(A)$ and $\Row(A) \subset \Dc_\delta(T_{\delta, s})$, with equality iff $\left.\Ac\right|_{T_{\delta, s}(\Cdxd)}$ is invertible.

To compute the condition number of $\left.\Ac\right|_{T_{\delta, s}(\Cdxd)}$, we may take an orthogonal basis matrix $N$ for $\Dc_\delta(T_{\delta, s})$ and analyze the singular values of $AN$.  We construct $N$ by considering that, for each $m = 1-\delta, \ldots, \delta -1$, $R_\dbar \Tc_s g_m^j \in \C^{\dbar \times s}$ has $\min\{s, \delta - \abs{m}\}$ non-zero columns; specifically, these are columns $[1, \delta - m +1)$ for $m \ge 0$ and $[\abs{m} + 1, \delta + 1) \mod s$ for $m < 0$.  We denote and enumerate these intervals by $J_m = \{j^m_1, \ldots, j^m_{\abs{J_m}}\}$ and set $N_m := \begin{bmatrix} e_{j^m_1} & \cdots & e_{j^m_{\abs{J_m}}} \end{bmatrix}$ such that $R_\dbar (\Tc_s g_m^j) N_m$ has no zero columns.  The identity $\circop(AB) = \circop(A) (I \kron B)$ gives \[\circop(R_\dbar (\Tc_s g_m^j) N_m) = \circop(R_\dbar \Tc_s g_m^j) (I_s \kron N_m),\] so setting $N := \diag(I_s \kron N_m)_{m = 1 - \delta}^{\delta - 1}$, the columns of $N$ form a basis for $\Row(A)$ as desired.  This result is summarized in \cref{prop:basis_for_tdelta_s}.

\begin{proposition}
  Fix $s, \delta, d \in \N$ satisfying $s < \delta, s \mid d, \delta \le \frac{d+1}{2}$.  Define $J_m := [1, \delta - m + 1)$ for $m \in [0, \delta)$ and $J_m := [\abs{m} + 1, \delta + 1) \mod s$ for $m \in [1-\delta, 0)$.  Further setting $N_m := \begin{bmatrix} e_{j^m_1} & \cdots & e_{j^m_{\abs{J_m}}} \end{bmatrix}$ and $N := \diag(I_s \kron N_m)_{m = 1 - \delta}^{\delta - 1}$, we have $\Col(N) = \Dc_\delta(T_{\delta, s}(\Cdxd))$. \label{prop:basis_for_tdelta_s}
\end{proposition}

To prove a result analogous to that of \cref{thm:meas_cond} for $AN$, we will need to show that the restriction matrix $N$ commutes well with the permutations used in the condition number analysis of \cref{sec:con_number}, preserving the block-circulant structures that made the analysis possible.  Thankfully it does; following the intuition of \eqref{eq:interleaved_meas}, referring to our expression of $A$ in \eqref{eq:A_block_ptych}, and making use of \cref{lem:interleave,lem:block_circ_right}, we can arrive at \begin{equation} \begin{aligned} A' := P^{(\dbar, D)} A \left(P^{(\dbar, 2 \delta - 1)} \otimes I_s\right)^* &= \circop^D\left(P^{(\dbar, D)} \begin{bmatrix} R_{\dbar} \Tc_s g_{1 - \delta}^1 & \cdots & R_{\dbar} \Tc_s g_{\delta - 1}^1 \\ \vdots & \ddots & \vdots \\ R_{\dbar} \Tc_s g_{1 - \delta}^D & \cdots & R_{\dbar} \Tc_s g_{\delta - 1}^D \end{bmatrix}\right). \\ &= \circop^D\left(P^{(\dbar, D)} (I_D \otimes R_{\dbar}) \begin{bmatrix}  \Tc_s g_{1 - \delta}^1 & \cdots &  \Tc_s g_{\delta - 1}^1 \\ \vdots & \ddots & \vdots \\  \Tc_s g_{1 - \delta}^D & \cdots &  \Tc_s g_{\delta - 1}^D \end{bmatrix}\right). \end{aligned} \label{eq:A_prime} \end{equation}  

\noindent This may be reduced further by applying \cref{lem:diag_kron_perm}, which gives us that, setting \begin{align*} P_1 &= \Pc(P^{(2 \delta - 1, \dbar)}, (s)_{j = 1}^{\dbar (2 \delta - 1)}) = P^{(2 \delta - 1, \dbar)} \kron I_s \\ P_2 &= \Pc(P^{(2 \delta - 1, \dbar)}, (\min\{s, \delta - \abs{m}\})_{m = 1 - \delta}^{\delta - 1}) \\ N' &= \diag(N_m)_{m = 1 - \delta}^{\delta - 1}, \end{align*} we will have $N = \diag(I_{\dbar} \kron N_m)_{m = 1 - \delta}^{\delta - 1} =  P_1 (I_{\dbar} \kron N') P_2^*.$  This gives \begin{equation} A' (I_{\dbar} \kron N') = P^{(\dbar, D)} A P_1 (I_{\dbar} \kron N') = P^{(\dbar, D)} A N P_2. \label{eq:pddanp2}\end{equation}  We refer to \eqref{eq:A_prime} to obtain \[A'(I_{\dbar} \kron N') = \circop^D\left(P^{(\dbar, D)} (I_D \otimes R_{\dbar}) \begin{bmatrix}  \Tc_s g_{1 - \delta}^1 & \cdots &  \Tc_s g_{\delta - 1}^1 \\ \vdots & \ddots & \vdots \\  \Tc_s g_{1 - \delta}^D & \cdots &  \Tc_s g_{\delta - 1}^D \end{bmatrix} N' \right),\] which, along with \eqref{eq:pddanp2} and \cref{cor:circ_diag_condition}, gives us \cref{thm:pty_con_number}.

\begin{theorem} \label{thm:pty_con_number}
  Take $A$ as in \eqref{eq:A_block_ptych}, $N$ and $N_m$ as in \cref{prop:basis_for_tdelta_s}, \[H = P^{(\dbar, D)} (I_D \otimes R_{\dbar}) \begin{bmatrix}  \Tc_s g_{1 - \delta}^1 & \cdots &  \Tc_s g_{\delta - 1}^1 \\ \vdots & \ddots & \vdots \\  \Tc_s g_{1 - \delta}^D & \cdots &  \Tc_s g_{\delta - 1}^D \end{bmatrix} \diag(N_m)_{m = 1 - \delta}^{\delta - 1}\] and $M_j = \sqrt{\dbar}(f_j^{\dbar} \otimes I_D)^* H$ for $j \in [\dbar]$.  The condition number of $AN$ is given by \[\dfrac{\max\limits_{i \in [\dbar]} \sigma_{\max} (M_i)}{\min\limits_{i \in [\dbar]} \sigma_{\min} (M_i)}.\]  In particular, $\left.\Ac\right|_{T_{\delta, s}(\C^{d \times d})}$ is invertible if and only if each of the $M_i$ are of full rank.
\end{theorem}

%% file: sections/ptychography/ptych_recovery.tex
In this section, we  focus on algorithms by which we can  estimate $x_0$ from $T_{\delta, s}(\mathcal{A}^{-1}(y))$.  We begin with \cref{sec:blocky_block}, which discusses some improvements over the results of \cite{our_paper} that can be made to the magnitude estimation step of \cref{alg:pr_basic} (line 2).  Indeed, these improvements, which we call ``blockwise'' magnitude estimation, were first implemented (without theoretical analysis) in the numerical study of \cite{our_paper}. While they empirically delivered better results, no proof was provided to quantify the improvement, and we remedy this in \cref{sec:blocky_block}. Having dealt with the magnitude estimation step in the ptychographic setting, we then handle the phase estimation step in \cref{sec:pty_phase_est}. Finally, 
\cref{sec:pty_std_rec} describes and proves the robustness bounds for an algorithm  analogous to that of \cite{our_paper}, albeit fully capable of handling ptychographic shifts $s>1$, and taking advantage of the results in \cref{sec:blocky_block} and \cref{sec:pty_phase_est} for improved magnitude and phase estimation.  

\subsection{Blockwise Magnitude  Estimation}
\label{sec:blocky_block}
For the moment, we restrict our discussion to the special case of dense shifts in our measurements, where $s = 1$ as in \cite{our_paper} and \cref{ch:meas}.  In \cite{our_paper}, a functional but rudimentary technique to calculate the magnitudes of the entries of $\ux$ from $X \approx T_\delta(\ux \ux^*)$ was proposed.  It relied, as we do here, on first computing $X = \Ac^{-1}(\Ac(x_0 x_0^*) + n)$, where $\ux \in \Cd$ is the ground truth objective vector, $\Ac : \Cdxd \to \R^{d D}$ is the linear measurement operator determined by the masks $\{m_j\}_{j \in [D]}$, as defined, for example, in \eqref{eq:meas_op} of \cref{sec:meas_intro}, and $n \in \R^{d D}$ is arbitrary noise.  Then $X = \uX + \Ac^{-1}(n)$, where $\uX = T_\delta(\ux \ux^*)$, and the magnitude of $(\ux)_i$ was estimated by simply taking $\abs{x_i} = \sqrt{X_{ii}} \approx \sqrt{\uX_{ii}} = \abs{\ux_i}$.  This technique works, as was proven in \cite{our_paper}, but it was also seen that empirically,  a  more sophisticated technique does a much better job. Next, we will prove stability bounds for the improved technique, which we now describe.  

We notice that taking $\abs{x_i} = X_{ii}$ is equivalent to taking $\abs{x_i}$ to be the rank-1 approximation of the $1 \times 1, i\th$ diagonal block matrix of $X$, namely $\begin{bmatrix} X_{ii} \end{bmatrix}$, since these diagonal blocks are equal to the diagonal blocks of the \emph{untruncated} $\ux \ux^*$ when there is no noise.  However, given the width of the diagonal band in $T_\delta(\Cdxd)$, we could just as easily take blocks of size up to $\delta \times \delta$ and calculate their top eigenvectors; this would give us $2 \delta - 1$ estimates for each entry's magnitude, so we can combine them by averaging them together.  
To denote these blocks, we will set $\abs{X}^{(\ell)} = \diag(\one_{[\delta]_\ell}) \abs{X} \diag(\one_{[\delta]_\ell})$\footnote{Here, and in the remainder of this section, we emphasize that all indices of objects in $\Cd$ and $\Cdxd$ are taken modulo $d$.} and $u^{(\ell)}$ to be the top eigenvector of $\abs{X}^{(\ell)}$, normalized such that $\norm{u^{(\ell)}}_2 = \norm{\abs{X}^{(\ell)}}_2$.  We then produce our estimate of the magnitudes by taking $\abs{x} = \frac{1}{\delta} \sum_{\ell = 1}^d u^{(\ell)}$.  

Before formally stating this algorithm, we observe how it may be generalized.  Firstly, we notice that this method can easily handle arbitrary block sizes $m$ for the blockwise, eigenvector-based magnitude estimations by simply taking $\abs{X}^{(\ell, m)} = \diag(\one_{[m]_\ell}) \abs{X} \diag(\one_{[m]_\ell})$,  $m \le \delta$ -- although this would require changing the denominator in the averaging step, using $\frac{1}{m} \sum_{\ell = 1}^d u^{(\ell, m)}$.  Proceeding  further, we can generalize this technique to use any collection $\{J_i\}_{i = 1}^N, J_i \subset [d]$ satisfying \begin{equation} \begin{gathered} [d] \subset \bigcup_{i = 1}^N J_i \quad \text{and} \quad \one_{J_i} \one_{J_i}^* \in T_{\delta}(\Cdxd) \end{gathered} \label{eq:T_delta_cover}.\end{equation}  We will call any collection satisfying \eqref{eq:T_delta_cover} a \emph{$(T_\delta, d)$-covering} (or just \emph{covering} when $T_\delta$ and $d$ are clear from context), and the process of estimating magnitudes of $\ux$ from $X$ with respect to a $(T_\delta, d)$-covering is described in \cref{alg:blocky_block}.\footnote{We remark that this definition and the recovery algorithm are very obviously extensible to the use of $T_{\delta, s}$ instead of $T_\delta$.  In fact, this is a \emph{restriction}, if we consider in \eqref{eq:T_delta_cover} that $T_{\delta, s} \subset T_\delta$.  The definition, therefore, of a $(T_{\delta, s}, d)$-covering, is made by analogy to \eqref{eq:T_delta_cover}.}  It is worth remarking that the ``averaging step,'' specified in line 3, is optimal in the least-squares sense.  We set $P_{J_i} = \diag(\one_{J_i})$ to be the orthogonal projections onto the coordinate subspace associated with $J_i$ and consider that the vectors $u^{(J_i)}$ (in line 2) represent estimates of the projections $P_{J_i} \abs{\ux}$.  The least squares solution to $P_{J_i} u = u^{(J_i)}$, or \[\colmatfl{P_{J_1}}{P_{J_N}} u = \colmatfl{u^{(J_1)}}{u^{(J_N)}}\] is obtained by taking the pseudoinverse.  In this case, we have \[\colmatfl{P_{J_1}}{P_{J_N}}^* \colmatfl{P_{J_1}}{P_{J_N}} = \sum_{i = 1}^N P_{J_i}^* P_{J_i} = \sum_{i = 1}^N P_{J_i} = \diag(\mu),\] with $\mu_j = \abs{\{i : j \in J_i\}}$ as in line 1 of \cref{alg:blocky_block}.  Considering that $P_{J_i} u^{(J_i)} = u^{(J_i)}$, we have \[u = \colmatfl{P_{J_1}}{P_{J_N}}^{\dag} \colmatfl{u^{(J_1)}}{u^{(J_N)}} =  \diag(\mu)^{-1} \colmatfl{P_{J_1}}{P_{J_N}}^* \colmatfl{u^{(J_1)}}{u^{(J_N)}} = D_\mu^{-1}\left(\sum_{i = 1}^N u^{(J_i)}\right).\]

\begin{algorithm}[htbp]
\renewcommand{\algorithmicrequire}{\textbf{Input:}}
\renewcommand{\algorithmicensure}{\textbf{Output:}}
\caption{Blockwise Magnitude Estimation}
\label{alg:blocky_block}
\begin{algorithmic}[1]
    \REQUIRE $X \in T_{\delta, s}(\H^d)$, typically assumed to be an approximation $X \approx T_{\delta, s}(\ux \ux^*)$.  A $(T_{\delta, s}, d)$-covering $\{J_i\}_{i \in [N]}$.
    \ENSURE An estimate $\abs{x}$ of $\abs{\ux}$.
    \STATE For $j \in [d]$, set $\mu_j = \abs{\{i : j \in J_i\}}$ to be the number of appearences the index $j$ makes in $\{J_i\}$.
    \STATE For $i \in [N]$, set $u^{(J_i)}$ to be the leading eigenvector of $\abs{X^{(J_i)}} = \diag(\one_{J_i}) \abs{X} \diag(\one_{J_i})$, normalized such that $\norm{u^{(J_i)}}_2 = \sqrt{\norm{X^{(J_i)}}_2}$.
    \STATE Return $\abs{x} = D_\mu^{-1}\left(\sum_{i = 1}^N u^{(J_i)}\right)$.
    \end{algorithmic}
\end{algorithm}

We will denote the output of \cref{alg:blocky_block} by $\abs{x} = \BlkMag(X, \{J_i\})$.  In an overloading of notation, when the covering consists of intervals of length $m$, in the sense that $J_i = [m]_i$ for $i = 1, \ldots, d$, we will also write it as $\BlkMag(X, \{[m]_i\}_{i \in [d]}) = \BlkMag(X, m)$.  To include the case of ptychography, where these intervals are shifted by more than 1, we write \begin{equation}\begin{gathered} \BlkMag(X, m) = \BlkMag(X, \{[m]_i\}_{i \in [d]}), \ \text{and} \\ \BlkMag(X, (m, s)) = \BlkMag(X, \{[m]_{1 + s (\ell-1)}\}_{\ell \in [\dbar]}), \end{gathered}\label{eq:bm_overload} \end{equation} where  $\dbar = \frac{d}{s}$ is an integer.  In this way, the improved magnitude estimation technique used in the numerical experiments of \cite{our_paper} is  $\abs{x} = \BlkMag(X, \delta)$, while $\abs{x_i} = \sqrt{X_{ii}} = \BlkMag(X, 1)$. We can  prove a bound on the error of the estimate.\bigskip

\begin{proposition} \label{prop:blocky_block}
  Let $\{J_i\}_{i \in [N]}$ be a $(T_{\delta, s}, d)$-covering, and suppose $\uX = T_{\delta, s}(\ux \ux^*)$ for some $\ux \in \Cd$.  Using the notation of \cref{alg:blocky_block} (in particular, $\mu_j$ is as in line 1, and $\uu^{(J_i)} = \abs{\diag(\one_{J_i}) \ux}$), given $X \in T_{\delta, s}(\H^d)$, we have that the output $\abs{x}$ satisfies
  \begin{equation}
  \begin{aligned}
    \BlkMag(\uX, \{J_i\}) &= \abs{\ux} \\
    \norm{\BlkMag(X, \{J_i\}) - \abs{\ux}}_2 &\le \dfrac{\max_j \mu_j}{\min_j \mu_j} \dfrac{1 + 2 \sqrt{2}}{\min_i \norm{\uu^{(J_i)}}_2} \norm{\uX - X}_F.
  \end{aligned}
  \label{eq:blocky_rec}
  \end{equation}
  As special cases for $T_\delta(\Cdxd)$, we have
  \begin{equation}
    \begin{aligned}
      \norm{\BlkMag(X, m) - \ux}_2 &\le \dfrac{1 + 2 \sqrt{2}}{\min_i \norm{\uu^{[m]_i}}_2} \norm{X - \uX}_F \\
      \norm{\BlkMag(X, \delta) - \ux}_2 &\le \dfrac{1 + 2 \sqrt{2}}{\min_i \norm{\uu^{[\delta]_i}}_2} \norm{X - \uX}_F \\
      \norm{\BlkMag(X, 1) - \ux}_2 &\le \dfrac{\norm{\diag(X - \uX)}_F}{\min_i \abs{\ux_i}} 
    \end{aligned}
    \label{eq:blocky_spec}
  \end{equation}
\end{proposition}

\begin{proof}[Proof of \cref{prop:blocky_block}]
  The first inequality of \eqref{eq:blocky_rec} is clear, since line 2 of \cref{alg:blocky_block} will always return $\uu^{(J_i)} = \one_{J_i} \circ \abs{\ux}$, so line 3 will give \[ \left(D_\mu^{-1}\left(\sum_{i = 1}^N \uu^{(J_i)}\right)\right)_j = \frac{1}{\mu_j} \sum_{i = 1}^N \abs{\ux}_j \one_{j \in J_i} = \abs{\ux}_j.\]  
  The second comes by writing \begin{equation} \norm*{D_\mu^{-1}(\sum_{i = 1}^N u^{(J_i)} - \uu^{(J_i)})}_2^2 \le \left(\frac{1}{\min_j \mu_j}\right)^2 \norm*{\sum_{i = 1}^N (u^{(J_i)} - \uu^{(J_i)})}_2^2. \label{eq:bb1} \end{equation}  From there, we consider that the $j\th$ term in the summation \begin{equation}\norm*{\sum_{i = 1}^N (u^{(J_i)} - \uu^{(J_i)})}_2^2 = \sum_{j = 1}^d \left(\sum_{i = 1}^N u^{(J_i)} - \uu^{(J_i)}\right)_j^2\label{eq:bb2}\end{equation} has at most $\max_k \mu_k$ nonzero summands, so, by $(\sum_{i=1}^n a_i)^2 \le n \sum_{i=1}^n a_i^2$, we have \begin{equation}\norm*{\sum_{i = 1}^N (u^{(J_i)} - \uu^{(J_i)})}_2^2 \le \max_j \mu_j \sum_{i = 1}^N (u^{(J_i)} - \uu^{(J_i)})^2.\label{eq:bb3}\end{equation}  We then apply 
  Lemma A.2 of \cite{our_paper}\footnote{Here, we use the substitution $\eta \norm{\x_0}_2 = \frac{\norm{X - X_0}_F}{\norm{\x_0}_2}$.} to get \begin{equation} \sum_{i = 1}^N (u^{(J_i)} - \uu^{(J_i)})^2 \le (1 + 2 \sqrt{2}) \dfrac{\norm{\uu^{(J_i)} \uu^{(J_i)*} - X^{(J_i)}}_F^2}{\norm{\uu^{(J_i)}}_2^2} \label{eq:bb4} \end{equation}  In this expression, we consider that, in the summation $\sum_{i = 1}^N \norm{\uu^{(J_i)} \uu^{(J_i)*} - X^{(J_i)}}_F^2,$ the term $(\uX_{ij} - X_{ij})^2$ appears at most $\max\{\mu_i, \mu_j\}$ times, such that $\sum_{i = 1}^N \norm{\uu^{(J_i)} \uu^{(J_i)*} - X^{(J_i)}}_F^2 \le \max_j \mu_j \norm{\uX - X}_F^2$.  Using this substitution, combining \eqref{eq:bb1}--\eqref{eq:bb4}, and taking the square root of both sides gives \[\norm*{D_\mu^{-1}(\sum_{i = 1}^N u^{(J_i)} - \uu^{(J_i)})}_2 \le \dfrac{\max_j \mu_j}{\min \mu_j} \dfrac{1 + 2 \sqrt{2}}{\min_i \norm{\uu^{(J_i)}}_2} \norm{\uX - X}_F\] as desired.
  The first two inequalities of \eqref{eq:blocky_spec} are immediate by observing that $\mu_j = m$ when the covering is $\{[m]_i\}_{i \in [d]}$.  The third comes from setting $\epsilon_i = X_{ii} - (\uX)_{ii}$ and writing
  \begin{align*}
    \norm{\abs{x} - \abs{\ux}}_2^2 &= \sum_{i = 1}^d \left(\sqrt{X_{ii}} - \abs{\ux}_i\right)^2 = \sum_{i = 1}^d \left(\sqrt{\abs{\ux}_i^2 + \epsilon_i} - \sqrt{\abs{\ux}_i^2}\right)^2 \\
    &= \sum_{i = 1}^d \left(\frac{((\abs{\ux}_i^2 + \epsilon_i) - \abs{\ux}_i^2)^2}{\sqrt{\abs{\ux}_i^2 + \epsilon_i} + \abs{\ux}_i}\right)^2 \le \sum_{i = 1}^d \frac{\epsilon_i^2}{\min_i{\abs{\ux}_i^2}} 
    = \frac{\norm{\diag(X - \uX)}_F^2}{\min_i \abs{\ux}_i^2}
  \end{align*}  
\end{proof}

One immediate benefit from \cref{eq:blocky_spec} is that the estimation error from $\BlkMag(X, 1)$ no longer scales poorly with $d^{1/4}$ as in Theorem 5 of \cite{our_paper}. 
In the error bound for $\BlkMag(X, m)$ in \cref{eq:blocky_spec}, we notice that the bound is \emph{strictly decreasing} with $m$, since, for $m_1 > m_2$, we have \[\min_i \norm{\ux^{[m_1]_i}}_2^2 \ge (m_1 - m_2) \min_i \abs{(\ux)_i}^2 + \min_i \norm{\ux^{[m_2]_i}}_2^2 \ge m_1 \min_i \abs{(\ux)_i}^2.\]  Also, considering \eqref{eq:blocky_rec}, it is clear that $\BlkMag(X, \delta)$ gives the absolute best bound over all $\{J_i\}$, since any $(T_\delta, d)$-covering $\{J_i\}_{i \in [N]}$ satisfies $J_i \subset [\delta]_\ell$ for some $\ell$.  This gives that $\min_i \norm{\uu^{(J_i)}}_2 \le \min_i \norm{\uu^{[\delta]_i}}_2$, and obviously $\frac{\max_j \mu_j}{\min_j \mu_j} \ge 1$, so the bound for $\BlkMag(X, \{J_i\})$ in \eqref{eq:blocky_rec} can never be better than that for $\BlkMag(X, \delta)$ in \eqref{eq:blocky_spec}. We also remark that two easy ways to ensure $\frac{\max_j \mu_j}{\min_j \mu_j} = 1$ is minimized are to take $J_0 \subset [\delta]_0$ and let $J_\ell = J_0 + \ell$ be a ``cyclic'' covering, or to let $\{J_i\}$ be a partition of $[d]$.  These strategies will be relevant in \cref{sec:pty_std_rec}, but in the case of $s = 1$, $\BlkMag(X, \delta)$ always has the optimal bound for magnitude estimation error.
\subsection{Phase Estimation for Ptychography}\label{sec:phase_est}
We now consider a simple algorithm for estimating the relative phase between the entries of $x_0$, from an estimate of $T_{\delta,s}(x_0 x_0^*)$ obtained, say, by inverting $\mathcal{A}$ on $T_{\delta,s}(\C^{d\times d})$. This algorithm is a  generalization of an analogous one in \cite{our_paper}, and is similar to another algorithm, albeit with slightly different bounds, that was very recently presented in \cite{melnyk2019phase}.  For simplicity, we restrict our attention to the case where $s$ divides $\delta$, and $\delta$ divides $d$, so that $\delta=s\bar{\delta}$ and $d=\delta\bar{d}$. We defer an analysis of the general case with arbitrary $s,\delta,d$, as well as a study of more sophisticated algorithms (e.g., based on semidefinite programming, as in \cite{BP_thesis}) to other work. The following simple, yet helpful, lemma holds. 
\begin{lemma}\label{lem:ptych_rel}
If $\delta,s, d \in \N$ are such that $s$ divides $\delta$, and $\delta$ divides $d$, then 
\[ T_{\delta,s}(\one_d\one_d^*) = T_{\delta/s}(\one_{d/s}\one_{d/s}^*)\otimes (\one_s\one_s^*).  \]
\end{lemma}
\begin{proof}[Proof of \cref{lem:ptych_rel}]
The proof consists of simply enumerating the non-zero indices of each of $T_{\delta,s}(\one_d\one_d^*)$  and $T_{\delta/s}(\one_{d/s}\one_{d/s}^*)\otimes (\one_s\one_s^*)$. In both cases, these are $\cup_{\ell\in[\bar{d}]_0} [\delta]_{s\ell+1}^2 $, by definition for the first and by the definition of the Kronecker product for the second.
\end{proof}
It then follows from the properties of Kronecker products that the eigenvalues of $T_{\delta,s}(\one_d\one_d^*)$ are simply the pairwise products of the eigenvalues of $T_{\delta/s}(\one_{d/s}\one_{d/s}^*)$ and $\one_s\one_s^*$. Equally importantly, the eigenvectors of $T_{\delta,s}(\one_d\one_d^*)$ are the pairwise Kronecker products of those of $T_{\delta/s}(\one_{d/s}\one_{d/s}^*)$ and $\one_s\one_s^*$. So, the entire eigendecomposition of $T_{\delta,s}(\one_d\one_d^*)$ is fully known, as $T_{\delta/s}(\one_{d/s}\one_{d/s}^*)$ is a circulant matrix, hence diagonalized by the discrete Fourier transform, and $\one_s\one_s^*$ is a rank-one matrix. In particular, the normalized eigenvector of $T_{\delta,s}(\one_d\one_d^*)$ corresponding to the leading eigenvalue is simply $\one_d/\sqrt{d}$. Now, consider that $\diag\left(\frac{x_0}{|x_0|}\right)$ is unitary and that \[T_{\delta,s}\left(\frac{x_0x_0^*}{|x_0 x_0^*|}\right) = \diag\left(\frac{x_0}{|x_0|}\right)T_{\delta,s}(\one_d\one_d^*)\diag\left(\frac{x_0}{|x_0|}\right)^*\] so $T_{\delta,s}\left(\frac{x_0x_0^*}{|x_0 x_0^*|}\right)$ and $T_{\delta,s}(\one_{d}\one_{d}^*)$ are similar. Given the eigendecomposition $T_{\delta,s}(\one_{d}\one_{d}^*)= V \Lambda V^*$, 
\begin{equation}T_{\delta,s}\left(\frac{x_0x_0^*}{|x_0 x_0^*|}\right) 
= 
\left(\diag\left(\frac{x_0}{|x_0|}\right)V\right)~\Lambda~\left( \diag\left(\frac{x_0}{|x_0|}\right)V\right)^* \label{eq:T_relations}\end{equation}
%
is itself an eigendecomposition. Together, these observations imply that the eigenvector corresponding to the leading eigenvalue of $T_{\delta,s}\left(\frac{x_0x_0^*}{|x_0 x_0^*|}\right) $ is  
$\diag\left(\frac{x_0}{|x_0|}\right)\one_d/\sqrt{d}$, i.e., it is the vector of phases $\frac{x_0}{|x_0|}$, up to a harmless normalization! This implies that in the absence of noise, we can obtain $\sgn(x_0)$ easily via an eigendecomposition. It remains to show that this procedure is robust to noise.
\begin{theorem}\label{thm:phase_stable}
Fix $s<\delta<d \in \N$ such that $\delta$ divides $d$ and $s$ divides $\delta$. Let $\widetilde{X}_0 = T_{\delta,s}\left(\frac{x_0x_0^*}{|x_0 x_0^*|}\right)$, and let $\widetilde{X}$ be as in line 2 from \cref{alg:pty_pr}. Let $\widetilde{x}=\sgn(v)$, with $v$ being the leading eigenvector of $\widetilde{X}$. If $\| \widetilde{X_0} - \widetilde{X}\|_F \leq \tilde\eta $
 for some $\tilde\eta>0$, then there exists a positive constant $C$ such that 
\[ \min_{\theta \in [0, 2\pi]} \| \sgn{(x_0)} - e^{\ii \theta} \tilde{x} \|_2 \leq C \tilde\eta 
\frac{{d}^{5/2}}{{\delta}^2} s
\] 
\end{theorem}
\begin{proof}[Proof of \cref{thm:phase_stable}]
Given an undirected graph $G=(V,E)$ with $d$ vertices, let $D$ be the diagonal matrix of its degrees and $W$ be its adjacency matrix, and define its connection Laplacian $L:=I-D^{-1/2}WD^{-1/2}$. Denoting by $0\leq \lambda_1 \leq ... \leq \lambda_d$,  the spectral gap associated with $G$ is  $\tau :=\lambda_2$. By Theorem 4 of \cite{our_paper}, noting that the graph associated with the ``adjacency matrix" $T_{\delta,s}{(\one_{d/s}\one_{d/s}^*)}$ is $(2\delta-s)$-regular, there exists a constant $C'>0$, such that
\[
\min_{\theta \in [0, 2\pi]} \| \sgn{(x_0)} - e^{\ii \theta} \tilde{x} \|_2 \leq C' \frac{\|\widetilde{X}-\widetilde{X}_0\|_F}{\tau \sqrt{2\delta-s}} \leq C'\frac{\tilde\eta}{\tau\sqrt{2\delta-s}}.
\]

It remains to bound $\tau$ from below and to that end we denote by $\nu_1\geq \nu_2 \geq ... \nu_{d/s}$ the eigenvalues of $T_{\delta/s}(\one_{d/s}\one_{d/s}^*)$.  We then invoke Lemma 2 of \cite{our_paper} to conclude that there exists a constant $C''>0$, such that $\min\limits_{j\neq 1}(\nu_1 - |\nu_j|) \geq C''\frac{(\delta/s)^3}{(d/s)^2}=C\frac{\delta^2}{d^2s}$, and the same conclusion holds for $T_{\delta,s}(\one_d \one_d^*)$ by \cref{lem:ptych_rel}. Now, since the eigenvector of $\widetilde{X}$ corresponding to its leading eigenvalue is an eigenvector of $L$ corresponding to its smallest eigenvalue, we have $\tau \geq C''\frac{\delta^2}{d^2s},$ so that 
\[
\min_{\theta \in [0, 2\pi]} \| \sgn{(x_0)} - e^{\ii \theta} \tilde{x} \|_2 \leq \frac{C'\tilde\eta d^2 s}{C''\delta^2\sqrt{2\delta-s}} \leq C \tilde\eta \frac{d^2}{\delta^{5/2}}s.
\]
%
%
\end{proof}
\begin{corollary}\label{cor:phase_stable}
Let $\widetilde{X}_0 = T_{\delta,s}\left(\frac{x_0x_0^*}{|x_0 x_0^*|}\right)$, and let $\widetilde{X}$ be as in line 2 from \cref{alg:pty_pr}. Let $\widetilde{x}=\sgn(v)$, with $v$ being the leading eigenvector of $\widetilde{X}$. If $\sigma_{min}:=\sigma_{min}(\mathcal{A}|_{T_{\delta,s}(\C^{d\times d})})$ and $\kappa:=\kappa(\mathcal{A}|_{T_{\delta,s}(\C^{d\times d})})$ denote the smallest singular value and condition number of $\mathcal{A}|_{T_{\delta,s}(\C^{d\times d})}$, respectively, then 
\[ \min_{\theta \in [0, 2\pi]} \| \sgn{(x_0)} - e^{\ii \theta} \tilde{x} \|_2 \leq C \cdot 
\frac{{d}^{2} \cdot s  }{{\delta}^{5/2} \cdot  } \cdot \frac{ \sigma_{min} ^{-1} \cdot \|n\|_2 }{\min\limits_{j}|(x_0)_j|^2} \] 
and
\[ \min_{\theta \in [0, 2\pi]} \| \sgn{(x_0)} - e^{\ii \theta} \tilde{x} \|_2 \leq C \cdot 
\frac{{d}^{2} \cdot s  }{{\delta}^{5/2} \cdot  } \cdot \frac{ \kappa \cdot \|X\|_F }{\SNR \cdot \min\limits_{j}|(x_0)_j|^2} \] 

\end{corollary}
\begin{proof}
We proceed as in the proof of Lemma 6 in \cite{our_paper}. Setting $N=X-X_0$ we have  
\begin{align*}
|(\widetilde{X}_0)_{j,k}-(\widetilde{X})_{j,k}| &=  \left| (\widetilde{X}_0)_{j,k} - \sgn\left( \frac	{X_{j,k}}{|(X_0)_{j,k}|}\right) \right| 
\\ &
\leq \left| (\widetilde{X}_0)_{j,k} - \frac{X_{j,k}}{|{(X_0)}_{j,k}|} \right| + \left| \frac{X_{j,k}}{|{(X_0)}_{j,k}|} - \sgn\left( \frac	{X_{j,k}}{|(X_0)_{j,k}|}\right) \right|
\\ &
\leq 2\left| (\widetilde{X}_0)_{j,k} - \frac{X_{j,k}}{|{(X_0)}_{j,k}|} \right|
= 2 \left|\frac{N_{j,k}}{|{(X_0)}_{j,k}|}\right|.
\end{align*}

This gives 
$\norm{\widetilde{X}-\widetilde{X}_0}_F \leq \frac{2\norm{N}_F}{\min\limits_{j,k}|(X_0)_{j,k}|}\leq  \frac{2\norm{N}_F}{\min\limits_{j}|(x_0)_j|^2},$ where $\norm{N}_F$ can be bounded by $\sigma_{min}(\mathcal{A})^{-1}\|n\|_2$ or by $\kappa(\mathcal{A})\frac{\norm{X}_F}{\SNR}$. Combining this with \cref{thm:phase_stable} yields the result. 
\end{proof}
%
\label{sec:pty_phase_est}

\subsection{Recovery Algorithm for Ptychography}
\label{sec:pty_std_rec}

We are now ready to prove that \cref{alg:pty_pr} stably produces an estimate of $\ux$;
%
%
%
%
 \Cref{thm:pty_rec} provides a bound on the accuracy of the output of \cref{alg:pty_pr}.  \smallskip



\begin{theorem} \label{thm:pty_rec}
 Let $d,\dbar,\delta,s\in \N$ be such that $s$ divides $\delta$, and $\delta$ divides $d$ with $d=\dbar s$.   Suppose we have a family of masks $\{m_j\}_{j \in [D]} \in \Cd$ of support $\delta$ and a $(T_{\delta, s}, d)$-covering $J_\ell:=\{[m]_{1 + s (\ell-1)}\}_{\ell \in [\dbar]}$. Denote by $\Ac$ the corresponding linear operator, and by $\sigma_{min}$ and $\kappa$ the smallest singular value and condition number of $\Ac|_{T_{\delta,s}(\C^{d\times d})}$.  Further let $\ux \in \Cd, n \in \R^{\dbar D}$ be arbitrary and set $\uX = T_{\delta,s}(\ux \ux^*), X = \Ac^{-1}(\Ac(\uX) + n)$.  Define  $\SNR = \frac{\norm{\Ac(X_0)}_2}{\norm{n}_2}$.  
Then the output $x$ of \cref{alg:pty_pr} satisfies \begin{equation}\begin{aligned} \mintheta \norm{x - \eit \ux}_2 &\le 
C \cdot  \sigma_{min} ^{-1} \cdot \left(  
\frac{{d}^{2} \cdot s  }{{\delta}^{5/2}   } \cdot \frac{ \|x_0\|_\infty  }{\min\limits_{j}|(x_0)_j|^2}  ~ + ~ \frac{1}{\min\limits_\ell \| \diag(\one_{J_\ell}) \ux \|_2} \right) \cdot \|n\|_2 
\\
\mintheta \norm{x - \eit \ux}_2 &\le 
C \cdot  \kappa \cdot \left(  
\frac{{d}^{2} \cdot s  }{{\delta}^{5/2}   } \cdot \frac{ \|x_0\|_\infty  }{\min\limits_{j}|(x_0)_j|^2}  ~ + ~ \frac{1}{\min\limits_\ell \| \diag(\one_{J_\ell}) \ux \|_2} \right) \cdot \frac{\|X_0\|_F}{\SNR}.
\end{aligned} \label{eq:pty_rec}\end{equation}
\end{theorem}


\begin{proof}[Proof of \cref{thm:pty_rec}]
We set $\tbx = \sgn(x), \tilde{x}_0 = \sgn(\ux)$ and use the triangle inequality to obtain
  \begin{align*}
    \mintheta \norm{x - \eit \ux}_2 &\le \norm{x_0}_\infty \mintheta \norm{\tbx - \eit \tilde{x}_0}_2 + \norm{\abs{x} - \abs{\ux}}_2.
  \end{align*}
We bound each of the summands separately. Specifically, we directly use \cref{cor:phase_stable} to bound $\mintheta \norm{\tbx - \eit \tilde{x}_0}_2$ and \cref{prop:blocky_block},   \cref{eq:blocky_rec},  to bound $\norm{\abs{x} - \abs{\ux}}_2$. For the latter, we use the fact that our choice of $d,s,\delta$ yields $\frac{\max_j \mu_j}{\min_j \mu_j} = 1$,
%
%
  %
and we use the stability of the linear system $\mathcal{A}$ 
to bound $\norm{X - \uX}_F$ above by $\sigma_{\min}^{-1} \norm{n}_2$ and $\kappa \frac{\norm{\uX}_F}{\SNR}$.  
\end{proof}

%

\if{
\subsection{Blockwise Vector Synchronization Method}
\label{sec:vecky_vec}
There is another method for retrieval from $X = \Ac^{-1}(y)$, based largely on the ideas behind the Blockwise Magnitude Estimation of \cref{sec:blocky_block}, that permits a very straightforward and satisfyingly general recovery method.  The idea here is to use something analogous to \cref{alg:blocky_block}, but which recovers the magnitudes \emph{and phases} of each block, and then stitches them together with the vector synchronization technique discussed in \cref{app:vec_sync}.  The blockwise, eigenvector-based step, called Blockwise Vector Estimation, is stated in \cref{alg:vecky_vec}, and the complete recovery process is stated in \cref{alg:pty_vec}.

Unfortunately, the theory behind this strategy is so far fairly weak: we are able to prove a recovery bound on it, but it suffers from tremendously poor scaling with problem dimension.  The proof is also very technical, so we defer it to \cref{app:vecky_vec}.  For the moment, we pose \cref{alg:pty_vec} as a theoretically unvetted, but intuitive and empirically promising technique, and defer a more optimistic view of its performance to the numerical studies of \cref{sec:ptych_num}.

\begin{algorithm}[htbp]
\renewcommand{\algorithmicrequire}{\textbf{Input:}}
\renewcommand{\algorithmicensure}{\textbf{Output:}}
\caption{Blockwise Vector Estimation}
\label{alg:vecky_vec}
\begin{algorithmic}[1]
    \REQUIRE $X \in T_{\delta, s}(\H^d)$, typically assumed to be an approximation $X \approx T_{\delta, s}(\ux \ux^*)$.  A connected $(T_{\delta, s}, d)$-covering $(\{J_i\}_{i \in [N]}, G)$.
    \ENSURE Estimates $x^{(J_i)}$ of $\diag(\one_{J_i}) \ux$.
    \STATE For $i \in [N]$, set $x^{(J_i)}$ to be the leading eigenvector of $X^{(J_i)} = \diag(\one_{J_i}) X \diag(\one_{J_i})$, normalized such that $\norm{x^{(J_i)}}_2 = \sqrt{\norm{X^{(J_i)}}_2}$.
    \STATE Return $\BlkVec(X, (\{J_i\}, G)) = \rowmatfl{x^{(J_1)}}{x^{(J_N)}}$.
\end{algorithmic}
\end{algorithm}

Despite the pessimistic result of \cref{prop:vecky_vec}, we expect this method to perform fairly well.  Recalling the discussion of \cref{sec:weighted_graph} concerning the advantages of using the magnitudes of $X$ in the angular synchronization process (namely, that this will prioritize accurate relative phases between larger entries of $x$), it seems like a potentially sensible step forward to let the blockwise magnitude estimation process of \cref{alg:blocky_block} simultaneously produce an estimate of the phases of $x$.  This process operates by finding the eigenvectors of dense, approximately rank-1 matrices, and the eigenvector perturbation bound used in the proof of \cref{prop:blocky_block} works just as well for complex matrices.  Synthesizing these blockwise estimates by calculating appropriate relative phases between them (e.g., by appealing to the vector synchronization process defined in \cref{app:vec_sync}) and averaging as in line 3 of \cref{alg:blocky_block} seems satisfyingly symmetric, but our lack of theoretical understanding of the performance of the vector synchronization process prevents us from achieving a stronger theoretical promise on this method.

\begin{algorithm}
\renewcommand{\algorithmicrequire}{\textbf{Input:}}
\renewcommand{\algorithmicensure}{\textbf{Output:}}
\caption{Phase Retrieval by Vector Synchronization}
\label{alg:pty_vec}
\begin{algorithmic}[1]
  \REQUIRE Measurements $y \approx \Ac(\ux \ux^*) \in \R^{\dbar D}$, where $\Ac$ is invertible on $T_{\delta, s}(\Cdxd)$ as in \eqref{eq:pty_meas_op}.  A connected $(T_{\delta, s}, d)$-covering $(\{J_i\}_{i \in [N]}, G')$.
  \ENSURE $x \in \Cd$ with $x \approx \eit \ux$ for some $\theta \in [0, 2 \pi]$.
  \STATE For $j \in [d]$, set $\mu_j = \abs{\{i : j \in J_i\}}$ to be the number of appearences the index $j$ makes in $\{J_i\}$.
  \STATE Set $V = \BlkVec(X, (\{J_i\}, G))$.
  \STATE Compute $\hZ$, the solution to \eqref{eq:ang_sync_sdp} with $L_{\rm vs}$ as defined in \eqref{eq:lap_vecky} using $V$ and $G'$ and set $\hz = \sgn(u),$ where $u$ is the top eigenvector of $\hZ$.
  \STATE Return $x = D_\mu^{-1} V \hz$.
\end{algorithmic}
\end{algorithm}

One remark that must be made before stating the vector synchronization-based recovery algorithm is that having a simple $(T_{\delta, s}, d)$-covering will not suffice for this method; we require a slightly stronger set of conditions on the sets $J_i \subset [d]$, which will allow vector synchronization to be performed on their blockwise estimates of $\ux$.  Namely, we define a \emph{connected $(T_{\delta, s}, d)$-covering} to be an ordered pair $(\{J_i\}, G')$ satisfying that $\{J_i\}$ is a $(T_{\delta, s}, d)$-covering and $G' = (V = [N], E)$ is a connected graph with \[E \subset \{(i, j) : J_i \cap J_j \neq \emptyset\},\] 
which will permit the blocks $x^{(J_i)}$ associated with each $J_i$ to have a meaningful sense of relative phase over a connected graph.  The synthesis process will be to take, for each $i \in [N]$, $x^{(J_i)}$ to be the top eigenvector of $X^{(J_i)} = \diag(\one_{J_i}) X \diag(\one_{J_i})^*$ normalized such that $\norm{x^{(J_i)}}_2 = \sqrt{\norm{X^{(J_i)}}_2}$ and to synchronize them by taking the solution $\hz$ of \eqref{eq:ang_sync} for \begin{equation} L_{\rm vs} = D_{\rm vs} - W_{\rm vs} \circ X_{\rm vs}, \label{eq:lap_vecky} \end{equation} defining $D_{\rm vs}, W_{\rm vs}, X_{\rm vs}$ as in \eqref{eq:vec_sync_mats}.  With $\mu_j$ defined as in line 1 of \cref{alg:blocky_block}, the final estimate of $\ux$ is rendered as $x = D_\mu^{-1} \rowmatfl{x^{(J_1)}}{x^{(J_n)}} \hz$, which simply averages together the phase-synced $x^{(J_i)} \hz_i$ terms.  With this, we state \cref{alg:vecky_vec,alg:pty_vec}.\footnote{In \cref{alg:pty_vec}, notice that $\BlkVec(X, (\{J_i\}, G))$ is the output of \cref{alg:vecky_vec}.}
}\fi

%% file: sections/ptychography/ptych_num.tex
\if{We now numerically evaluate the new structures and results enumerated and studied in this chapter.  In \cref{sec:pn_blocky}, we compare the blockwise magnitude estimation technique stated in \cref{sec:blocky_block} to the previous magnitude estimation strategy of reading the magnitudes of $x$ directly from the main diagonal of $X \approx T_{\delta, s}(x x^*)$.  \cref{sec:pty_cond} considers the condition number $\kappa(\left.\Ac\right|_{T_{\delta, s}(\Cdxd)})$ of the measurement operator in the ptychographic case, as well as the spectral gap $\tau_G$ associated with the unweighted graph whose adjacency matrix is $T_{\delta, s}(\one_d \one_d^*) - I_d$.  

We note that, in \cref{sec:pty_cond
}, we exclusively use masks $m_j$ that are randomly generated.  Unfortunately, despite the structure uncovered in \cref{sec:con_number_ptych}, we have not yet discovered a way to import the notion of a local Fourier measurement system, which produced such a satisfyingly simple and physically realizable family of masks for which the condition number was predictably controlled, to the case of general shifts $s > 1$.  As of this dissertation, the author is unaware of any deterministic construction $\{m_j\}_{j \in [s(2 \delta - s)]}$ that scales with $s, \delta,$ and $d$ to consistently produce invertible measurement operators $\Ac$.
}\fi

\begin{figure}
  \centering
  \begin{subfigure}[b]{0.3\textwidth}
    \centering
    \scalebox{0.9}{
    \begin{tikzpicture}[ampersand replacement=\&,baseline=-\the\dimexpr\fontdimen22\textfont2\relax]
    \matrix (m)[matrix of math nodes,left delimiter=(,right delimiter=)]
            {
              * \& * \& * \&   \&   \&   \& * \& * \\
              * \& * \& * \&   \&   \&   \&   \&   \\
              * \& * \& * \& * \& * \&   \&   \&   \\
                \&   \& * \& * \& * \&   \&   \&   \\
                \&   \& * \& * \& * \& * \& * \&   \\
                \&   \&   \&   \& * \& * \& * \&   \\
              * \&   \&   \&   \& * \& * \& * \& * \\
              * \&   \&   \&   \&   \&   \& * \& * \\
            };

            \begin{pgfonlayer}{myback}
              \fhighlight[blue!30]{m-1-1}{m-1-1}
              \fhighlight[blue!30]{m-2-2}{m-2-2}
              \fhighlight[blue!30]{m-3-3}{m-3-3}
              \fhighlight[blue!30]{m-4-4}{m-4-4}
              \fhighlight[blue!30]{m-5-5}{m-5-5}
              \fhighlight[blue!30]{m-6-6}{m-6-6}
              \fhighlight[blue!30]{m-7-7}{m-7-7}
              \fhighlight[blue!30]{m-8-8}{m-8-8}
            \end{pgfonlayer}
    \end{tikzpicture}}
    \caption{$\BlkMag(X, 1)$}
  \end{subfigure}
  \begin{subfigure}[b]{0.3\textwidth}
    \centering
    \scalebox{0.9}{
    \begin{tikzpicture}[ampersand replacement=\&,baseline=-\the\dimexpr\fontdimen22\textfont2\relax]
    \matrix (m)[matrix of math nodes,left delimiter=(,right delimiter=)]
            {
              * \& * \& * \&   \&   \&   \& * \& * \\
              * \& * \& * \&   \&   \&   \&   \&   \\
              * \& * \& * \& * \& * \&   \&   \&   \\
                \&   \& * \& * \& * \&   \&   \&   \\
                \&   \& * \& * \& * \& * \& * \&   \\
                \&   \&   \&   \& * \& * \& * \&   \\
              * \&   \&   \&   \& * \& * \& * \& * \\
              * \&   \&   \&   \&   \&   \& * \& * \\
            };

            \begin{pgfonlayer}{myback}
              \fhighlight[blue!30]{m-1-1}{m-3-3}
              \fhighlight[blue!30]{m-3-3}{m-5-5}
              \fhighlight[blue!30]{m-5-5}{m-7-7}
              \fhighlight[blue!30]{m-7-1}{m-8-1}
              \fhighlight[blue!30]{m-1-7}{m-1-8}
              \fhighlight[blue!30]{m-7-7}{m-8-8}
              \fhighlight[blue!50]{m-1-1}{m-1-1}
              \fhighlight[blue!50]{m-3-3}{m-3-3}
              \fhighlight[blue!50]{m-5-5}{m-5-5}
              \fhighlight[blue!50]{m-7-7}{m-7-7}
            \end{pgfonlayer}
    \end{tikzpicture}}
    \caption{$\BlkMag(X, (3, 2))$}
  \end{subfigure}
  \begin{subfigure}[b]{0.3\textwidth}
    \centering
    \scalebox{0.9}{
    \begin{tikzpicture}[ampersand replacement=\&,baseline=-\the\dimexpr\fontdimen22\textfont2\relax]
    \matrix (m)[matrix of math nodes,left delimiter=(,right delimiter=)]
            {
              * \& * \& * \&   \&   \&   \& * \& * \\
              * \& * \& * \&   \&   \&   \&   \&   \\
              * \& * \& * \& * \& * \&   \&   \&   \\
                \&   \& * \& * \& * \&   \&   \&   \\
                \&   \& * \& * \& * \& * \& * \&   \\
                \&   \&   \&   \& * \& * \& * \&   \\
              * \&   \&   \&   \& * \& * \& * \& * \\
              * \&   \&   \&   \&   \&   \& * \& * \\
            };

            \begin{pgfonlayer}{myback}
              \fhighlight[blue!30]{m-1-1}{m-2-2}
              \fhighlight[blue!30]{m-3-3}{m-4-4}
              \fhighlight[blue!30]{m-5-5}{m-6-6}
              \fhighlight[blue!30]{m-7-7}{m-8-8}

            \end{pgfonlayer}
    \end{tikzpicture}}
    \caption{$\BlkMag(X, (2, 2))$}
  \end{subfigure}
  \caption{Blocks used for blockwise magnitude estimation in $T_{3, 2}(\C^{8 \times 8})$}
  \label{fig:blocky_pic}
\end{figure}
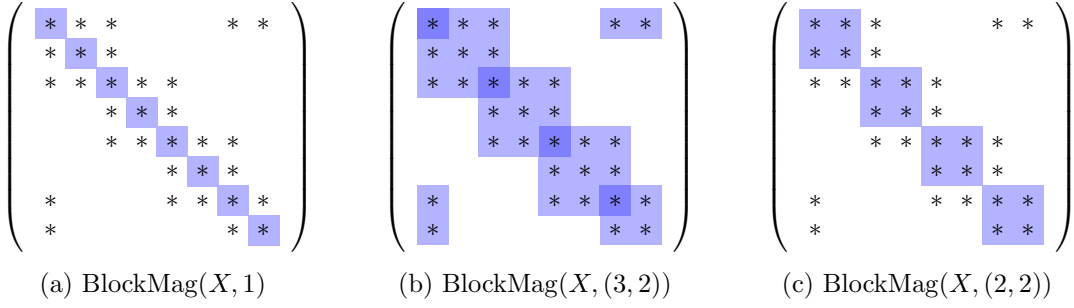

\if{
\subsection{Numerical Study of Magnitude Estimation Techniques}
\label{sec:pn_blocky}
The blockwise magnitude estimation technique of \cref{sec:blocky_block} was first introduced during the numerical analysis of \cref{alg:phaseRetrieval1} presented in \cref{sec:NumEval}.  This method had been invented as a mere intuition on how to slightly improve the numerical performance of \cref{alg:phaseRetrieval1}, and because it appeared to remove a few decibels of relative error from the reconstruction process as a whole, it was used in the evaluations.  A theoretical analysis of this technique, however, was not developed until \cref{sec:blocky_block}, when bounds on the relative error in the magnitudes were proven, and indeed the blockwise, eigenvector-based method was shown to have a sharper theoretical bound.  Nonetheless, a direct numerical comparison between this strategy and the diagonal strategy was never made, so we present this comparison here.

Recall from the discussion of \cref{sec:blocky_block} that the blockwise magnitude estimation technique consists of applying \cref{alg:blocky_block} to obtain $\BlkMag(X, \{J_i\}_{i \in [N]})$, where $J_i \subset [d]$ forms a $(T_{\delta, s}, d)$-covering, defined in \cref{eq:T_delta_cover}.  It was remarked that taking $\abs{x}_i = \sqrt{\abs{X_{ii}}}$, as in line 5 of \cref{alg:phaseRetrieval1}, was actually a special case of $\BlkMag$, being equivalent to $\BlkMag(X, 1)$, using the notation of \cref{eq:bm_overload}, while the ``maximally dense'' method described in \cref{sec:MagEstImpNumerical} is equivalent to $\BlkMag(X, \delta)$.  Therefore, a comparison of these two techniques is, in some sense, merely an evaluation of $\BlkMag$ over different parameters.  We remark that, again according to \eqref{eq:bm_overload}, these techniques are trivially extensible to the ptychographic case of $s > 1$ by taking $\BlkMag(X, 1)$ and $\BlkMag(X, (\delta, s))$.

A third strategy, worthy of consideration, is briefly mentioned at the end of \cref{sec:blocky_block}.  In this case, we let $\{J_i\}$ be a partition of $[d]$.  Of course, taking $J_i = [1]_i$, as in $\BlkMag(X, 1)$, is one example of such a partition, but in \crefrange{fig:bb_chart1}{fig:bb_chart3}, the ``Part'n'' data points refer to calculating $\BlkMag(X, (s\floor{\frac{\delta - 1}{s}}, s))$, which is in some sense the ``largest'' partition.  To get a sense of how $\BlkMag(X, 1),$ $\BlkMag(X, (\delta, s)),$ and $\BlkMag(X, (s\floor{\frac{\delta - 1}{s}}, s\floor{\frac{\delta - 1}{s}}))$ scale with $s$, examples are illustrated in \cref{fig:blocky_pic}.
}\fi

\if{

With this, we present the results of this numerical experiment in \cref{fig:bb_charts}, where we have compared the relative error in magnitude estimation of these three techniques as a function of signal-to-noise ratio ($\SNR$) for $s = 1, 3$, and $5$, when $d = 60$ and $\delta = 6$.  This gives us a comparison in the dense case (when $s = 1$), an overlap of $50\%$ ($s = 3$), and a minimal overlap of $83.3\%$ ($s = 5$).

For each $\SNR$, we randomly generated $128$ objective vectors $\ux \in \CN(0, I_d)$, and computed $\uX = T_{\delta, s}(\ux \ux^*)$.  We then randomly drew a Hermitian perturbation $N = T_{\delta, s}(N' + N'^*)$ where $N'_{ij} \iid \CN(0, 1)$, and scaled it to ensure $\norm{\uX}_F / \norm{N}_F = \SNR$.  We then wrote $X = \uX + N$ and calculated $\abs{x_{\diag}} = \BlkMag(X, 1),b$ $\abs{x_{\delta}} = \BlkMag(X, (\delta, s))$, and $\abs{x_{\rm part}} = \BlkMag(X, (s\floor{\frac{\delta - 1}{s}}, s\floor{\frac{\delta - 1}{s}}))$, along with the relative errors $\frac{\norm{\abs{x_{\cdot}} - \abs{\ux}}_2}{\norm{\ux}_2}$ for each.
These relative errors were averaged over the $128$ trials and plotted against $\SNR$ for each technique in \cref{fig:bb_charts}.

\begin{figure}
  \centering
  \begin{subfigure}[b]{.49\textwidth}
    \centering
    \includegraphics[width=\textwidth,trim={.1in 2.5in .8in 2.4in}]{figs/bb_chart1}
    \caption{$s = 1$}
    \label{fig:bb_chart1}
  \end{subfigure}
  \begin{subfigure}[b]{.49\textwidth}
    \centering
    \includegraphics[width=\textwidth,trim={0.1in 2.5in .8in 2.4in}]{figs/bb_chart2}
    \caption{$s = 3$}
    \label{fig:bb_chart2}
  \end{subfigure}
  \begin{subfigure}[b]{.49\textwidth}
    \centering
    \includegraphics[width=\textwidth,trim={.1in 2.5in .8in 2.4in}]{figs/bb_chart3}
    \caption{$s=5$}
    \label{fig:bb_chart3}
  \end{subfigure}
  \caption{Relative error in magnitude estimation vs. SNR.  $d = 60, \delta = 6, s \in \{1, 3, 5\}$}
  \label{fig:bb_charts}
\end{figure}

The results mostly affirm the intuition of \cref{sec:MagEstImpNumerical} and the theory of \cref{prop:blocky_block}: the larger block sizes show consistently stronger performance than the diagonal-only recovery method originally put forth.  In fact, once the $\SNR$ becomes usable (at $10^1$ or $10^2$ -- at $\SNR = 10^{-1}$ or $10^0$, the measurements are clearly useless), the relative error on the blockwise magnitude estimates is reduced by roughly a factor of five!  Not too surprisingly, the partition method produces results comparable to those of the ``full blocks'' of $\BlkMag(X, (\delta, s))$, since the block sizes differ by at most $\min\{s, \delta - s\}$, and the gain that $x_{\rm part}$ experiences from using off-diagonal entries to inform its magnitude estimates is on the same order as the gain experienced by $x_\delta$.

On the other hand, one observation that is \emph{not} expected is that the performance of none of these methods appears to deteriorate at all with increasing $s$.  This may be explained by supposing that it is the actual \emph{eigenvector strategy}, rather than averaging over several estimates, that gives this method its strength.  indeed, this intuition would also somewhat explain the comparability between $x_{\rm part}$ and $x_\delta$, since the block sizes for each are roughly equal, and $x_\delta$'s main advantage is that each magnitude arises from an average of several estimates, rather than a single estimate as in the case of a partition.  However, none of these intuitions are formalized any further at the moment.

\subsection{Conditioning and Spectral Gap}
\label{sec:pty_cond}

In this section, we consider the measures we have that gauge the degree to which the recovery methods of \cref{sec:ptych_recov}, especially \cref{alg:pty_pr}, magnify noise in the measurements $\Ac(\ux \ux^*) + n$.  Specifically, we are going to numerically evaluate $\kappa$, the condition number of the linear operator $\Ac$ defined in \eqref{eq:pty_meas_op}, and $\tau_G$, the spectral gap of the graph that will appear in the angular synchronization phase (specifically, line 3 of \cref{alg:pty_pr}).

As already mentioned, there are so far no known constructions of masks $\{m_j\}_{j \in [D]}$ that have theoretical conditioning bounds (or even a theoretical guarantee of invertability), so \cref{fig:pty_cond} describes the distribution of $\kappa(\Ac)$ over masks given by $m_j \iid \CN(0, I_\delta), j \in [s(2 \delta - s)]$, much like \cref{fig:arbchart} in \cref{sec:rand_fam}.  In this experiment, we have fixed $d = 60$ and $\delta = 13$, but we vary $s$ so that we can see whether and how much the level of overlap affects the distribution of condition numbers.  For each value of $s$, we have drawn $1024$ sets of masks and calculated the condition number $\kappa$ of each set.  Note that, to compare the distributions of $\kappa$ between values of $s$, in lieu of presenting a histogram with four overlaid sets of bins, we have opted for a line plot which traces the heights that the bins would otherwise take.

\begin{figure}[htb]
  \centering
  \includegraphics[width=.6\textwidth,trim={.2in 2.5in 1in 2.5in}]{figs/pty_arbchart}
  \caption{Distribution of condition numbers.  $d = 60, \delta = 13$.}
  \label{fig:pty_cond}
\end{figure}

Seen side by side, the distributions appear roughly equal -- they each have their peak between $10^3$ and $10^4$ and possess relatively little probability mass after $10^5$.  It is difficult to comment on whether these results are predictable, however: considering the expressions obtained for $A$ in \eqref{eq:A_prime} and its condition number in \cref{thm:pty_con_number}, where $g_m^j$ is a vector of products of Gaussian random variables, and the $g_m^j$ are dependent across $m$, there is no obvious intuition for the singular values of $A$.  Nonetheless, it is reassuring to see that the distribution of $\kappa$ does not radically deteriorate when we introduce shifts $s > 1$, and somewhat interesting to see that, after slightly worsening when $s$ goes through $2$ and $6$, at $s = 12$ (the maxiumum allowed shift for $\delta = 13$), the distribution comes back to a peak more commensurate with that of $s = 1$, although it does come with a thicker tail.

\Cref{fig:tau_ptych} is a very simple display of the spectral gap $\tau_G$ of the graphs associated with the subspaces $T_{\delta, s}$ for different values of $s$.  To be precise, we define $G_{\delta, s}^d = (V_{\delta, s}^d = [d], E_{\delta, s}^d)$ to be the graph on $d$ vertices such that the adjacency graph of $G_{\delta, s}^d$ is $A_G(d, \delta, s) := T_{\delta, s}(\one_d \one_d^*) - I_d$.  For simplicity of notation, we will refer to $\tau_{G_{\delta, s}^d}$ as $\tau_G(d, \delta, s)$, or $\tau_G$ when the arguments are superfluous.  Recall from \cref{sec:Perturb,sec:ang_sync_improve} that the spectral gap of this graph is an important component of the error bound that is proven for any of the angular synchronization techniques studied in this dissertation -- the bound is inversely proportional to $\tau_G$ for the faster, eigenvector-based recovery method, and to $\sqrt{\tau_G}$ for the SDP recovery, so we want to have a spectral bound that is as large as possible.

However, $\tau_G$ is, in some sense, a measure of ``how connected'' the graph $G$ is, and is in fact often called the \emph{algebraic connectivity} of a graph in the literature \cite{deabreu2006algebraicconnectivity,fiedler1973algebraic_connectivity}.  Indeed, both of these works remark that $\tau_G$ is monotonically decreasing as edges are removed from a graph, so the introduction of larger shifts $s > 1$ can only serve to make the error bounds of \cref{thm:SpecGraphPertBound,thm:improved_spec_pert} worse.  \Cref{fig:tau_ptych} numerically illustrates the dependence of $\tau_G(d, \delta, s)$ on each of its variables.  For these charts, we constructed the Laplacian $L(d, \delta, s) = \diag(A_G(d, \delta, s) \one_d) - A_G(d, \delta, s)$ for each $(d, \delta, s)$ and simply calculated the smallest eigenvalue of $P L P^*$, where $P$ is an orthogonal projection onto $\one_d^\perp = \Col(L(d, \delta, s))$.  Notice, in each case, that the choices of parameters $(d, \delta, s)$ may seem somewhat erratic; this is simply because it is required that $\dbar = d / s$ be an integer.  In any case, the data points we have produced are sufficient to illustrate the most pressing trends.

\begin{figure}[htb]
  \centering
  \begin{subfigure}[b]{.49\textwidth}
    \centering
      \includegraphics[width=\textwidth,trim={.2in 175pt .9in 175pt}]{figs/tau_ptych}
    \caption{Spectral gap of $G_{\delta, s}^d$ vs. $d$.  $\delta = 16$.  $s \in \{1, 4, 8, 12, 14, 15\}$.}
    \label{fig:tau_ptych_d}
  \end{subfigure}
  \hfill
  \begin{subfigure}[b]{.49\textwidth}
    \centering
      \includegraphics[width=\textwidth,trim={.2in 181pt .9in 175pt}]{figs/tau_ptych_delta_ann}
    \caption{Spectral gap of $G_{\delta, s}^d$ vs. $\delta$.  $d = 240$.  $\text{Overlap} \approx 100\%, 50\%, 25\%, \frac{1}{\delta}$.}
    \label{fig:tau_ptych_delta}
  \end{subfigure}
  \caption{Dependence of $\tau_G$ on $d, \delta$, and $s$.}
  \label{fig:tau_ptych}
\end{figure}

These results are qualitatively similar to those for the condition number: the behavior deteriorates in the direction we expect, but not catastrophically.  In \cref{fig:tau_ptych_d}, for example, we can see that, fixing $\delta$ and changing $d$, $\tau_G$ asymptotically behaves like $\bigO(1 / d^2)$ for each level of overlap; recalling from \cref{lem:EigGap} that, for $s = 1$, we had $\tau_G \approx \bigO(\delta^3 / d^2)$, this is a reasonable observation.  Nicely enough, as the shifts $s$ get larger, the spectral gap doesn't appear to shrink too horrendously: between $s = 1$ and $s = \delta - 1 = 15$, the minumum and maximum possible shifts, $\tau_G$ only decreases by roughly a factor of $16$.

Interestingly, on the other hand, $s$ appears to play a larger role in shaping how $\tau_G$ varies with $\delta$.  \Cref{fig:tau_ptych_delta} demonstrates this dependence, where we have now fixed $d = 240$, and we see that the relationship $\tau_G \approx \bigO(\delta^3)$ does \emph{not} hold over different values of $s$.  For $s = \delta - 1$, the maximal shift (and minimal overlap), $\tau_G$ doesn't even appear to level out at $\bigO(\delta^2)$.  From a rough reading of \cref{fig:tau_ptych_delta}, we see that going from $s = 1$ to $s = \delta - 1$ at $\delta = 16$ costs about a factor of $2^4$ in $\tau_G$, whereas the same operation -- full overlap to minimum overlap -- at $\delta = 81$ costs a factor of $2^6$.  This increased ``spectral cost'' may be accounted for by considering that, with higher values of $\delta$, the difference between the graphs $G_{\delta, 1}^d$ and $G_{\delta, \delta - 1}^d$ is a far greater number of edges, so the graph's ``connectivity'' may be expected to suffer more.  For the present moment, we rest with these observations and defer a formalization of these results to future study.
}\fi


%% file: sections/appendix/appendix.tex
\section{Appendix}
\input{con_number_appendix.tex}

%% file: sections/appendix/con_number_appendix.tex
\subsection{Interleaving Operators and Circulant Structure} 
\label{sec:interlemma}

To set the stage for the proof of \cref{thm:meas_cond}, we introduce a certain collection of permutation operators and study their interactions with circulant and block-circulant matrices.  The structure we identify here will be of much use to us in unraveling the linear systems we encounter in our model for phase retrieval with local correlation measurements.
%
For $\ell, N_1, N_2 \in \N, v \in \C^{\ell N_1}, k \in [\ell],$ and $H \in \C^{\ell N_1 \times N_2}$, we define the block circulant operator $\circop^{N_1}$ by
\begin{align*}
  \circop_k^{N_1}(v) &= \begin{bmatrix} v & S^{N_1} v & \cdots & S^{(k - 1)N_1} v \end{bmatrix} \\
  \circop_k^{N_1}(H) &= \begin{bmatrix} H & S^{N_1} H & \cdots & S^{(k - 1) N_1}H \end{bmatrix},
\end{align*}
where, as with $\circop(\cdot)$, when we omit the subscript we define $\circop^{N_1}(H) = \circop_\ell^{N_1}(H)$ and $\circop^{N_1}(v) = \circop_\ell^{N_1}(v)$.  We now proceed with the following lemmas; the first establishes the inverse of $P^{(d, N)}$.

\begin{lemma} \label{lem:interleave_inverse}
  For $d, N \in \N,$ we have \[(P^{(d, N)})^{-1} = P^{(d, N) *} = P^{(N, d)}.\]
\end{lemma}

\begin{proof}[Proof of \cref{lem:interleave_inverse}]
Simply take $v \in \C^{d N}$ and calculate, for $i \in [d], j \in [N]$,
  \begin{align*}
    (P^{(d, N)} P^{(N, d)} v)_{(i - 1) N + j} &= (P^{(d, N)} (P^{(N, d)} v))_{(i - 1) N + j} 
    = (P^{(N, d)} v)_{(j - 1) d + i} 
    = v_{(i - 1) N + j},
  \end{align*}
  with these equalities coming from the definition in \eqref{eq:interleave_def}.  
\end{proof}

We now observe some useful ways in which the interleaving operators commute with the construction of circulant matrices.

\begin{lemma}\label{lem:interleave}
  Suppose $V_i \in \C^{k \times n}, v_{ij} \in \C^k, w_j \in \C^{k N_1}$ for $i \in [N_1], j \in [N_2]$ and
  \begin{gather*}
    M_1 = \begin{bmatrix} \circop(V_1) \\ \vdots \\ \circop(V_{N_1}) \end{bmatrix},\quad
    M_2 = \begin{bmatrix} \circop^{N_1}(w_1) & \cdots & \circop^{N_1}(w_{N_2}) \end{bmatrix},\ \text{and} \\
    M_3 = \begin{bmatrix} \circop(v_{11}) & \cdots & \circop(v_{1 N_2}) \\ \vdots & \ddots & \vdots \\ \circop(v_{N_1 1}) & \cdots & \circop(v_{N_1 N_2}) \end{bmatrix}.\end{gather*}
  Then
  \begin{align}
    P^{(k, N_1)} M_1 &= \circop^{N_1}\left(P^{(k, N_1)} \begin{bmatrix} V_1 \\ \vdots \\ V_{N_1} \end{bmatrix}\right) \label{eq:M_1} \\
    M_2 P^{(k, N_2)*} &= \circop^{N_1}\left(\begin{bmatrix} w_1 & \cdots & w_{N_2} \end{bmatrix}\right) \label{eq:M_2} \\
    P^{(k, N_1)}M_3P^{(k, N_2)*} &= \circop^{N_1}\left(P^{(k, N_1)} \begin{bmatrix} v_{11} & \cdots & v_{1 N_2} \\ \vdots & \ddots & \vdots \\ v_{N_1 1} & \cdots & v_{N_1 N_2} \end{bmatrix}\right). \label{eq:M_3}
  \end{align}
\end{lemma}

\begin{proof}[Proof of lemma \ref{lem:interleave}]
  We index the matrices to check the equalities.  For \eqref{eq:M_1}, we take $(a, b, \ell, j) \in [d] \times [N_1] \times [k] \times [n]$ and have 
  \begin{align*}
    (P^{(k, N_1)} M_1)_{(a-1)N_1 + b, (\ell - 1) n + j} &= (M_1)_{(b - 1) k + a, (\ell - 1)n + j} 
    = \begin{bmatrix} S^{\ell - 1} V_1 \\ \vdots \\ S^{\ell - 1} V_{N_1} \end{bmatrix}_{(b - 1)k + a, j} 
    \\ &
    = (S^{\ell - 1}V_b)_{a, j} = (V_b)_{a + \ell - 1, j}
  \end{align*}
  and
  \begin{align*}
    \circop^{N_1}\left(P^{(k, N_1)} \begin{bmatrix} V_1 \\ \vdots \\ V_{N_1} \end{bmatrix}\right)_{(a-1)N_1 + b, (\ell - 1) n + j} &= \left(P^{(k, N_1)} \begin{bmatrix} V_1 \\ \vdots \\ V_{N_1} \end{bmatrix}\right)_{(a - 1)N_1 + b + (\ell-1)N_1, j} \\
    &= 
    (V_b)_{a + \ell - 1, j}
  \end{align*}
  For \eqref{eq:M_2}, we take $(a, b, j) \in [k] \times [N_2] \times [k N_1]$ and have
 \[(P^{(k, N_2)} M_2^*)_{(a - 1)N_2 + b, j} = (M_2)_{j, (b - 1) k + a} = (w_b)_{j + (a - 1)N_1}\]
  and
  \[\left(\circop^{N_1}\left(\begin{bmatrix} w_1 & \cdots & w_{N_2} \end{bmatrix}\right)\right)_{j, (a - 1)N_2 + b} = (S^{N_1(a - 1)} w_b)_j = (w_b)_{j + N_1(a - 1)},\]
and  \eqref{eq:M_3} follows immediately by combining \eqref{eq:M_1} and \eqref{eq:M_2}.
\end{proof}

\Cref{lem:interkron} introduces useful identities relating  interleaving operators to kronecker products.

\begin{lemma}\label{lem:interkron}
  For $v \in \C^N, V = \rowmat{V}{1}{\ell} \in \C^{N \times \ell}, A = \rowmat{A}{1}{m} \in \C^{d \times m}$, and $B_i \in \C^{m \times k}, i \in [\ell]$, we have
  \begin{align}
    P^{(d, N)} (v \kron A) &= A \kron v
    \label{eq:interkron_vec} \\
    P^{(d, N)} (V \kron A) &= \rowmat{A \kron V}{1}{\ell}
    \label{eq:interkron_mat} \\
    P^{(d, N)} (V \kron A) P^{(\ell, m)} &= A \kron V
    \label{eq:interkron_swap} \\
    (V \kron A) \diagmat{B}{1}{\ell} &= \rowmatfun{V_@ \kron A B_@}{1}{\ell}
    \label{eq:kron_diag}
  \end{align}
\end{lemma}

\begin{proof}[Proof of \cref{lem:interkron}]
  For \eqref{eq:interkron_vec}, we see that, for $i, j, k \in [d] \times [N] \times [m]$, we have
  \begin{align*}
    (P^{(d, N)} v \otimes A)_{(i - 1) N + j, k} &= (v \otimes A)_{(j - 1) d + i, k} 
    = v_j A_{i k}, \ \text{while} \\
    (A \kron v)_{(i - 1) N + j, k} &= A_{i k} v_j,
  \end{align*}
  and \eqref{eq:interkron_mat} follows by considering that $V \otimes A = \rowmatfun{V_@ \kron A}{1}{\ell}.$  To get \eqref{eq:interkron_swap}, we trace the positions of columns, considering that $(V \kron A) e_{(i - 1) m + j} = V_j \kron A_i$.  From \eqref{eq:interkron_mat}, we observe that $P^{(d, N)} (V \kron A) e_{(i - 1) m + j} = A_j \kron V_i,$ so \begin{align*}
    P^{(d, N)} (V \kron A) P^{(m, \ell)} e_{(j - 1) \ell + i} &= P^{(d, N)} (V \kron A) e_{(i - 1) m + j} = A_j \kron V_i = (A \kron V) e_{(j - 1) \ell + i}.
  \end{align*}
 As for \eqref{eq:kron_diag}, we remark that \begin{gather*} (V \kron A) \diagmat{B}{1}{\ell}
    = (V \kron A) \rowmatfun{e^\ell_@ \kron B_@}{1}{\ell}  \\ = \rowmatfun{(V \kron A) (e_@^\ell \kron B_@)}{1}{\ell} = \rowmatfun{V_@ \kron A B_@}{1}{\ell}.\end{gather*}
\end{proof}
\noindent The following lemma on the Kronecker product is standard (e.g., Theorem 13.26 in \cite{laub2004matrix}).

\begin{lemma}\label{lem:kronvec}
  We have $\vec(A B C) = (C^T \kron A) \vec(B)$ for any $A \in \C^{m \times n}, B \in \C^{n \times p}, C \in \C^{p \times k}.$  In particular, for $a, b \in \Cd, \vec a b^* = \conj{b} \kron a$, and \begin{equation}\label{id:kronsimp} \vec E_{jk} (\vec E_{j' k'})^* = E_{k k'} \kron E_{j j'}. \end{equation}
\end{lemma}
The next lemma covers the standard result concerning the diagonalization of circulant matrices, as well as a generalization to block-circulant matrices.
\begin{lemma}
  For any $v \in \Cd$, we have \begin{equation} \circop(v) = F_d \diag(\sqrt{d} F_d^* v) F_d^* = \sqrt{d} \sum_{j = 1}^d (f_j^{d *} v) f_j^d f_j^{d *} \label{eq:circ_dft_diag} \end{equation}
  Suppose $V \in C^{k N \times m}$, then $\circop^N(V)$ is block diagonalizable by \begin{equation} \circop^N(V) = \left(F_k \otimes I_N\right) \left(\diag(M_1, \ldots, M_k)\right) \left(F_k \otimes I_m\right)^*,  \label{eq:circ_dft_blk} \end{equation} where \begin{equation} \sqrt{k}\left(F_k \otimes I_N\right)^* V = \begin{bmatrix} M_1 \\ \vdots \\ M_k \end{bmatrix}, \quad \text{or} \quad M_j = \sqrt{k} (f_j^k \otimes I_N)^* V \label{eq:M_ell}\end{equation} \label{lem:circ_diag}
\end{lemma}

\begin{proof}[Proof of lemma \ref{lem:circ_diag}]
  The diagonalization in \eqref{eq:circ_dft_diag} is a standard result: see, e.g., Theorem 7 of \cite{gray2006circulant}.

  To prove \eqref{eq:circ_dft_blk}, we set $V_i$ to be the $k \times m$ blocks of $V$ such that $V^* = \rowmat{V^*}{1}{k}$ and begin by observing that, for $u \in \C^k$ and $W \in \C^{m \times p}$, the $\ell\th$ $k \times p$ block of $\circop^N(V)(u \otimes W)$ is
   \[\left(\circop^N(V)(u \otimes W)\right)_{[\ell]} = \sum_{i = 1}^k u_i (S^{N (i - 1)}V)_\ell W = \sum_{i = 1}^k u_i V_{\ell - i + 1} W.\]  Taking $u = f_j^k$ and $W = I_m$, this gives \begin{align*} \left(\circop^N(V)(f_j^k \otimes I_m)\right)_{[\ell]} &= \frac{1}{\sqrt{k}}\sum_{i = 1}^k \omega_k^{(j - 1) (i - 1)} V_{\ell - i + 1} I_m = \frac{1}{\sqrt{k}} \omega_k^{(j - 1) (\ell - 1)} \sum_{i = 1}^k \omega_k^{-(j - 1)(i - 1)} V_i \\ &= (f_j^k)_\ell \left(\sqrt{k} (f_j^k \otimes I_N)^* V \right) = (f_j^k)_\ell M_j. \end{align*}  This relation is equivalent to the statement of the lemma, i.e., having \[\circop^N(V) (f_j^k \otimes I_m) = (f_j^k \otimes M_j) = (f_j^k \otimes I_N) M_j.\] 
\end{proof}

Lemma \ref{lem:circ_diag} immediately gives the following corollary regarding the conditioning of $\circop^N(V)$, with which we return to considering spanning families of masks.

\begin{corollary} \label{cor:circ_diag_condition}
  With notation as in lemma \ref{lem:circ_diag}, the condition number of $\circop^N(V)$ is \[\dfrac{\max\limits_{i \in [k]} \sigma_{\max} (M_i)}{\min\limits_{i \in [k]} \sigma_{\min} (M_i)}.\]
\end{corollary}

\subsection{Lemmas on Block Circulant Structure}
\label{sec:ptych_lem}
We begin with \cref{lem:circ_transpose}, which describes the transposes of block circulant matrices.  For this lemma and the remainder of this section, the reader is advised to recall the definitions of $R_k$ and $P^{(d, N)}$ from \cref{sec:notation} and \eqref{eq:interleave_def}, as well as $\Pc(P, \{k_i\})$ and $\Tc_N$ from \cref{sec:pty_new_op}.
\begin{lemma}
\label{lem:circ_transpose} Given $k, N, m \in \N$ and $V \in \C^{kN \times m}$, we have \[\circop^N(V)^* = \circop^m\left( (R_k \otimes I_m) \Tc_N(V) \right).\] 
\end{lemma}
\begin{proof}[Proof of \cref{lem:circ_transpose}]
  Suppose $V_i$ are the $N \times m$ blocks of $V$, such that $V = \left[V_1^T \cdots V_k^T\right]^T$.  Indexing blockwise, we have $\circop^N(V)_{[ij]} = V_{i - j + 1}$, so that $\circop^N(V)^*_{[ij]} = V_{j - i + 1}^*$.  In other words,
  \[
  \circop^N(V)^* = \begin{bmatrix} V_1^* & V_2^* & \cdots & V_N^* \\ V_N^* & V_1^* & \cdots & V_{N - 1}^* \\ \vdots & & \ddots & \vdots \\ V_2^* & V_3^* & \cdots & V_1^* \end{bmatrix} = \circop^m((R_k \otimes I_m) \Tc_N(V) )
 . \]
\end{proof}

\Cref{lem:block_circ_right,lem:diag_kron_perm} provide identities for a few block matrix structures that will be of interest.
\begin{lemma}\label{lem:block_circ_right}
  Given $N_1, N_2, k, m \in \N$ and $V_i \in \C^{k N_1 \times m}$ for $i \in [N_2]$, we have \[\begin{bmatrix} \circop^{N_1}(V_1) & \cdots & \circop^{N_1}(V_{N_2}) \end{bmatrix} (P^{(k, N_2)} \otimes I_m)^* = \circop^{N_1}\left(\begin{bmatrix} V_1 & \cdots & V_{N_2}\end{bmatrix}\right).\]
\end{lemma}

\begin{proof}[Proof of \cref{lem:block_circ_right}]
  We quote \eqref{eq:M_2} from \cref{lem:interleave} and consider that $P^{(k, N_2)} \otimes I_m$ is a permutation that changes the blockwise indices of $m \times p$ blocks (or, acting from the right, $p \times m$ blocks) exactly the way that $P^{(k, N_2)}$ changes the indices of a vector.
\end{proof}

\begin{lemma} \label{lem:diag_kron_perm}
  Given $k, n \in \N$ and $V_j \in \C^{m_j \times n_j}$ for $j \in [n]$ and setting $M = \sum_{j = 1}^n m_j, N = \sum_{j = 1}^N n_j,$ and $D = \diag(I_k \kron V_j)_{j = 1}^n \in \C^{k M \times k N}$, we have \[D = \diagcornmat{I_k \kron V_1}{I_k \kron V_n} = P_1 (I_k \otimes \diag(V_j)_{j = 1}^n) P_2^*\] where $P_1 = \Pc(P^{(n, k)}, (m_{j \mod_1 n})_{j = 1}^{kn})$ and $P_2 = \Pc(P^{(n, k)}, (n_{j \mod_1 n})_{j = 1}^{kn})$.
\end{lemma}

\begin{proof}[Proof of \cref{lem:diag_kron_perm}]
  We immediately reduce to the case $m_j = n_j = 1$ for all $j$ by observing that $P_1$ and $P_2$ will act on blockwise indices precisely as $P^{(n, k)}$ acts on individual indices.  Here, we replace $V_j$ with $v_j \in \C$, and note that $\diag(V_j)_{j = 1}^n = \diag(v)$.  Hence, we need only remark that \[\left(\diag(I_k \otimes v_\ell)_{\ell = 1}^n\right)_{((i_1 - 1)k + i_2) ((j_1 - 1)k + j_2)} = \left\{\begin{array}{r@{,\quad}l} v_{i_1} & i_1 = j_1 \ \text{and} \ i_2 = j_2 \\ 0 & \text{otherwise} \end{array}\right.,\] while 
  \begin{align*} (P^{(n, k)} (I_k \otimes \diag(v)) P^{(n, k)*})_{((i_1 - 1)k + i_2) ((j_1 - 1)k + j_2)} 
 & = (I_k \otimes \diag(v))_{((i_2 - 1)n + i_1) ((j_2 - 1)n + j_1)} 
  \\ &
  = \left\{\begin{array}{r@{,\quad}l} v_{i_1} & i_1 = j_1 \ \text{and} \ i_2 = j_2 \\ 0 & \text{otherwise} \end{array}\right..\end{align*}
\end{proof}